\documentclass[12pt,leqno]{amsbook}  
\usepackage{minitoc}
\usepackage{amsmath,amsfonts,amstext,amsthm,amssymb,pxfonts}
\usepackage[colorlinks,citecolor=red,pagebackref,hypertexnames=false]{hyperref}
\usepackage[backrefs]{amsrefs}
  
\numberwithin{equation}{section}
\numberwithin{section}{chapter}

\allowdisplaybreaks
\swapnumbers
 
\theoremstyle{plain} 
\newtheorem{lemma}[equation]{Lemma} 
\newtheorem{proposition}[equation]{Proposition} 
\newtheorem{theorem}[equation]{Theorem} 
 
\newtheorem{conjecture}[equation]{Conjecture}

\newtheorem{smallBallConjecture}[equation]{Small Ball Conjecture}
\newtheorem{smooth}[equation]{Smooth Small Ball Conjecture}
\newtheorem{restrictedSmallBallConjecture}[equation]{Restricted Small Ball Conjecture}
\newtheorem{talagrand}[equation]{Talagrand's Theorem}
\newtheorem{schmidt}[equation]{Schmidt's Theorem}
\newtheorem{halasz}[equation]{Hal{\'a}sz' Theorem}
\newtheorem{hyperbolic}[equation]{Hyperbolic Sup Norm Conjecture}
\newtheorem{ellOne}[equation]{$L^1$ Norm Conjecture}
\newtheorem{product}[equation]{Product Rule in Dimension 2}
\newtheorem{Sharp}[equation]{Sharpness of the Hyperbolic Sup Norm Conjecture}
\newtheorem{khintchine}[equation]{Khintchine Inequalities}
\newtheorem{lpi}[equation]{Littlewood Paley Inequalities}
\newtheorem{cww}[equation]{Chang Wilson Wolff Inequality}
\newtheorem{halton}[equation]{Halton's Theorem}
\newtheorem{weak}[equation]{Weak $ L ^{1} $ Bound for the Square Function}
\newtheorem{decomp}[equation]{Calder\'on Zygmund Decomposition}
\newtheorem{good}[equation]{Good $ \lambda $ Inequality} 
\newtheorem{brownian}[equation]{Small Ball Problem for the Brownian Sheet}

\theoremstyle{definition}
\newtheorem{definition}[equation]{Definition} 

\theoremstyle{remark}
\newtheorem*{remark}{Remark}
\newtheorem*{Acknowledgment}{Acknowledgment}

\def\norm#1.#2.{\lVert#1\rVert_{#2}}
\def\Norm#1.#2.{\bigl\lVert#1\bigr\rVert_{#2}}
\def\NOrm#1.#2.{\Bigl\lVert#1\Bigr\rVert_{#2}}
\def\NORm#1.#2.{\biggl\lVert#1\biggr\rVert_{#2}}
\def\NORM#1.#2.{\Biggl\lVert#1\Biggr\rVert_{#2}}

\def\ip#1,#2,{\langle #1,#2\rangle}
\def\Ip#1,#2,{\bigl\langle#1,#2\bigr\rangle}
\def\IP#1,#2,{\Bigl\langle#1,#2\Bigr\rangle}

\def\mid{\,:\,}

\def\abs#1{\lvert#1\rvert}
\def\Abs#1{\bigl\lvert#1\bigr\rvert}
\def\ABs#1{\Bigl\lvert#1\Bigr\rvert}

\def\XXint#1#2#3{{\setbox0=\hbox{$#1{#2#3}{\int}$}
     \vcenter{\hbox{$#2#3$}}\kern-.5\wd0}}

\def\name #1 #2{\operatorname{#1}(#2)}

\begin{document}

\title[Small Ball Problems] {Small  Ball  and Discrepancy Inequalities  }

\author[M.~T.~Lacey]{Michael T. Lacey}

\address{Michael Lacey\\
School of Mathematics\\
Georgia Institute of Technology\\
Atlanta,  GA 30332 USA}

\email{lacey@math.gatech.edu}

 \dedicatory{Dedicated to the Memory of Walter Philipp, Teacher and Steadfast Friend}

%
\maketitle

\tableofcontents

\chapter*{Preface} 

We discuss an inequality for three dimensional Haar functions motivated by 
questions in a range of areas.  These are 
\begin{itemize}
\item Irregularity of Distributions of points in the unit cube, relative 
to boxes in the standard coordinate basis. 

\item Chung's Law for the Brownian Sheet, or equivalently, sharp estimates 
for the probability that the Brownian Sheet has a small sup norm in the unit cube. 

\item  Lower bounds on the number of $ L ^{\infty }$ balls of small radius needed 
to cover certain compact classes of functions with bounded mixed 
derivative in three dimensions. 
\end{itemize}
Of these three questions, the first admits the easiest description, and 
has the longest history, beginning with van Aardenne-Ehrenfest 
\cites{MR0015143,MR0032717}, with 
significant contributions by a variety of authors over many years.  
See the first chapter of Beck and Chen \cite{MR903025}.  Our methods are 
influenced by many of these contributions; the reader will find references to them 
in the pages below.  Indeed, these notes are our effort to understand the 
famous contribution of J{\'o}zsef Beck \cite{MR1032337} to the irregularities 
of distribution in three dimensions, and its connection with other questions in analysis.  
Along the way, we will simplify and extend his argument, in a manner that 
raises hopes that one could resolve the issue in three dimensions.

The latter two problems listed above
have a more sophisticated description, indeed one that admits 
an abstract formulation.  The relationship between them is rather precise, 
and well known, \cites{MR94j:60078,MR2001c:60059}.  

These topics  are unified by their methods of  proofs.  In its simplest manifestation, 
this is a particular inequality about Haar functions in three dimensions, a 
question which can be viewed as just beyond the reach of Littlewood Paley theory. 
We take this question as our main focus, as doing so will permit us to develop the necessary analytical tools with 
some efficiency.  We establish a  partial result in the direction of the 
main conjecture in the subject, Theorem~\ref{t.bl}.  Afterwords, we discuss the other subjects above.

In the subject of Irregularities of Distribution, the principal new result is 
an extension of the result of Beck already cited, namely Theorem~\ref{t.beck}. 
 The entire subject is also of interest in two dimensions; 
we include this theory in our notes, as it is the foundation from which one must 
generalize.  The two dimensional case is substantially easier, and all important 
elements of that theory have been developed see 
\cites{MR95k:60049,MR96c:41052,MR0319933,MR637361} among other references listed 
in the paper below.

The central methods of this paper are those of Harmonic Analysis: Riesz products; 
Littlewood Paley inequalities;  conditional expectation arguments; and  product theory. 
These notes are written with a focus on these issues. (This is the area of expertise 
of the author.) 
We have written a separate chapter recalling some of these basic issues in a 
separate chapter, see Chapter~\ref{background}.   As our subject touches 
a range of issues, we have also included background material on Irregularities 
of Distributions, Approximation Theory, and Probability Theory.  These are 
offered for the convenience of the reader, with the caveat that the author is 
not an expert in these subjects.

\emph{Notation.}  The language and notation of probability and expectation is 
used throughout.  Thus, 
\begin{equation*}
\mathbb E f = \int _{[0,1] ^{d}} f (x) \; dx
\end{equation*}
and $ \mathbb P (A)= \mathbb E \mathbf 1 _{A}$.  This serves to keep formulas 
simpler.  As well, certain conditional expectation arguments are essential to 
us.  We use the notation 
\begin{equation*}
\mathbb P (B\,|\, A)= \mathbb P (A) ^{-1} \mathbb P (A\cap B)\,. 
\end{equation*}
For a sigma field $ \mathcal F$, 
\begin{equation*}
\mathbb E (f\,|\, \mathcal F)
\end{equation*}
is the conditional expectation of $ f$ given $ \mathcal F$.  In all instances, 
$ \mathcal F$ will be generated by a finite collection of atoms $ \mathcal F
_{\textup{atoms}}$,  in which case 
\begin{equation*}
\mathbb E (f\,|\, \mathcal F)=\sum _{A\in \mathcal F _{\textup{atoms}}} 
\mathbb P (A) ^{-1} \mathbb E (f \mathbf 1 _{A}) \cdot \mathbf 1 _{A} \,. 
\end{equation*}

We suppress many constants which do not affect the arguments in essential ways. 
$ A \lesssim B$ means that there is an absolute constant so that $ A \le K B$. 
Thus $ A \lesssim 1$ means that $ A$ is bounded by an absolute constant. 
And $ A \simeq B$ means $ A \lesssim B \lesssim A$.

\begin{Acknowledgment}
Walter Philipp, my thesis advisor who passed away unexpectedly in 
the summer of 2006, introduced me to this topic while I was in 
graduate school.  
Vladimir Temlyakov lead me through the theory that had been developed since 
graduate school days. 
I report on joint work with Dmitry Bilyk. We have benefited from 
several conversations with Mihalis Kolountzakis and Vladimir Temlyakov
on this subject.  A substantial part of 
this manuscript was written while in residence at the University of Crete. 
\end{Acknowledgment}

\chapter{The Small Ball Problem}

 \section{The Principal Conjecture}

 In one dimension, the class of dyadic intervals are $\mathcal D {} \coloneqq {}\{ [j2^k,(j+1)2^k)\mid j,k\in \mathbb Z\} $. 
 Each dyadic interval has a left and right half, indicated below, which are also dyadic.  Define the
 Haar functions 
 \begin{equation*}
h_I \coloneqq -\mathbf 1 _{I_{\textup{left}}}+ \mathbf 1 _{I_{\textup{right}}}
\end{equation*}
 Note that this is an $ L ^{\infty }$ normalization of these functions, which we will 
 keep through out these notes.  This will cause some formulas to look a little odd 
 to readers accustomed to an $ L ^{2}$ normalization for Haar functions.

 In dimension $d $, a \emph{dyadic rectangle} is a product of dyadic intervals, thus an 
element of 
 $\mathcal D^d $.   A Haar function associated to $R $ we take to the be product of the Haar functions associated 
 with each side of $R $, namely  
 \begin{equation*}
 h_{R_1\times\cdots\times R_d }(x_1,\ldots,x_d) {} \coloneqq {}\prod _{j=1}^d h _{R_j}(x_j). 
 \end{equation*}
 This is the usual `tensor' definition.\footnote{Note that we are not claiming that 
 these functions form a basis.}

 We will concentrate on rectangle with a fixed volume, and 
 consider a local problem. 
 This is the `hyperbolic'
 assumption, that pervades the subject. 
 Our concern is the following Theorem and Conjecture concerning a 
 \emph{lower bound} on the $ L ^{\infty }$ norm of sums of hyperbolic Haar functions:

 \begin{talagrand}\label{j.talagrand}
    For dimensions $d\ge2 $, we have 
 \begin{equation}  \label{e.talagrand}
2 ^{-n} \sum _{\abs{ R}= 2 ^{-n} } \abs{ \alpha(R) } {}\lesssim{} n ^{\frac12(d-2) } 
\NOrm \sum _{\abs R \ge 2 ^{-n}} \alpha (R) h_R .\infty . 
 \end{equation}
 Here, the sum on the right is taken over all rectangles with area \emph{at least } $ 2
 ^{-n}$. 
 \end{talagrand}

 \begin{smallBallConjecture} \label{small} For dimension $ d\ge 3$ we have the inequality 
 \begin{equation}\label{e.Talagrand}
2 ^{-n} \sum _{\abs{ R}= 2 ^{-n} } \abs{ \alpha(R) } 
{}\lesssim{} n ^{\frac12(d-2) } \NOrm \sum _{\abs{ R}\ge 2 ^{-n} } \alpha (R) h_R .\infty . 
\end{equation}
 \end{smallBallConjecture}

 This conjecture is, by one square root, better than the trivial estimate available from Cauchy Schwartz,
 see \S~\ref{s.trivial}.  As well, see that section for 
 an explaination as to why the conjecture is sharp. 
 The motivations for the conjecture are indirect, a subject we return to 
 in the discussion of functions with $L^2 $ mixed partials below, \S~\ref{s.realapplications}.
 Nevertheless, we have begun with this conjecture as it provides the 
 quickest path to the essential technical aspects behind the various conjectures 
 of these notes.

 The result in the case of $d=2 $ is that of Talagrand \cite{MR95k:60049}. 
 We will  give the easier proof  
of Temlyakov \cite{MR96c:41052}, which proof resonates with the ideas of  Roth 
\cite{MR0066435},
Schmidt \cite{MR0319933}, and Hal{\'a}sz \cite{MR637361}.  Compare \S~\ref{s.talagrand}
 and \S~\ref{s.schmidt}.

For many applications of interest, one can restrict attention to this version of the 
conjecture 

\begin{restrictedSmallBallConjecture}
\label{c.restricted}
We have the inequality (\ref{e.talagrand}), or (\ref{e.Talagrand}), 
in the case where the coefficients $ \alpha (R)\in \{0,\pm 1\}$, for $\abs{ R}= 2 ^{-n} $
and about 
half of the $ \alpha (R)\neq 0$. Namely, under these assumptions on the 
coefficients $ \alpha (R)$ we have the inequality 
\begin{equation}\label{e.restricted}
\NOrm \sum _{\abs{ R}= 2 ^{-n}} \alpha (R) h_R. \infty . 
\gtrsim n ^{d/2}\,.
\end{equation}
\end{restrictedSmallBallConjecture}

It is possible that the proof would  simplify considerably---and be of interest---if 
one in addition assumed that $ \abs{ \alpha (R)}\equiv1$.  But some of the applications 
may not be available in this case.

The principal point of these notes is to expound on the three dimensional case, 
providing a partial resolution of this case.  We extend and simplify 
 an approach of J.~Beck \cite{MR1032337}, establishing this result.

\begin{theorem}\label{t.bl} In dimension $ d=3$,  
there is a small positive $ \epsilon >0$ for which we have the estimate 
\begin{equation}  \label{e.bkl}
2 ^{-n} \sum _{\abs{ R}=2 ^{-n} } \abs{ \alpha_R } {}\lesssim{} 
n ^{1 - \epsilon  }   \NOrm \sum _{\abs{ R}=2 ^{-n} } \alpha_r h_R .\infty . 
 \end{equation}
\end{theorem}
 
 This result is due to Bilyk
 and Lacey \cite{bl}.  Beck  \cite{MR1032337} established this inequality 
 with  $ n ^{-\epsilon }$ replaced by  
a term logarithmic in $ n$.\footnote{J.~Beck 
 did not state the result this way, as the principal concern of that paper is on 
 the question of irregularities of distribution.  See \S~\ref{s.discrep}.}

 The organization of the proof, at the highest level, and outlined in 
 \S~\ref{s.Short}, is that of J{\'o}zsef Beck \cite{MR1032337}.
 At the same time, both the exact construction and subsequent 
 details are in many respects easier than in Beck's paper. 
 In particular, the 
 construction in that section is a Riesz product construction, following the lines 
 of \S~\ref{s.talagrand}.  But, the product, with our current understanding, 
 must be taken to be `short,' a dictation to us from the 
 third dimension:  The `product rule'  \ref{p.product} does not hold 
 in dimension three.  This unfortunate, and 
critical fact,  forces the 
 definition of `strongly distinct' on us. See  Definition~\ref{d.distinct}.

 Critically, J{\'o}zsef Beck observed that in the case of that the `strongly distinct' 
  does \emph{not} hold, there is a \emph{gain over naive estimates}. 
  See Lemma~\ref{l.SimpleCoincie} and Theorem~\ref{t.NSD}. 
  We will refer to any instance of this phenomena as the \emph{Beck Gain}. 
  The 
  simplest instance of this is discussed in detail in \S~\ref{s.beck}.  
  Here, we obtain a better range of results, and a larger gain, than Beck. 
 
  Beck's insight is that this gain permits one to carry out 
  a proof, provided the Riesz product is sufficiently short, so short that the 
  combinatorial explosion generated by the expansion of the Riesz product does not 
  overwhelm the gain.

  Beck's gain has other  surprising implications, namely in \S~\ref{s.norm}
  we see that  hyperbolic sums of Haar functions obey a 
 range of sub-gaussian estimates,\footnote{This observation is not essential to 
 our main theorem.}
 not predicted by the general theory in \S~\ref{s.exp}. 
 This section employs a conditional expectation argument to permit an effective 
 application of the Beck gain.

 Concerning the 
 value of $ \epsilon $ for which our Theorem holds, it is computable, 
 but we do not carry out this step, as the particular $ \epsilon $ 
we would obtain is certainly not  optimal.  Instead, the point of this proof is 
that the methods pioneered by J{\'o}zsef Beck are more powerful than originally suspected. 
We expect more efficient organizations of the proof will yield quantifiable 
and substantive improvements to the results of this paper.

 \section{The Trivial Bounds}  \label {s.trivial}
 
  The inequality (\ref{e.talagrand}) with an extra square root of $ n$ is easy to prove. 
  
 \begin{lemma}\label{l.hyperbolic}
  It is the case that 
 \begin{equation*}
 \sum _{\abs R=2 ^{-n}} \abs{\alpha(R)}\cdot \abs R \lesssim 
 n ^{\tfrac 12 (d-1)} \NOrm \sum _{\abs R=2 ^{-n}  } \alpha(R) h_R .\infty. \,.
 \end{equation*}
 \end{lemma}
 
 \begin{proof}  Each point 
 $x\in [0,1]^d $, is in at most $  n ^{d-1} $ possible rectangles.  This is the essential 
point dictated by the hyperbolic nature of the problem.  
 Using this, and the Cauchy--Schwartz inequality, we have 
 \begin{align*}
 \sum _{\abs R=2 ^{-n}}
 \abs{ \alpha_R }\cdot \abs R&{}={} 
 	\NOrm \sum _{\abs R=2 ^{-n} }
 \abs{ \alpha_R } {\mathbf 1}_{R} .1. 
 \\&{}\lesssim 
 n^{\tfrac 12 (d-1)} \NORm\Biggl[ \sum _{\abs R=2 ^{-n}}
 \abs{ \alpha_R }^2 {\mathbf 1}_{R}\Biggr]^{1/2}  .1.
 \\&
 {}\lesssim  n^{\tfrac 12 (d-1)}\NOrm \sum _{\abs R=2 ^{-n} } \alpha(R) h_R .2. 
 \\
 &{}\lesssim  n^{\tfrac 12 (d-1)}\NOrm \sum _{\abs R=2 ^{-n} } \alpha(R) h_R .\infty. 
 \end{align*}
 \end{proof} 
 
Let us also see that the Small Ball Conjecture is sharp. 
Indeed, we take the $ \alpha (R)$ to be random choices of signs.  It is immediate 
that 
\begin{equation*}
  2 ^{-n}\sum _{\abs R=2 ^{-n} }
 \abs{ \alpha_R } \simeq n ^{d-1}\,. 
\end{equation*}
We now turn to properties of Rademachers outlined in Chapter~\ref{background}. 
On the other hand, for fixed $ x\in [0,1] ^{d}$ we have 
\begin{equation*}
\mathbb E \ABs{\sum _{\abs R=2 ^{-n} } \alpha(R) h_R (x)} \simeq n ^{\tfrac 12  (d-1)}\,.
\end{equation*}
 It is also well known that sums of Rademacher random variables obey a sub--Gaussian 
 distributional estimate.  The supremum of such sums admit easily estimated upper bounds. 
 In particular, it is enough to test the $ L ^{\infty }$ norm of the sum at a 
grid of $ 2 ^{nd}$ points in the unit cube, hence we have 
\begin{align*}
\mathbb E \NOrm \sum _{\abs R=2 ^{-n} } \alpha(R) h_R .\infty.
& \lesssim 
\sqrt {\log 2 ^{nd}} \cdot 
\sup _{x} \mathbb E \ABs{\sum _{\abs R=2 ^{-n}  } \alpha(R) h_R (x)}  
\\
& \lesssim n ^{ d/2}\,. 
\end{align*}
Comparing these two estimates shows that the Small Ball Conjecture is sharp. 

 The Small Ball Conjecture could be substantially resolved if one could directly 
 show that in the random case that this estimate is sharp.

 \section{Proof of Talagrand's Theorem} \label{s.talagrand}

 We follow the approach of V.~Temlyakov \cite{MR96c:41052} to the stronger inequality 
 (\ref{e.Talagrand}) in the case of $ d=2$, and invite the reader to compare 
 this argument to the proof of Schmidt's Theorem in \S~\ref{s.schmidt}.

  The decisive point in two dimensions is that one has a `product rule.'  Let us formalize 
it as this proposition, and leave the proof to the reader. 

\begin{product}\label{p.product} Let $ R,R'$ be two dyadic rectangles of the 
same area.  Then, 
\begin{equation*}
h_R \cdot h _{R'}\in \bigl\{ 0\,, \, 1 _{R} \,,\, h _{R\cap R'} \bigr\}. 
\end{equation*}
More generally, let $ R_1,R_2,\dotsc, R_k$ be dyadic rectangles of equal area 
and distinct lengths in e.\thinspace g. their first coordinates.  Then 
\begin{equation*}
\prod _{j=1} ^{k} h_ {R_j} \in \bigl\{ 0\,,\, \pm h _{R_1\cap \cdots \cap R_k} \bigr\}\,.
\end{equation*}
\end{product}

 The proof of (\ref{e.Talagrand}) is by duality.  Fix 
 \begin{equation*}
H=\sum _{\abs{ R}\ge 2 ^{-n}} \alpha (R) h_R \,. 
\end{equation*}
 We will construct a function $\Psi $ with $L^1 $ norm at most $1 $, for which 
 the inner product 
 \begin{equation*}
 \ip H,\Psi,= 2 ^{-n-1}\sum _{\abs R=2^{-n} } \abs{ \alpha (R)}\,. 
 \end{equation*}
 This clearly implies the Theorem.  Moreover, the function $\Psi $ is defined as a 
 Riesz product.

Our Riesz product is 
 \begin{align*}
 \Psi &  \coloneqq \prod _{s=1}^n (1+\tfrac12 \psi_s )\,, 
 \\
 \psi _s& =\sum _{R\,:\,\abs{R_1}=2^{-s}, \abs {R_2}=2^{-n+s} }  
\operatorname {sgn} (\alpha (R))  h_R 
 \end{align*}
 Of course $ \Psi $ is non--negative. Moreover, it has $ L^1$ norm one: Expanding the 
 product, the leading term is $ 1$.  All products of $ \psi _s$ are, by
 Proposition~\ref{p.product}, a sum of Haar functions, hence have mean zero.

The Proposition also implies that 
\begin{align*}
\ip H, \Psi , & = \sum _{s=1}^{n} \ip H, \psi _n, 
= 2 ^{-n-1} \sum _{\abs{ R}= 2 ^{-n}} \abs{ \alpha (R)}\,.
\end{align*}
The proof is complete.

\begin{remark}\label{r.notRandom}  If one considers the case of $ \alpha (R) \equiv 1$, 
it is clear that the $ L^\infty $ norm is achieved--or nearly achieved--on a set of measure 
approximately $ 2 ^{-c n}$. That is, the supremum is achieved on a very thin set. 
Experience shows that  Riesz products are very useful in such situations. 
\end{remark}

\begin{remark}\label{r.riesz} Traditionally, a Riesz product is of the form 
\begin{equation*}
\prod _{k=1} ^{\infty } (1+ \cos 4 ^{k}x) \,. 
\end{equation*}
By a well known heuristic, the functions $ \cos 4 ^{k }x $ behave 
as independent random variables, so we don't make a distinction between the 
classical Riesz product and the Riesz products of our proofs. 
Using Riesz products as above has a long history in the subject of irregularities 
of distributions. 
\end{remark}

\section{Exponential Moments} \label{s.exp} 
 
We state a  
distributional estimate for sums of hyperbolic Haars which 
 shapes  the potential forms of approach to the Small Ball Conjecture. 
However, while the estimates we describe here are in general sharp, they admit 
certain improvements,  for small $ p$; see \S~\ref{s.norm}. 

Background on these issues are developed on Chapter~\ref{background}.

\begin{theorem}\label{t.distributional} In dimension $ d\ge2$ we have the 
estimate below, phrased in terms of the exponetial Orlicz Lebesgue spaces.  
\begin{equation}\label{e.distributional}
\NOrm \sum _{\abs{ R}= 2 ^{-n}} \alpha (R) h_R . \operatorname {exp} (L ^{2/ (d-1)}) . 
\lesssim 
\NOrm \Bigl[ \sum _{\abs{ R}= 2 ^{-n}} \alpha (R) ^2 \mathbf 1 _{R} \Bigr] ^{1/2} 
. \infty . \,. 
\end{equation}
\end{theorem}

\begin{remark}\label{r.notInLit} 
The estimates above, specialized to hyperbolic sums in dimension $ 3$ or higher, 
are better than those that appear in the literature associated to the Discrepancy function. 
\end{remark}

Here we are using a typical definition of the exponential integrability classes, 
as given in \S~\ref{s.orlicz}.  
This definition could be for instance 
\begin{equation}\label{e.expL}
\norm X . \operatorname {exp}(L ^{\alpha }).
 \simeq \sup _{p\ge1} p ^{-1/\alpha } \norm X.p. 
\end{equation}
The equivalence holding on any probability space.

Of principal relevance to us is the three dimensional case, where the estimate 
above asserts that the hyperbolic sums are exponentially integrable.

\begin{proof}
The tool is the vector valued Littlewood Paley inequality, with sharp 
rate of growth in the constants as $p \to \infty  $.  As such the proof is 
a standard one, see \cites{MR850744,MR1439553}.

Applying the one dimensional Littlewood Paley inequality in the coordinate 
$ x _1$ we see that 
\begin{equation*}
\NOrm \sum _{\abs{ R}= 2 ^{-n}} \alpha (R) h_R . p . 
\lesssim 
\sqrt p 
\NOrm  
\Bigl[  \sum _{r_1=1} ^{n} 
\Abs{ \sum _{\substack{\abs{ R}= 2 ^{-n}\\ \abs{ R_1}=  2 ^{-r_1} }} \alpha (R) h_R  } ^2 
\Bigr] ^{1/2}  .p. 
\end{equation*}
If we are in dimension $ 2$, note that 
\begin{equation} \label{e.=}
\Abs{ \sum _{\substack{\abs{ R}= 2 ^{-n}\\ \abs{ R_1}=  2 ^{-r_1} }} \alpha (R) h_R  } ^2 
=\sum _{\substack{\abs{ R}= 2 ^{-n}\\ \abs{ R_1}=  2 ^{-r_1} }} \abs{ \alpha (R)} ^2   \mathbf 1
_{R}
\end{equation}
so our proof is complete in this case.

In the higher dimensional case, 
the key point is to observe that the last term can be viewed as an $ \ell ^2 $ 
space valued function.  Then, the Hilbert space analog of the Littlewood Paley 
inequalities applies to the second coordinate, to give us 
\begin{equation*}
\NOrm \sum _{\abs{ R}= 2 ^{-n}} \alpha (R) h_R . p . 
\lesssim 
 p 
\NOrm  
\Bigl[  \sum _{r_1=1} ^{n} \sum _{r_2=1} ^{n} 
\Abs{ \sum _{\substack{\abs{ R}= 2 ^{-n}\\ \abs{ R_j}=  2 ^{-r_j}\,, \ j=1,2 }} 
\alpha (R) h_R  }  ^2 \Bigr] ^{1/2}  .p. 
\end{equation*}
Observe that we have a full power of $ p$, due to the two applications 
of the Littlewood Paley inequalities.   And if $ d=3$, then analog of (\ref{e.=}) 
holds, completing the proof in this case.

In the case of dimension $ d\ge 4$ note that we can continue applying 
the Littlewood Paley inequalities inductively.  They need only be used $ d-1$ 
times due to the hyperbolic assumption. Thus, we have the inequality   
\begin{equation*}
\NOrm \sum _{\abs{ R}= 2 ^{-n}} \alpha (R) h_R . p. 
\lesssim  p ^{ (d-1)/2}
\NOrm \Bigl[ \sum _{\abs{ R}= 2 ^{-n}} \alpha (R) ^2 \mathbf 1 _{R} \Bigr] ^{1/2} 
. p .\,, \qquad 2\le p < \infty  \,. 
\end{equation*}
The implied constant depends upon dimension; the main point we are interested in 
is the rate of of growth of the $ L ^{p}$ norms. 
Assuming that the Square Function of the sum is bounded in $ L ^{\infty }$, the $ L ^{p}$ 
norms can only grow at the rate of $ p ^{(d-1)/2}$, which completes the proof. 
\end{proof}

\begin{remark}\label{r.zygmund} It is a thesis of A.\thinspace Zygmund that 
when one is concerned with product domain questions, the relevant estimates 
are governed by the effective number of parameters involved.  This thesis 
in the hyperbolic setting, says that relevant estimates should be those of 
$ d-1$ parameters in dimension $ d$.  We have just seen one instance of this. 
While it is known that this thesis does not hold in full generality, the hyperbolic 
setting is simple enough that it should hold for most, if not all, questions of interest.
\end{remark}

\section{Definitions and Initial Lemmas for Dimension Three}   \label{s.definition}

The principal difficulty in three and higher dimensions is that the product of Haar functions 
is not necessarily a Haar function. On this point, we have the following proposition 
which does not admit any essential extension. 

\begin{proposition}\label{p.productsofhaars}
Suppose that  $R_1,\ldots,R_k$ are rectangles such that there is no choice of $1\le j<j'\le k$ and no choice of 
coordinate $1\le{} t\le d$ for which we have $R _{j,t}=R _{j',t}$.  Then, for a choice of sign $\varepsilon\in 
\{\pm1\}$
we have 
\begin{equation}  \label{e.productsofhaars}
\prod _{j=1}^k h_R=\varepsilon h_S, \qquad S=\bigcap _{j=1}^k R_k. 
\end{equation}
\end{proposition}

\begin{proof}
Expand the product as 
\begin{equation*}
\prod _{m=1}^\ell  h_{R_m} (x_1,\dotsc, x_d) = \prod _{m=1}^\ell  \varepsilon _m  \prod _{t=1} ^{d}
h_{R_{m,t}} (x_t)
\end{equation*}
Here $ \varepsilon _m\in \{\pm 1\}$.  Our assumption is that for each $ t$, there 
is exactly one choice of $ 1\le m_0\le \ell $ such that 
$ R _{m_0,t}=S_t$.  And moreover, since the minimum value of $ \abs{ R _{m,t}}$ is obtained 
exactly once, for $ m\neq m_0$, we have that 
$ h _{R _{m,t}}$ is constant on $ S_t$.  Thus, in the $ t$ coordinate, the product is 
\begin{equation*}
\varepsilon _{m_0}  h _{S_t} (x_t)
\prod _{1\le m\neq m_0\le \ell } \varepsilon _{m} h _{R _{m,t}} (S_t)\,.
\end{equation*}
This proves our Lemma. 
\end{proof}

Let $\vec r\in \mathbb N^d$ be a partition of $n$, thus $\vec r=(r_1,r_2,r_3)$, 
where the  $r_j$ are non negative integers and $\abs{ \vec r} \coloneqq \sum _t r_t=n$.  
Denote all such vectors at $ \mathbb H _n$. (`$ \mathbb H $' for `hyperbolic.') 
For vector $ \vec r $ let $ \mathcal R _{\vec r} $ be all dyadic rectangles 
$ R$ such that for each coordinate $ k$, $ \lvert  R _k\rvert= 2 ^{-r_k} $.

\begin{definition}\label{d.rfunction} 
We call a function $f$ an \emph{$\mathsf r$ function  with parameter $ \vec r$ } if 
\begin{equation}
\label{e.rfunction} f=\sum_{R\in \mathcal R _{\vec r}}
\varepsilon_R\, h_R\,,\qquad \varepsilon_R\in \{\pm1\}\,. 
\end{equation}
We will use $f _{\vec r} $ to denote a generic $\mathsf r$ function.   
A fact used without further comment is that $ f _{\vec r} ^2 \equiv 1$. 

\end{definition}

\begin{definition}\label{d.distinct}  For vectors $ \vec r_j \in \mathbb N ^{3}$, 
say that $ \vec r_1,\dotsc,\vec r_J$ are \emph{strongly distinct } 
iff for coordinates $ 1\le t\le 3$ the integers 
$ \{ r _{j,t}\mid 1\le j \le J\}$ are distinct.  The product of strongly distinct 
$ \mathsf r$ functions is also an $ \mathsf r$ function.
\end{definition}

The $\mathsf r$ functions we are interested in are:
\begin{equation}\label{e.fr}
f _{\vec r} \coloneqq \sum _{R\in \mathcal R _{\vec r}} \operatorname {sgn}(a (R)) \, h_R 
\end{equation}

\section{J{\'o}zsef Beck's  Short Riesz Product}\label{s.Short}

Let us define relevant parameters by 
\begin{gather} \label{e.q}
q= a  n ^{\varepsilon}  \,, \qquad   b =\tfrac 14
\\
\label{e.rho} 
\widetilde \rho=a  q ^{b}  n^{-1}\,, \qquad \rho =  {\sqrt q} n ^{-1}.
\end{gather}
Here, $ a $ are small positive constants, we 
use the notation of $ b=1/4$ throughout, so as not to obscure  those 
aspects of the argument that that dictate this choice of $ b$. 
 $ \widetilde \rho $ is  a `false' $ L^2$
 normalization for the sums we consider, while the larger term $ \rho $ is the 
 `true' $ L ^{2}$ normalization.
Our `gain over the trivial estimate' in the Small Ball Conjecture is $ q ^{b}
=n ^{\varepsilon /4}$.  $ 0<\varepsilon <1$ is a small constant.  It certainly 
can't be more than $ 1/6$ in view of (\ref{e.partialX2}) though 
there are other more severe restrictions on the size of $ \varepsilon $; the exact determination of 
what we could take $ \varepsilon $ equal to in this proof doesn't seem to be worth 
calculating.

In Beck's paper, the value of $ q=q _{\textup{Beck}} 
= \tfrac {\log n} {\log \log n}$ was much smaller than our value of $ q$.  
The point of this choice is that $ q
_{\textup{Beck}}^{q _{\textup{Beck}}} \simeq n$, with the term $ q ^{q}$ 
controlling many of the combinatorial issues concerning the 
expansion of the Riesz product.\footnote{Specifically, $ q ^{Cq}$ 
is a naive bound for the number of admissible graphs, as defined in \S~\ref{s.nsd}.}
With our 
\emph{substantially  larger value of $ q$},
we need to introduce additional tools to control the 
combinatorics.  These tools are 
\begin{itemize}
\item A Riesz product that will permit us to implement various 
conditional expectation arguments.  
\item Attention to $ L ^{p} $ estimates of various sums, and their growth rates in $ p$.
\item  Systematic use of the Littlewood Paley inequalities, with the sharp exponents in $ p$.
\end{itemize}

Divide the integers $ \{1,2,\dotsc,n\}$ into $ q$ disjoint intervals 
$ I_1,\dotsc, I_q$, and let $ \mathbb A _t \coloneqq \{\vec r\in \mathbb H _n 
\mid r_1\in I_t\}$.  Let 
\begin{equation} 
\label{e.G_t}
F_t =  \sum _{\vec r\in \mathbb A _t}  f _{\vec r}\,. 
\end{equation}
The Riesz product is now a `short product.' 
\begin{equation*}
\Psi \coloneqq \prod _{t=1} ^{q} (1+  \widetilde  \rho F_t) \,.
\end{equation*}
Note the subtle way that the false $ L^2$ normalization enters into the product. 
It means that the product is, with high probability, positive.  And of course, for 
a positive function $ F$, we have $ \mathbb E F=\norm F.1.$, with expectations being 
typically easier to estimate. 
This heuristic is made precise below. 

We need to decompose the product $ \Psi $ into 
\begin{equation}
\label{e.FG}  \Psi =1+\Psi ^{\textup{sd}}+\Psi ^{\neg  }\,,
\end{equation}
where the two pieces are the `strongly distinct' and `not strongly distinct' pieces.  
To be specific, for integers $ 1\le u\le q$, 
let 
\begin{equation*}
\Psi ^{\textup{sd}}_k
 \coloneqq  \widetilde \rho ^{k}
 \sum _{1\le v_1< \cdots < v_k\le q} \; 
 \sideset{}{^{\textup{sd}}}\sum _{\vec r_t\in \mathbb A _{v_t}}  \prod _{t=1} ^{u} f _{\vec r _t}
\end{equation*}
 where $\sideset{}{^{\textup{sd}}}\sum $ is taken to 
 be over all $ k$ tuples of vectors  $\{(\vec r_1
 ,\dotsc, \vec r_k)\in \prod _{t=1} ^{k} \mathbb A _{v_t} \} $  such that: 
 \begin{equation}  \label{e.distinct}
 \text{ the vectors  $\{\vec r _{t }\mid 1\le{} t\le k\} $ are strongly 
 distinct. }
 \end{equation}
Then define 
\begin{equation} \label{e.zCsd}
\Psi ^{\textup{sd}} {} \coloneqq {}
\sum _{u=1} ^{q}  \Psi ^{\textup{sd}} _u
\end{equation}

With this definition, it is clear that we have 
\begin{equation} \label{e.gain>trivial}
\ip H_n, \Psi ^{\textup{sd}}, 
=\ip H_n, \Psi ^{\textup{sd}} _{1}, 
\gtrsim   q ^{b} \cdot n ^{ 1}\,,
\end{equation}
so that $ q ^{b}$ is our `gain over the trivial estimate.'

The bulk of the proof is taken up with the proof of the technical estimates 
below.  The main point of the Lemma is the last estimate, (\ref{e.sd1}), which 
with (\ref{e.gain>trivial}) above proves Theorem~\ref{t.bl}.

\begin{lemma}\label{l.technical} We have these estimates: 
\begin{align}
\label{e.<0}
\mathbb P ( \Psi <0) &\lesssim \operatorname {exp} (-A q ^{1-2b} )\,;
\\  \label{e.expq2b}
\norm \Psi . 2. & \lesssim \operatorname {exp} (a' q ^{2b})\,;
\\
\label{e.E1}
 \mathbb E  \Psi  & = 1 \,;
\\ 
\label{e.CL1}
\norm  \Psi .1. &\lesssim 1\,;
\\ \label{e.neq1}
\norm \Psi ^{\neg } .1. &\lesssim 1\,;
\\
\label{e.sd1}
\norm \Psi ^{\textup{sd}}.1. &\lesssim 1\,.
\end{align}
Here, $ 0<a'<1$, in (\ref{e.expq2b}), is a small constant, decreasing to zero as $ a $ in (\ref{e.q}) 
goes to zero; and $ A>1$, in (\ref{e.<0}) is a large constant, tending to infinity as $ a$ in (\ref{e.q}) 
goes to zero. 
\end{lemma}

\begin{proof}
We give the proof of the Lemma, assuming our main inequalities proved in the 
subsequent sections.  In particular, the first two estimates of our Lemma are 
substantial, as they reflect the influence of the 
non trivial sub--gaussian estimates of \S~\ref{s.norm}. 

\smallskip 

\emph{Proof of (\ref{e.<0}).} The main tool is the  distributional 
estimate (\ref{e.partialX2}). Observe that 
\begin{align*}
\mathbb P (\Psi <0)&\le \sum _{t=1} ^{q} \mathbb P ( \widetilde \rho \,  F_t < -1)
\\
& = \sum _{t=1} ^{q}  \mathbb P (\rho F_t < - a ^{-1} q ^{1/2-b})
\\
& \lesssim q \operatorname {exp} (c a ^{-1} q ^{1-2b})\,. 
\end{align*}

\smallskip

\emph{Proof of (\ref{e.expq2b}).} The proof of this is detailed enough that 
we postpone it to Lemma~\ref{l.2b} below.

\smallskip

\emph{Proof of (\ref{e.E1}).} 
Expand the product in the definition of $ \Psi $.  
The leading term is one.  Every other term is a product 
\begin{equation*}
\prod _{k\in V} \widetilde \rho \, F_k \,,  
\end{equation*}
where $ V$ is a non-empty subset of $  \{1 ,\dotsc, q \}$.  This product is in turn 
a product of $ \mathsf r$ functions.  Among this product, the maximum in the first coordinate 
is unique. This fact tells us that the expectation of this product of $ \mathsf r$ functions
is zero.  So the expectation of the product above is zero. The proof is complete.

\smallskip

\emph{Proof of (\ref{e.CL1}).} We use the first two estimates of our Lemma. 
Observe that 
\begin{align*}
\norm \Psi .1. & = \mathbb E \Psi -2 \mathbb E \Psi \mathbf 1 _{\Psi <0}
\\
& \le 1+ 2\mathbb P (\Psi <0) ^{1/2} \norm \Psi .2.
\\
& \lesssim 1+ \operatorname {exp} ( -A q ^{1-2b}/2+ a' q ^{ 2b})\,.
\end{align*}
We have taken $ b=1/4$ so that $ 1-2b=2b$.  For sufficiently small $ a$ in  (\ref{e.q}),
we will have $ A \gtrsim a'$. 
We see that (\ref{e.CL1}) holds.\footnote{Here of course we are strongly using the 
fact $ \Psi $ is positive with high probability.}

Indeed,   Lemma~\ref{l.2b} proves a uniform estimate, namely 
\begin{equation*}
\sup _{V\subset \{1 ,\dotsc, q\}}
\mathbb E \prod _{v \in V}  (1+ \widetilde \rho F_t) ^2 \lesssim 
\operatorname {exp} (a' q ^{2b})\,. 
\end{equation*}
Hence, the argument above proves 
\begin{equation}\label{e.CCLL}
\sup _{V\subset \{1 ,\dotsc, q\}} 
\NOrm \prod _{v \in V}  (1+ \widetilde \rho F_t) .1. \lesssim 1\,. 
\end{equation}

\smallskip 
 
\emph{Proof of (\ref{e.neq1}).}  The primary facts are (\ref{e.CCLL}) 
and Theorem~\ref{t.NSD}; we use the notation devised for that Theorem.

Note that the Inclusion Exclusion principle gives us the identity 
\begin{equation*}
\Psi ^{\neg} = \sum _{\substack{V\subset \{1 ,\dotsc, q\}\\  \abs{ V}\ge 2}} 
(-1) ^{\abs{ V}+1}
\operatorname {Prod} (\operatorname {NSD} (V)) 
\cdot \prod _{t\in \{1 ,\dotsc, q\}-V} (1+ \widetilde \rho F_t)\,. 
\end{equation*}
We use the triangle inequality, the estimates of Lemma~\ref{l.2b}, 
H\"older's inequality, with indices $ 1+1/q ^{2b}$ and $ q ^{2b}$, and 
the estimate of (\ref{e.NSD}) in the calculation below.
Notice that we have 
\begin{align*}
\sup _{V\subset \{1 ,\dotsc, q\}} 
\NOrm \prod _{v \in V}  (1+ \widetilde \rho F_t) .1+q ^{2b}. 
&\le 
\sup _{V\subset \{1 ,\dotsc, q\}} 
\NOrm \prod _{v \in V}  (1+ \widetilde \rho F_t) .1. ^{(1+q ^{-2b})/(1-q ^{-2b})}
\\& \qquad \times  
\NOrm \prod _{v \in V}  (1+ \widetilde \rho F_t) .2. ^{q ^{-2b}/(1+q ^{-2b})}
\\
& \lesssim 1\,. 
\end{align*}
And recall that $ q ^{2b}= q ^{1/2}$ is a small power of $ n$.  So the $ L ^{p}$ 
norms that we need on terms arising from $ \operatorname {NSD} (V)$ below
are for moderate values of $ p$, namely we only need $ p \le q ^{2b} $.
This is a key reason why we can 
control the combinatorial explosion associated with our short Riesz product. 

We estimate 
\begin{align*}
\norm \Psi ^{\neg}.1. 
&\le 
\sum _{\substack{V\subset \{1 ,\dotsc, q\}\\  \abs{ V}\ge 2}} \NOrm 
\operatorname {Prod} (\operatorname {NSD} (V)) 
\cdot \prod _{t\in \{1 ,\dotsc, q\}-V} (1+ \widetilde \rho F_t) .1. 
\\
& \le 
\sum _{\substack{V\subset \{1 ,\dotsc, q\}\\  \abs{ V}\ge 2}} \norm 
\operatorname {Prod} (\operatorname {NSD} (V)) . q ^{2b}. \cdot 
\NOrm  \prod _{t\in \{1 ,\dotsc, q\}-V} (1+ \widetilde \rho F_t) .1+ q ^{-2b}. 
\\
& \lesssim \sum _{v=2} ^{q} 
[ q ^{C'}  n ^{-1/6} ] ^{v}
\\
& \lesssim q ^{C''} n ^{-1/6}
\\
& \lesssim n ^{-\varepsilon '} \lesssim 1\,.
\end{align*}

\smallskip 
\emph{Proof of (\ref{e.sd1}).}  This follows from (\ref{e.neq1}) and (\ref{e.CL1}) 
and the identity $ \Psi =1+\Psi ^{\textup{sd}}+\Psi ^{\neg}$ and the triangle 
inequality. 

\end{proof}

\section{The Beck Gain in the Simplest Instance} \label{s.beck}

Beck   considered sums of products of 
$\mathsf r$ functions that are \emph{not} strongly distinct,
and observed that the 
$ L^2$ norm of the same are smaller than one would naively expect.  
This is what we call the \emph{Beck Gain.} 
A product of $ \mathsf r$ functions 
will not be strongly distinct if the product involves two or more  
vectors which agree in  one or more  coordinates.    
In this section, we study the sums of products of two $ \mathsf r$ functions 
which are not strongly distinct. 
A later section, \S~\ref{s.nsd}, will study the general case. 

In this section, and again in \S~\ref{s.nsd}, we will use this notation. 
For a subset $ \mathbb C \subset \mathbb H _n ^{k}$, let 
\begin{equation}\label{e.prod}
\operatorname {Prod} (\mathbb C ) 
\coloneqq 
\sum _{ (\vec r_1 ,\dotsc, \vec r_k) \in \mathbb C } 
\prod _{j=1} ^{k} f _{\vec r_j}\,. 
\end{equation}
In this section, we are exclusively interested in $ k=2$. 

Let 
$
\mathbb C (2) \subset  \mathbb H _n ^{2}
$
consist of all pairs of distinct $ \mathsf  r$  vectors 
$ \{\vec r_1,   \vec r _{2}\} $ for which 
$
 r _{1,2}= r _{2,2}\,
$
J.\thinspace Beck calls such terms `coincidences' and we will continue to use that term. 
We need norm estimates on 
the sums of  products of such $ \mathsf r$ vectors.  

\begin{lemma}\label{l.SimpleCoincie}[\textup{\textbf{The  Simplest Instance of the Beck Gain.}}] 
We have these 
estimates for arbitrary subsets $ \mathbb C \subset \mathbb C (2)$
\begin{align}\label{e.Simple2}
\norm \operatorname {Prod} (\mathbb C ) .p.
& \lesssim   p ^{5/4}  n ^{7/4}\,.
\end{align}
Moreover, if we have $ \mathbb C = \mathbb C (2) \cap \mathbb A _s \times \mathbb A _t$
for some $ 0\le s,t \le q$ we have 
\begin{align}\label{e.Simple2q}
\norm \operatorname {Prod} (\mathbb C ) .p.
& \lesssim   p ^{3/2}  n ^{3/2} q ^{-1/2} \,.
\end{align}
Finally, we have the estimate 
\begin{equation}\label{e.Xsimple}
\NOrm \sum _{\substack{\vec r \neq \vec s\in \mathbb A _s\\ r_1=s_1 }} 
f _{\vec r} \cdot f _{\vec s} .p. 
\lesssim p ^{3/2} n ^{3/2} q ^{-1}\,. 
\end{equation}
\end{lemma}

We will use the second estimate of the Lemma, which we do \emph{not} claim for 
{arbitrary} subsets of $ \mathbb C (2)$. 
This estimate appears to be sharp, in that the collection $ \mathbb C (2)$ 
has three free parameters, and the estimates is in terms of $ n ^{3/2}$. 
Note that for $ p \simeq n$ we have 
\begin{equation*}
\norm \operatorname {Prod} (\mathbb C_2 ) .n. 
\simeq 
\norm \operatorname {Prod} (\mathbb C_2 ) .\infty . \,. 
\end{equation*}
And the latter term can be as big as $ n ^{3}$, which matches the bound above.

The proof of the Lemma requires we pass through an intermediary collection 
of four tuples of $ \mathsf r$ vectors. 
Let $ \mathbb B (4) \subset \mathbb H _n ^{4}$ be four tuples of distinct vectors 
$ (\vec r,\vec s,\vec t,\vec u)$ for which (i) $ r_1=s_1$ and $ t_1=u_1$; 
and (ii) in the second and third coordinate  two of the vectors agree.

\begin{proof}  
The method of proof is probably best explained by considering first the case 
of $ p=2$.  Observe that 
\begin{equation*}
\norm \operatorname {Prod} (\mathbb C ).2. ^2 
=\mathbb E \operatorname {Prod} (\mathbb B ) \,,
\end{equation*}
where $ \mathbb B = \mathbb C \times \mathbb C \cap \mathbb B (4) $.
Indeed, the main point is that in order for 
\begin{equation*}
\mathbb E f _{\vec r_1} \cdot  f _{\vec r_2} \cdot  f _{\vec r_3} \cdot  f _{\vec r_4}
\neq 0
\end{equation*}
there is a coincidence among the four vectors in each coordinate.  But this 
is the definition of $ \mathbb B (4)$. 
Thus the case $ p=2$ follows immediately from Lemma~\ref{l.2}.

Now, let us consider $ 4\le p  \le n $, as the inequalities we prove are  trivial for $ p>n$. 
Let $ K _{3/2}$ be the best constant in the inequality 
\begin{equation*}
N (p) \coloneqq \sup \norm \operatorname {Prod} (\mathbb C ).p.
\le K p ^{3/2} n ^{3/2}\,. 
\end{equation*}
Here the supremum is over all choices of $ n$ and $ \mathsf r$
functions.  We give an \emph{a priori} estimate of $ K _{3/2}$. 
We  define $ K _{7/4}$ similarly.

Each pair $ (\vec r, \vec s)\in \mathbb C $ must be distinct in the 
first and third coordinates.  Therefore, we can apply the 
Littlewood Paley inequalities in these coordinates to estimate 
\begin{equation*}
N (p) \coloneqq \norm \operatorname {Prod} (\mathbb C ).p. 
\lesssim 
 p \NOrm \Bigl[ \sum _{a,b} \ABs{ 
\sum _{\substack{(\vec r, \vec s)\in \mathbb C \\ \max \{r_1, s_1\}=a 
\\ \max \{r_3, s_3\}=b}} f _{\vec r} \cdot f _{\vec s} 
} ^2  \Bigr] ^{1/2} .p.\,.
\end{equation*}
Here, we have a full power of $ p$, as we apply the Littlewood Paley inequalities 
twice.  Observe that 
\begin{equation*}
\sum _{a,b} \ABs{ 
\sum _{\substack{(\vec r, \vec s)\in \mathbb C \\ \max \{r_1, s_1\}=a 
\\ \max \{r_3, s_3\}=b}} f _{\vec r} \cdot f _{\vec s} 
} ^2= \sharp \mathbb C + \sum _{i\neq j\in \{1,2,3,4\}}
\operatorname {Prod} (\mathbb C _{i,j})+
 \operatorname {Prod} (\mathbb B _{\textup{max}} )\,.
\end{equation*}
The term $  \sharp \mathbb C$ arises from the diagonal of the square. 
The terms $ \mathbb C _{i,j}$ are  
\begin{equation*}
\begin{split}
\mathbb C _{i,j}\coloneqq \{ (\vec r_1,\vec r_2,\vec r_3,\vec r_4)\in \mathbb C 
\times  \mathbb C 
\mid  \vec r_i=\vec r_j\,, \textup{and the other two vectors are distinct}\}
\end{split}
\end{equation*}
Note that by definition, $\mathbb C _{1,2}=\mathbb C _{3,4}=\emptyset $. 
The term $ \mathbb B _{\textup{max}}$  is 
\begin{equation*}
\begin{split}
\mathbb B _{\textup{max}}\coloneqq \{ (\vec r_1,\vec r_2,\vec r_3,\vec r_4)\in \mathbb C 
\times  \mathbb C 
\mid   \textup{the maximum in the first and third coordinates occur twice}\}
\end{split}
\end{equation*}

Then, we can estimate by the triangle inequality, 
and the sub-additivity of $ x\mapsto \sqrt x$, 
\begin{align} \label{e.Cinduct}
\begin{split}
 \norm \operatorname {Prod} (\mathbb C ).p. 
&\lesssim   p (\sharp \mathbb C)  ^{1/2} + 
 p \sum _{i<j\in \{1,2,3,4\}}
\norm \operatorname {Prod} (\mathbb C _{i,j}) . p/2. ^{1/2} 
+ p \norm \operatorname {Prod} (\mathbb B _{\textup{max}} ). p. 
\\
& \lesssim  p ^{3/2} + p n ^{1/2} N (p/2) ^{1/2}+ 
p  \norm \operatorname {Prod} (\mathbb B _{\textup{max}} ). p. \,. 
\end{split}
\end{align}

Using the estimate (\ref{e.2}), we see that 
\begin{align*}
\norm \operatorname {Prod} (\mathbb C ).p. 
& \lesssim  p n ^{3/2}+ p n ^{1/2} N (p/2) ^{1/2}+ p ^{5/4} n ^{7/4}
\\
& \lesssim p ^{5/4} n ^{7/4}+ p n ^{1/2} N (p/2) ^{1/2} 
\end{align*}
This implies that 
\begin{align*}
K _{7/4} &\lesssim 1+ \sup _{2\le p \le n} 
p ^{-5/4} n ^{7//4} p n ^{1/2} \bigl( K p ^{5/4}n ^{7/4} \bigr) ^{1/2} 
\\
& \lesssim 1 + K _{5/7} ^{1/2} \,. 
\end{align*}
Clearly, this implies $ K _{7/4} \lesssim 1$. 

Using the estimate (\ref{e.2max}), the proof  that $ K _{3/2} \lesssim 1$ 
is entirely similar. 

\end{proof}

Recall that 
$ \mathbb B (4) \subset \mathbb H _n ^{4}$ be four tuples of distinct vectors 
$ (\vec r,\vec s,\vec t,\vec u)$ for which (i) $ r_1=s_1$ and $ t_1=u_1$; 
and (ii) in the second and third coordinate  two of the vectors agree.

\begin{lemma}\label{l.2}    For any subset $ \mathbb B \subset \mathbb B (4)$
\begin{equation}\label{e.2}
\norm \operatorname {Prod} (\mathbb B  ) .p. \lesssim  \sqrt p \, n ^{7/2}\,. 
\end{equation}
Moreover, for $ \mathbb B \subset \mathbb B (4)\cap ( \mathbb A _s \times \mathbb A _t )^{2}$, 
for any choice of $ 0\le s\neq t\le q$, we have 
\begin{equation}\label{e.2q}
\norm \operatorname {Prod} (\mathbb B  ) .p. \le c \sqrt p \, n ^{7/2} q ^{-2}\,. 
\end{equation}
If we do not consider arbitrary subsets, the estimates improve.  We have the 
the estimates 
\begin{align}\label{e.2'}
\norm \operatorname {Prod} (\mathbb B (4)  ) .p. & \lesssim    p\, n ^{3}\,,
\\ \label{e.2a'}
\norm \operatorname {Prod} (\mathbb B  (4) 
\cap ( \mathbb A _s \times \mathbb A _t )^{2}) .p. & \lesssim  p n ^{3}\,. 
\end{align}
Finally, define 
\begin{equation*}
\begin{split}
\mathbb B _{\textup{max}} & \coloneqq 
\{ (\vec r_1, \vec r_2, \vec r_3, \vec r_4) \in
\mathbb B  (4)
\cap ( \mathbb A _s \times \mathbb A _t )^{2} \,:\, 
\\ & \qquad \textup{the maximum in second and third coordinates occur twice}
\}
\end{split}
\end{equation*}
Then, we have the estimate 
\begin{equation}\label{e.2max}
\norm \operatorname {Prod} (\mathbb B _{\textup{max}} ) .p.  \lesssim  p n ^{3}\,. 
\end{equation}
\end{lemma}

This Lemma, with exponents on $ n$ being $ n ^{7/2}$ appears in Beck's paper 
\cite{MR1032337}, in the case of $ p=2$.  The $ L ^{p}$ variants, 
following from consequences of Littlewood Paley inequalities, are important for us. 

The first group of estimates are recorded, as it is interesting that they 
apply to arbitrary subsets of $ \mathbb B (4)$.  We will rely upon the 
second group of  estimates.  Pointed out to us by Mihalis Kolountzakis, these 
estimates are better for all ranges of $ p \le n$.

\begin{proof}
We discuss (\ref{e.2}) explicitly, and note as we go the improvements needed 
to get the estimate (\ref{e.2q}).

The proof is a case analysis, 
depending upon the number of $\{ \vec r, 
\vec s, \vec t, \vec u\} $ at which the maximums occur in the second and third coordinates. 
We proceed immediately to the cases.

Let $ \mathbb B _1 \subset \mathbb B $ consist of those four--tuples 
$ \{ \vec r, \vec s, \vec t, \vec u\}$ for which 
\begin{equation*}
r_2=t_2=\max \{r_2,s_2,t_2,u_2\}\,, \qquad 
r_3=t_3=\max \{r_3,s_3,t_3,u_3\}\,.
\end{equation*}
This collection is  empty,   for necessarily we must have $ r_1=s_1=t_1=u_1$, 
but then   $ \vec r=\vec s$, as the parameters of all vectors is $ n$. 
This violates the definition of $ \mathbb B $. 

\smallskip 

Let $ \mathbb B _{3} \subset \mathbb \mathbb B  $ consist of those four--tuples 
$ \{ \vec r, \vec s, \vec t, \vec u\}$ for which 
\begin{equation*}
r_2=t_2=\max \{r_2,s_2,t_2,u_2\}\,, \qquad 
r_3=u_3=\max \{r_3,s_3,t_3,u_3\}\,.
\end{equation*}
That is, the maximal values involve \emph{three} distinct vectors.  These 
four vectors can be depicted as 
\begin{equation*}
\vec r=\left( \begin{array}{c}
r_1 \\ r_2 (\Box) \\ r_3 \end{array}\right) 
\,,\quad 
\vec s= 
\left( \begin{array}{c}
r_1\\ \ast \\ \Box   \end{array}\right) 
\,,\quad 
\vec t=
\left( \begin{array}{c}
t_1 \\ r_2 \\ \Box  \end{array}\right) 
\,,\quad  
\vec u=
\left( \begin{array}{c}
t_1\\ \Box  \\ r_3 \end{array}\right) 
\end{equation*}
A $ \Box$ denotes a parameter which is determined by other choices.  
It is essential to note that choices of $ r_1$ and $ r_3$ determine the value 
of $ r_2$ (hence the $ \Box$ in the middle coordinate for $ \vec r$),
and so the vector $ \vec r$.
  The only free parameters 
are (say) $ s_2$, denoted by an $ \ast$ above.

But, note that we must then have $ \abs{ \vec s}=s_1+s_2+s_3<n$.  Therefore this 
case is empty.

\smallskip 

Let $ \mathbb B _{4}$ be those four tuples  four tuples $ \{ \vec r, 
\vec s, \vec t, \vec u\}\in \mathbb B$ such that $ s_2=t_2 $ and $ r_3=u_3$. 
That is there are \emph{four} vectors involved in the maximums of the second 
and third coordinates.  These four vectors can be represented as 
\begin{equation} \label{e.4vec}
\vec r=\left( \begin{array}{c}
r_1 \\ \Box  \\ r_3\end{array}\right) 
\,,\quad 
\vec s= 
\left( \begin{array}{c}
r_1\\ s_2 \\ \Box   \end{array}\right) 
\,,\quad 
\vec t=
\left( \begin{array}{c}
t_1 \\ s_2 \\ \Box  \end{array}\right) 
\,,\quad  
\vec u=
\left( \begin{array}{c}
t_1\\ \Box  \\ r_3 \end{array}\right) 
\end{equation}

The next argument proves (\ref{e.2}). 
Let $ \mathbb B _{4} (a,a',b)$ be those four tuples 
$ \{ \vec r, 
\vec s, \vec t, \vec u\}\in \mathbb B$ such that 
\begin{equation*}
  r_1=s_1=a\,, \quad 
  t_1=u_1=a'\, , \quad s_2,t_2=b\,.
\end{equation*}
The point to observe is that 
\begin{equation*}
\norm \operatorname {Prod} (\mathbb B _{4} (a,a',b)).p. \le C \sqrt p\,  \sqrt n\,. 
\end{equation*}
As there at most $ \lesssim n ^3  $ choices for $ a,a'$ this proves the Lemma.  
(And, in the case of (\ref{e.2q}), there are at most $ n(n/q) ^2 $ choices for 
these three parameters.)

Indeed, we have not specified $ r_3=u_3$.  Since all vectors are distinct, 
$ a\neq a'$, and in considering the norm above, we ignore $ \vec s$ and $ \vec t$, 
as they are completely specified by the datum $ (a,a',b)$.  The product 
$ f _{\vec r} \cdot f _{\vec u}$, in the second coordinate, is equal in distribution to 
a Rademacher function.  And then the estimate above follows.  
The proof of (\ref{e.2}) and (\ref{e.2q}) are finished. 

\medskip 

We turn to the proof of (\ref{e.2'}) and (\ref{e.2a'}), arguing similarly. 
Let $ \mathbb B _{4} (a,a')$ be those four tuples 
$ \{ \vec r, 
\vec s, \vec t, \vec u\}\in \mathbb B$ such that 
\begin{equation*}
  r_1=s_1=a\,, \quad 
  t_1=u_1=a'\,.
\end{equation*}
The point to observe is that 
\begin{equation*}
\norm \operatorname {Prod} (\mathbb B _{4} (a,a')).p. \le C  p\,  n\,. 
\end{equation*}
As there at most $ \lesssim n  ^2  $ choices for $ a,a'$ this proves the Lemma.  
(And, in the case of (\ref{e.2a'}), there are at most $ (n/q) ^2 $ choices for 
these two parameters.)

The point is that $ \operatorname {Prod} (\mathbb B _{4} (a,a'))$ 
splits into a product.  Namely, 
\begin{align*}
\operatorname {Prod} (\mathbb B _{4} (a,a'))
&= \operatorname {Prod} (\mathbb B _{4,1} (a,a')) 
\cdot 
\operatorname {Prod} (\mathbb B _{4,2} (a,a'))
\\
\operatorname {Prod} (\mathbb B _{4,1} (a,a'))
& \coloneqq 
\bigl\{ \{\vec r, \vec u\}
\mid  \textup{there exists
$\{\vec s, \vec t\}\in \mathbb B ^2 $ with $ \{\vec r,\vec s, \vec t,\vec u\}
\in \mathbb B_4 (a,a')$}
\bigr\}\,,
\\
\operatorname {Prod} (\mathbb B _{4,2} (a,a'))
& \coloneqq 
\bigl\{ \{\vec s, \vec t\}
\mid  \textup{there exists
$\{\vec r, \vec r\}\in \mathbb B ^2 $ with $ \{\vec r,\vec s, \vec t,\vec u\}
\in \mathbb B_4 (a,a')$}
\bigr\}\,.
\end{align*}

We estimate 
\begin{equation*}
\norm 
\operatorname {Prod} (\mathbb B _{4} (a,a')) .p. 
\lesssim 
\norm 
\operatorname {Prod} (\mathbb B _{4,1} (a,a')) .2p.
\cdot 
\norm 
\operatorname {Prod} (\mathbb B _{4,2} (a,a')) .2p.\,. 
\end{equation*}
Both of the last two norms are at most $ \lesssim \sqrt p \cdot n ^{1/2}$, 
which will finish the proof.  

That is the estimate is 
\begin{equation} \label{e.mk}
\norm 
\operatorname {Prod} (\mathbb B _{4,1} (a,a')) .2p.
\lesssim \sqrt p \cdot n ^{1/2} \,. 
\end{equation}
We may assume without loss of generality that $ a>a'$.  The pairs in 
$\operatorname {Prod} (\mathbb B _{4,1} (a,a')) $ consist of the 
two vectors $ \vec r$ and $ \vec u$ in (\ref{e.4vec}). 
These two vectors are parameterized by 
$ u_2$, say.  Since $a=r_1<a'=u_1 $, and $ r_3=u_3$, the hyperbolic assumption implies $ u_2$ 
is the maximal coordinate.  Therefore, the Littlewood Paley inequality applies.  

The proof of (\ref{e.2a'}) is exactly the same, just noting that 
$ a,a'$ can only take $ (n/q) ^2 $ values in that case.

\medskip 

We turn to the proof of the estimate (\ref{e.2max}).  Here, it suffices 
to prove that 
\begin{equation}\label{e.2MAX}
\norm \operatorname {Prod} ( \mathbb B (4)\cap (\mathbb A _s \times \mathbb A _t) ^2 
- \mathbb B _{\textup{max}}).p. \lesssim p  n ^{3}\,. 
\end{equation}

This last collection of four tuples of vectors can be further subdivided into 
finite number of  collections, $ \mathbb B_j'$, for $ 1\le j \le 6$ .
Take $ \mathbb B _1'$ to be a subset of 
 four tuples $ (\vec r, \vec s,\vec t,\vec u) 
\in \mathbb B (4)\cap (\mathbb A _s \times \mathbb A _t) ^2 
- \mathbb B _{\textup{max}} $ with 
\begin{equation*}
\vec r=\left( \begin{array}{c}
r_1 \\ r_2  \\ r_3\end{array}\right) 
\,,\quad 
\vec s= 
\left( \begin{array}{c}
s_1\\ r_2 \\  s_3  \end{array}\right) 
\,,\quad 
\vec t=
\left( \begin{array}{c}
t_1 \\ t_2 \\ s_3  \end{array}\right) 
\,,\quad  
\vec u=
\left( \begin{array}{c}
s_1\\ t_2 \\ t_3 \end{array}\right) 
\end{equation*}
Here we assume that $ r_1$ is the unique maximal integer among $ \{r_1,s_1,t_1,u_1\}$. 
Note that $ \vec r$ and $ \vec s$ have a coincidence in the second coordinate; 
$ \vec s,\vec t$ have a coincidence in the first coordinate; and $ \vec s,\vec u$ 
have a coincidence in the third coordinate.  
The other collections $ \mathbb B '_j$ differ in the location of the maximums 
in either the first and third coordinates, and the particular patterns of 
coincidences. 

It is important to observe that we necessarily have $ r_1>t_1>s_2=u_2$. 
And we will apply the Littlewood Paley inequality in the $ r_1$ and $ t_1$ 
variables.  Clearly, we can apply the Littlewood Paley inequality 
in $ r_1$ to get the estimate 
\begin{equation*}
\Norm \operatorname {Prod} (\mathbb B _1') .p. 
\lesssim 
\sqrt p \NOrm \Bigl[ \sum _{a}  \operatorname {Prod} (\mathbb B _1' (a)) ^2 
\Bigr] ^{1/2} .p. 
\end{equation*}
Here, $ \mathbb B _1 (a)$ is the collection of all four tuples $ \{\vec r,
\vec s,\vec t,\vec u\}\in \mathbb B _1'$ with $ r_1=a$. 

Next, we use the triangle inequality in the values of $ r_2$ and $ s_3$. 
Note that with $ r_1,r_2,s_3$ specified, the values of $ r_3$ and $ s_1$ 
are then forced.  Let $ \mathbb B _1' (a,b,c)$ be the 
pairs of vectors $ \{\vec t,\vec u\}$ for which there are vectors $ \{\vec r, 
\vec s\}$ with $ \{\vec r,\vec s,\vec t,\vec u\}\in \mathbb B _1'$, with in addition 
\begin{equation*}
r_1=a\,, \quad r_2=b\,, \quad s_3=b\,. 
\end{equation*}
By the triangle inequality, we can estimate 
\begin{equation*}
\Norm \operatorname {Prod} (\mathbb B _1') .p. 
\lesssim 
\sqrt p \cdot n ^{5/2} \sup _{a,b,c} 
\Norm \operatorname {Prod} (\mathbb B _1' (a,b,c)) .p. \,. 
\end{equation*}
Now, among the pairs of vectors in $ \mathbb B _1' (a,b,c)$ 
have only one free parameter, which can be taken to be the maximum 
in the first coordinate.  Thus, by the Littlewood Paley inequality we see that 
\begin{equation*}
\Norm \operatorname {Prod} (\mathbb B _1') .p. 
\lesssim 
p n ^{3}\,. 
\end{equation*}
The analysis of the other possible forms of the collections $ \mathbb B '_j$
proceeds along similar lines. We omit the details.

\end{proof}

There is another corollary to the proof above required at a later stage of the proof. 
For an integer $ a$, let 
$ \mathbb B_a (4) \subset \mathbb H _n ^{4}$ be four tuples of distinct vectors 
$ (\vec r,\vec s,\vec t,\vec u)$ for which (i) $ r_1=s_1$ and $ t_1=u_1$; 
and (ii) in the second   coordinate   we have $ r _{2}=t _{2}=a$; 
and (iii) two of the four vectors agree in the third coordinate.

\begin{lemma}\label{l.Ba}  For any integer $ a$, and 
subset $ \mathbb B \subset \mathbb B _{a} (4)$ we have 
\begin{equation}\label{e.Ba}
\norm \operatorname {Prod} (\mathbb B ). p. \lesssim p n ^{5/2} \,. 
\end{equation}
Moreover, for $ \mathbb B_a \subset \mathbb B (4)\cap ( \mathbb A _s \times \mathbb A _t )^{2}$, 
for any choice of $ 0\le s\neq t\le q$, we have 
\begin{equation}\label{e.2Bq}
\norm \operatorname {Prod} (\mathbb B  ) .p. \le c  p \, n ^{5/2} q ^{-2}\,. 
\end{equation}
\end{lemma}

The point of this estimate is that we reduce the number of parameters of $ \mathbb B (4)$ 
by one, and gain a full power of $ n$ in the size of the $ L ^{p}$ norm, 
as compared to the estimate in (\ref{e.2q}). 

\begin{proof}
In the proof of Lemma~\ref{l.2}, in the analysis of the terms $ \mathbb B _3$ 
and $ \mathbb B _4$ we used the triangle inequality over the term $ b=r_2=t_2$. 
Treating this coordinate as fixed, we gain a term  $ n ^{-1}$ in the previous 
proof, hence proving the Lemma above.  The additional powers of $ q$ are 
obtained by using the fact that the first coordinates can only vary over 
a set of size $ \simeq n /q$. 

\end{proof}

A further sub-case of the inequality (\ref{e.Simple2}) demands attention.  
Using the notation of Lemma~\ref{l.SimpleCoincie}, let 
\begin{equation}\label{e.Cadef}
\mathbb C _{2,b} \coloneqq 
\{ (\vec r_1, \vec r_2)\in \mathbb C _{2} \mid r _{1,1}=b\}\,,
\qquad 1\le a \le n\,. 
\end{equation}
Thus, this collection consists of pairs of distinct vectors, with 
a coincidence in the second coordinate, and the first coordinate of $ \vec r_1$ 
is fixed. 
Note that these collections of variables have two free parameters. 
At $ L^2$ we find a $ 1/4$ gain over the `naive' estimate. 

\begin{lemma}\label{l.Simple2a} For any $ b$ 
and any subset $ \mathbb C \subset \mathbb C _{2,b}$ we have the estimates 
\begin{equation}\label{e.Simple2a}
\norm \operatorname {Prod} (\mathbb C ) .p. 
\lesssim   p \cdot  n ^{5/4}\,, \qquad 2 \le p<\infty \,.
\end{equation}
Moreover, if $ \mathbb C \subset  \mathbb A _s \times \mathbb A _t )$, 
for any choice of $ 0\le s\neq t\le q$, we have 
\begin{equation}\label{e.Simple2aC}
\norm \operatorname {Prod} (\mathbb C ) .p. 
\lesssim   p \cdot  n ^{5/4} q ^{-1/2}\,, \qquad 2 \le p<\infty \,.
\end{equation}
\end{lemma}

\begin{proof}
As in the proof of Lemma~\ref{l.SimpleCoincie}, we begin with the case $ p=2$. 
Observer that 
\begin{align*}
\norm \operatorname {Prod} (\mathbb C ) .2. ^2 
&= \mathbb E \operatorname {Prod} (\mathbb B )\,,
 \end{align*}
 where $ \mathbb B = \mathbb C _{2,b} \times \mathbb C _{2,b}\cap \mathbb B _{b} (4)$, 
 with the last collection defined in Lemma~\ref{l.Ba}. 
Therefore, the Lemma in this case follows from that Lemma.  

More generally, no pair of vectors in $ \mathbb C _{2,b} (2)$ can have a 
coincidence in the third coordinate, so we can use the Littlewood Paley 
inequalities in that coordinate to estimate 
\begin{equation*}
\norm \operatorname {Prod} (\mathbb C ) .p. 
\lesssim \sqrt p 
\NOrm \Bigl[ \sum _{c} \ABs{ \sum _{\substack{ (\vec r_1, \vec r_2) \in \mathbb C \\
\max \{r _{1,3}, r _{2,3}\} =c} } 
f _{\vec r_1} \cdot f _{\vec r_2} 
} ^2  \Bigr] ^{1/2} .p. 
\end{equation*}

Observe that 
\begin{align} \label{e.++}
\sum _{c} \ABs{ \sum _{\substack{ (\vec r_1, \vec r_2) \in \mathbb C \\
\max \{r _{1,3}, r _{2,3}\} =c} } 
f _{\vec r_1} \cdot f _{\vec r_2} 
} ^2= \sharp \mathbb C 
+ \sum _{i<j\in \{1,2,3,4\}} \operatorname {Prod} (\mathbb C _{i,j})
+ \operatorname {Prod} (\mathbb B )\,. 
\end{align}
Similar to before, we define the collections $ \mathbb C _{i,j}$ as follows. 
\begin{equation*}
\begin{split}
\mathbb C _{i,j}\coloneqq \{ (\vec r_1,\vec r_2,\vec r_3,\vec r_4)\in \mathbb C \times 
\mathbb C 
\mid  \vec r_i=\vec r_j\,, \textup{and the other two vectors are distinct}\}
\end{split}
\end{equation*}

In this case, observe that four of these collections are empty, namely 
\begin{equation*}
\mathbb C _{1,2}=\mathbb C _{2,3}=\mathbb C _{1,4}=\mathbb C _{2,3}=
\mathbb C _{2,4}\emptyset\,. 
\end{equation*}
The only non-empty collection is $ \mathbb C _{1,3}$. 
Yet, in $ \mathbb C _{1,3}$, the vectors $ \vec r_2$ and $ \vec r_4$ have 
a coincidence in the second coordinate.  Thus, 
Lemma~\ref{l.SimpleCoincie} 
applies to $ \mathbb C _{1,3}$, so that we have the estimate 
\begin{equation} \label{e.13}
\norm \operatorname {Prod} (\mathbb C _{1,3}).p. \lesssim 
p ^{5/4} n ^{7/4} \,. 
\end{equation}

Let us prove (\ref{e.Simple2a}). 
Combining these observations with (\ref{e.++}) and Lemma~\ref{l.Ba} we see that 
\begin{align*}
 p ^{-1/2}\norm \operatorname {Prod} (\mathbb C ) .p. 
&\lesssim  n + 
\norm \operatorname {Prod} (\mathbb C _{1,3}).p/2. ^{1/2} 
+ \norm \operatorname {Prod} (\mathbb B ).p/2. ^{1/2} 
\\ 
& \lesssim n +  p ^{5/8} n ^{7/8}+ p ^{1/2} n ^{5/4}\,. 
\end{align*}
Concerning the right hand side, note that for $ 2<p< n ^{3}$, 
we have $ p ^{5/8} n ^{7/8}< p ^{1/2} n ^{5/4}$. Hence we have proved 
\begin{equation*}
\norm \operatorname {Prod} (\mathbb C ) .p. 
\lesssim p n ^{5/4}\,, \qquad 1<p < n ^{3}\,. 
\end{equation*}
Yet, for $ p \gtrsim n$ the $ L ^{p}$ norm above is comparable to the $ L ^{\infty }$ 
norm, so we have finished the proof of (\ref{e.Simple2a}).  

The case of (\ref{e.Simple2aC}) is left to the reader. 

\end{proof}

\section{Norm Estimates Particular to the Hyperbolic Assumption} \label{s.norm}

The result of Theorem~\ref{t.distributional} admits
an  improvement, which we state in the 
a form adapted to our Riesz product.   These improvements are 
subtle consequences of the detailed information we have about the Beck Gain. 

\begin{theorem}\label{t.better} Using the notation of (\ref{e.rho}) and (\ref{e.G_t}), 
we have this estimate, valid for all $ 1\le t\le q$. 
\begin{equation}\label{e.better}
\norm  \rho F_t  .p.  \lesssim  \sqrt p \,  \,, 
\qquad 1\le p \le 
c  n ^{1/3} \,. 
\end{equation}
As a consequence, we have the distributional estimate 
\begin{equation}\label{e.partialX2}
\mathbb P ( \rho G _t >x ) \lesssim \operatorname {exp} (- c x ^2 )\,, 
\qquad x< c  n ^{1/6}\,.
\end{equation}
Here $ 0<c<1$ is an absolute constant. 
\end{theorem}

\begin{remark}\label{r.notneeded} It is perhaps worth emphasizing that we do not 
need this Theorem to deduce our main result, Theorem~\ref{t.bl} on the Small Ball 
Conjecture in three dimensions.\footnote{If one does not use the result above, 
a smaller value of $ b= \tfrac13$ is required.}
Nevertheless, we will use the result above.  
And we find the proof to be a compelling application of the Beck Gain. 
\end{remark}

\begin{remark}\label{r.futher?}  
There are limits to  validity to these kinds of inequalities:  Recall that one has $ \ell ^{\infty } \simeq \ell ^{
\log N}$.  Thus, for appropriate $ F_t$ we would have 
\begin{equation*}
\norm \rho F_t . \infty . \simeq \norm \rho  F_t . 3n . \simeq  n/ \sqrt q   \,   .  
\end{equation*}
Hence, the sub--gaussian bound above can't hold for this range of $ p$, 
unless $ q \simeq n$, but then the sub--gaussian estimate is immediate.  
\end{remark}

\begin{proof}

Recall that 
\begin{equation*} 
F_t =  \sum _{\vec r\in \mathbb A _t}  f _{\vec r}\,. 
\end{equation*}
where  $ \mathbb A _t \coloneqq \{\vec r\in \mathbb H _n 
\mid r_1\in I_t\}$, and $ I_t$  in an interval of integers of length $ n/q$, 
so that $ \sharp \mathbb A _t \simeq   \rho ^2 $, with $ \rho $ defined in (\ref{e.rho}).

Apply the Littlewood Paley inequality  in the first coordinate. 
This results in the estimate 
\begin{align*}
\norm \rho F_t .p. 
& \lesssim \sqrt p 
\NOrm \Bigl[\sum _{s\in I_j} \ABs{ \rho \sum _{\vec r\,:\, r_1=s} f _{\vec r}} ^2  
\Bigr] ^{1/2} .p. 
\\
& \lesssim  \sqrt p  
\norm   1+  \Gamma _t   .p/2. ^{1/2} 
\\
& \lesssim  \sqrt p \Bigl\{ 1 + \norm \Gamma _t   .p/2. ^{1/2}\Bigr\}
\\
\Gamma _t & \coloneqq \rho ^2 \sum _{\substack{\vec r\neq \vec s\in \mathbb A _t\\ 
r_1=s_1}} f _{\vec r} \cdot f _{\vec s}\,. 
\end{align*}
Of course the terms $ \Gamma _t$ are controlled by the estimate in (\ref{e.Simple2q}). 
In particular, we have 
\begin{equation}\label{e.FFFjjj}
\norm \Gamma _t.p. \le C p  ^{3/2} n ^{-1/2}\,. 
 \end{equation}
Hence (\ref{e.better}) follows.
 
\smallskip  

The second distributional inequality is a well known consequence of the norm 
inequality.  Namely, one has the inequality below, valid for all $ x$:
\begin{equation*}
\mathbb P (\rho F_t >x ) \le C ^{p} p ^{p/2} x ^{-p}\,, 
\qquad 1\le p \le  c n ^{1/3}\,. 
\end{equation*}
If $ x$ is as in (\ref{e.partialX2}), we can take $ p \simeq x ^{2}$ to prove the 
claimed exponential squared bound. 
 \end{proof}

\begin{remark}\label{r.beyond1/5} The proof above does permit better than `naive' estimates 
for $ \norm  \rho F_t .p.$ for a range of $ p> n ^{1/3}$. 
The estimate we have is 
\begin{equation*}
\norm \rho F_t.p. \lesssim \min  \{ p\,,\ 
\sqrt p (1+ p ^{3/2} n ^{-1/2})\}\,.
\end{equation*}
The first estimate is 
from Theorem~\ref{t.distributional} while the second estimate is from the proof above. 
The minimum will be the second estimate provided $ p \lesssim n ^{1/2}$. 
Thus, for $ n ^{1/3} < x <  n ^{1/2}$ one can achieve an estimate that is better 
than from that of Theorem~\ref{t.distributional}. 
\end{remark}

We now prove a central estimate of the proof.

\begin{lemma}\label{l.2b} The estimate (\ref{e.expq2b}) holds.  Moreover, we have 
\begin{equation}\label{e.22bb}
\sup _{V\subset \{1 ,\dotsc, q\}}
\mathbb E \prod _{v \in V}  (1+ \widetilde \rho F_t) ^2 \lesssim 
\operatorname {exp} (a' q ^{2b})\,. 
\end{equation}
Here, $ \widetilde \rho $ is as in (\ref{e.rho}), and $ a'$ is a fixed constant times 
$ 0<a<1$, the small constant that enters into the definition of $ \widetilde \rho $. 
\end{lemma}

\begin{remark}\label{r.condExpect} A conditional expectation argument is essential 
to this proof.  This Lemma is also proved in Beck's paper. Yet, due to a more complicated 
Riesz product, the use of our line of reasoning was not available to him.
\end{remark}

\begin{proof}
The supremum over $ V$ will be an immediate consequence of the proof below, and
so we don't address it specifically.

Let us give the initial, essential observation.
We expand
\begin{equation*}
\mathbb E \prod _{v=1} ^{q} (1+ \widetilde \rho F_t) ^2
=
\mathbb E \prod _{v=1} ^{q} (1+ 2\widetilde \rho F_t+ (\widetilde \rho F_t) ^2  )\,.
\end{equation*}
Hold the $ x_2$ and $ x_3$ coordinates fixed, and let $ \mathcal F$ be the sigma
field generated by $ F_1 ,\dotsc, F_{q-1}$.  We have
\begin{equation} \label{e.;p}
\begin{split}
\mathbb E (1+ 2\widetilde \rho F_q+ (\widetilde \rho F_q) ^2
\,\big|\, \mathcal F) &=1+\mathbb E ((\widetilde \rho F_q) ^2
\,\big|\, \mathcal F)
\\
&=1+   a ^2  q ^{2b-1}+ \widetilde \rho ^2  \Gamma _q\,,
\\ \text{where}\,\,
\Gamma _t & \coloneqq \sum _{\substack{ \vec r\neq \vec s \in
\mathbb A _{t}\\ {r_1=s_1} }} f _{\vec r} \cdot f _{\vec s},
\end{split}
\end{equation}
 Then, we see that
\begin{align} \nonumber
\mathbb E \ \prod _{v=1} ^{q} (1+
2\widetilde \rho F_t+ (\widetilde \rho F_t) ^2  ) &= \mathbb E
\Bigl\{  \prod _{v=1} ^{q-1} (1+ 2\widetilde \rho F_t+ (\widetilde
\rho F_t) ^2  )\,. \times \mathbb E (1+ 2\widetilde \rho F_t+
(\widetilde \rho F_t) ^2 \,\big|\, \mathcal F) \Bigr\}
\\
&\le  \label{e.;;}
(1+a ^2 q ^{2b-1})
\mathbb E \prod _{v=1} ^{q} (1+ 2\widetilde \rho F_t+ (\widetilde \rho F_t) ^2  )
\\  \label{e.;;;;}
& \qquad + \Abs{\mathbb E
\prod _{v=1} ^{q} (1+ 2\widetilde \rho F_t+ (\widetilde \rho F_t) ^2  )
\cdot
 \widetilde \rho ^{2} \Gamma _q }
\end{align}
This is the main observation: one should induct on (\ref{e.;;}),
while treating the term in (\ref{e.;;;;}) as an error, as the `Beck
Gain' estimate (\ref{e.Simple2q}) applies to it.

Let us set up notation to implement this line of approach.  Set
\begin{equation*}
N (V;r) \coloneqq
\NOrm
\prod _{v=1} ^{V} (1+ \widetilde \rho F_t) .r.  \,, \qquad   V=1 ,\dotsc, q\,.
\end{equation*}
We will use the trivial inequality available from the exponential 
moments
\begin{align*}
N (V; 4)&\le
\prod _{v=1} ^{V} \norm 1+ \widetilde \rho F_t .4V.
\\
& \le ( 1 +  C q ^{b-1/2} V ) ^{V}
\\
& \le (Cq) ^{Cq}\,.
\end{align*}

This of course is a terrible estimate, but we now use interpolation, noting that
\begin{equation}\label{e.killq^q}
N (V;2(1-1/q) ^{-1} )\le N (V;2) ^{1-1/q} \cdot N (V; 4) ^{1/q}\,.
\end{equation}

We see that (\ref{e.;;}), (\ref{e.;;;;}) and (\ref{e.killq^q})
give us the inequality
\begin{equation}\label{e.==}
\begin{split}
N (V+1; 2) & \le (1+a ^2 q ^{2b-1}) ^{1/2}  N (V; 2)
+ C  \cdot  N (V; 2 (1-1/q) ^{-1} )  \cdot
\norm  \widetilde \rho ^{2} \Gamma _q . q.
\\
& \le
(1+a ^2 q ^{2b-1}) ^{1/2}  N (V; 2)
+ C  N (V;2) ^{1-1/q} \cdot N (V; 4) ^{1/q}
\norm  \widetilde \rho ^{2} \Gamma _q . q.
\\
& \le
(1+a ^2 q ^{2b-1}) ^{1/2}  N (V; 2)
+ C  q ^{C} n ^{-1/2} N (V;2) ^{1-1/q} \,.
\end{split}
\end{equation}
In the last line we have used the the inequality (\ref{e.Simple2q}).

Of course we only apply this as long as $ N (V; 2)\ge 1$.  Assuming this is true
for all $ V\ge 1$, we see that
\begin{align*}
N (q;2)&\lesssim  (1+a ^2  q ^{2b-1}+ C q ^{C} n ^{-1/2}) ^{q}
\\
& \lesssim \operatorname e ^{ a' q ^{2b}}\,.
\end{align*}
Here  of course we need $ C q ^{C} n ^{-1/2}\le  a q ^{2b-1} $, which we
certainly have for large $ n$. 

\end{proof}

\section{The Beck Gain} \label{s.nsd}

Let us state the main result of this section. 
Given $ V\subset \{ 1 ,\dotsc, q\}$ let 
\begin{align*}
\operatorname {NSD} (V)
\coloneqq 
\Bigl\{ \{\vec r_j \mid j\in V\}\in \operatorname \times _{j\in V} \mathbb A _{j} 
\; \big|\; 
&\textup{for each $ j\in V$, there is a choice of $ j'\in V- \{j\}$}
\\
& \qquad \textup{and $ \ell =2,3$ so that $ r _{j,\ell }=r _{j',\ell }$}
\Bigr\} \,.
\end{align*}
That is, we take tuples of $ \mathsf r$ vectors, indexed by $ V$, requiring that 
each $ \vec r_j$ be in a coincidence.    
Such sums admit a favorable estimate  on their $ L^2$ norms. 

\begin{theorem}\label{t.NSD}[\textup{\textbf{The Beck Gain.}}] 
There are positive constants $ C_0, C_1, C_2, C_3, \eta   $ for which we have the estimate 
\begin{equation}\label{e.NSD}
\rho ^{\abs{ V}} \Norm 
\operatorname {Prod} (\operatorname {NSD} (V)) .p. 
\lesssim  [C_0 \abs V ^{C_1}  p ^{C_2}  q ^{C_3} n ^{-\eta }] ^{\abs{ V}} 
\,, \qquad V\subset \{1 ,\dotsc, q\}\,\ 
 \,.
\end{equation}
\end{theorem}

\begin{remark}\label{r.gaininV} The novelty in this estimate is that 
we find that  (a) the gain is proportional to the number of vertices , 
and (b) the gain also holds in $ L ^{p}$ norms.  In application, 
$ p \lesssim  q ^{2b}= \sqrt q \simeq n ^{\epsilon '}$, so the polynomial 
growth in $ p$  and in $ q$ is acceptable to us.   Beck \cite{MR1032337}
found a gain in $ L ^{2}$ norm of order $ n ^{-1/4}$, for all $ V$. 
Such a small gain of course forces a much shorter Riesz product. 
\end{remark}

\begin{remark}\label{r.eta}  It is disappointing that we cannot identify a 
reasonable value of $ \eta >0$, which is in large measure, the amount of the gain. 
Yet, the goal of this proof is to have a relatively simple method of proof. 
  Obviously, a finer understanding of this estimate, 
among other issues, will be central to future progress on the range of 
questions discussed in these notes. 
\end{remark}

The proof of this Theorem requires a careful analysis of the 
variety of ways that a product can fail to be strongly distinct.  
That is, we need to understand the variety of ways that coincidences 
can arise, and how coincidences can contribute to a smaller  norm.

  It is important 
at the outset to recognize that  patterns of coincidences can be quite complex, 
a point best illustrated by a few examples of such patterns. 
Consider the specific product 
\begin{equation} \label{e.PROD}
\prod _{j=1} ^{7} \sum_{ \substack{\vec r_j \in \mathbb A _{j}\\ }} 
 f _{\vec r_j}\,,  
\end{equation}
and the ways that summands in such a product could fail to be strongly distinct.
One could consider those terms in which the first three choices of $ \vec r_j$
agree in the second coordinate: 
\begin{equation*}
r_ {1,2}= r _{2,2}= r _{3,2}\,,
\end{equation*}
while imposing no restriction on the remaining four vectors $ \vec r_4, \vec r_5, 
\vec r _6, \vec r_7$.  Note that 
\begin{equation} \label{e.Xprod}
\prod _{j=1} ^{7} \sum_{ \substack{\vec r_j \in \mathbb A _{j}\\
r_ {1,2}= r _{2,2}= r _{3,2}}} 
 f _{\vec r_j}
 =\Biggl[ \prod _{j=1} ^{3} \sum_{ \substack{\vec r_j \in \mathbb A _{j}\\
r_ {1,2}= r _{2,2}= r _{3,2}}}  f _{\vec r_j}\Biggr]
\cdot 
\Biggl[ \prod _{j=4} ^{7} \sum_{ \substack{\vec r_j \in \mathbb A _{j}\\ }} 
 f _{\vec r_j}
 \Biggr]
\end{equation}
That is, we have a product of terms, with a `simple' coincidence in the 
first term, and no restriction on the sum in the second. In this instance, 
we would take $ V= \{1,2,3\}$. 

Similarly, a pattern of coincidences could be 
\begin{equation*}
r_ {1,2}= r _{2,2}= r _{3,2}\,,
\qquad 
r_ {4,3}= r _{5,3}= r _{6,3}= r _{7,3}\,. 
\end{equation*}
As in the first case, the corresponding sum would break into a product. 
And the $ L ^{1}$ norm would be substantially smaller, due to the presence of 
two sets of `simple' coincidences. 

Yet, one could have a more complicated set of coincidences, such as 
\begin{equation*}
r_ {1,2}= r _{2,2}= r _{3,2}\,,
\qquad 
r_ {1,3}= r _{4,3}= r _{5,3} \,, 
\qquad 
r _{2,3}= r _{6,3}= r _{7,3}\,.
\end{equation*}
Here, the first and second vectors are both involved in two distinct sets of coincidences. 
This case, as it turns out, are also substantially smaller in $ L ^{1}$ norm 
than the first case, due to the `overlapping' coincidences. 

Following Beck, we will use the language of Graph Theory to describe these 
general patterns of coincidences.

Before passing to the general description of these results, the reader should keep 
forefront in their minds these points: 
\begin{itemize}
\item Coincidences can only occur in the second and third coordinates, due to the 
specific way we form our products.  
\item  Our graphs will have as vertices the   integers $ j\in \{1,2 ,\dotsc, q\}$, the index 
of the product in (\ref{e.PROD}). 
\item  Edges in the  graph represent a coincidence between two vectors.  
Edges come in two different types,  or colors, associated to coincidence in the second 
or  third coordinates. 
\item  Equality is transitive, so the edges in e.\thinspace g.\thinspace 
the second coordinate will naturally decompose into \emph{cliques}. 
\item  As we work in three dimensions, a clique in the second coordinate, 
and a clique in the third coordinate can contain at most one common vertex, as 
two common vertices  would imply that our product contains two equal vectors. 
This case is specifically excluded from our consideration. 

\item The presence of an edge will mean that we enforce a coincidence of that type 
in the products we consider.  The absence of an edge will mean that no such condition 
is assumed---not that equality is forbidden.
This will permit product formulas such as (\ref{e.Xprod}) above hold. 

\item  A graph is naturally associated to sums of products of $ \mathsf r$ functions.
We seek effective $ L ^{p}$ norms on these sums.  Larger cliques, and more overlapping 
cliques serve to reduce the number of parameters, and give \emph{smaller} norms. 
\end{itemize}

\subsection*{Graph Theory Nomenclature}

We adopt familiar nomenclature from Graph Theory.\footnote{There is no graph theoretical fact 
that we need, rather the use of this language is just a convenient way to 
do some bookkeeping.}
The class of graphs 
that we are interested satisfy particular properties.  
A \emph{graph} $ G$ is the triple of $ (V (G), E_2, E_3)$, of the 
\emph{vertex set} $ V (G)\subset \{1 ,\dotsc, q\}$, and \emph{ edge sets $ E_2$ 
and $ E_3$, of color $ 2$ and $ 3$ respectively}.  Edge sets are 
are subsets of 
\begin{equation*}
E _{j}\subset V (G) \times V (G) - \{ (k,k)\;|\;  k\in V (G)\}\,.
\end{equation*}
Edges are symmetric, thus if $ (v,v')\in E_j$ then necessarily $ (v',v)\in E_j$. 

A \emph{clique of color $ j$} is a maximal subset $ Q\subset V (G)$ 
such that for all $ v\neq v'\in Q$ we have $ (v,v')\in E_j$.  By \emph{maximality}, 
we mean that no 
strictly larger set of vertices  $ Q'\supset Q$ satisfies this condition.  

Call a graph $ G$ \emph{admissible} iff 
\begin{itemize}
\item  The edges sets, in both colors, decompose into a union of cliques. 
\item  Any two cliques  $ Q_2$ in color $ 2$ and clique $ Q_3$ in color $ 3$ 
can contain at most one common vertex. 
\item Every vertex is in at least one clique. 
\end{itemize}

 A graph $ G$ is \emph{connected } iff for any two vertices  
in the graph, there is a path that connects them.
A \emph{path} in the graph $ G$ is a sequence of vertices  $ v_1 ,\dotsc, v_k$ 
with an edge of \emph{either color,} spanning adjacent vertices , that is   $ (v _{j}, v _{j+1}) 
\in E _2 \cup E_3 $.

\subsection*{Reduction to Admissible Graphs}

Given admissible graph $ G$ on vertices $ V$, we set $ X (G)$ 
to be those tuples of $ \mathsf r$ vectors 
\begin{equation*}
\{\vec r _v \mid v\in V\}\in \prod _{v\in V} \mathbb A _v
\end{equation*}
so that if $ (v, v')$ is an edge of color $ j$ in $ G$, then $ r _{v,j}= r _{v',j}$. 

We will prove the Lemma  below in the following two sections.  
\begin{lemma}\label{l.admissible<} For an admissible graph $ G$ on vertices $ V$ we have 
the estimate below for positive, finite constants $ C_0, C_1, C_2, C_3$: 
\begin{equation}\label{e.equiv}
\rho ^{\abs{ V}} \norm \operatorname {Prod} (X (G)).1. 
 \le  [C_0\abs{ V} ^{C_1} p ^{C_2} q ^{C_3}   n ^{-\eta }] ^{ \abs{ V}}\,, 
\qquad 2<p< \infty  \,.
\end{equation}
\end{lemma}

Let us give the proof of Theorem~\ref{t.NSD} assuming this Lemma. 
Our tool is the Inclusion Exclusion Principle, but to apply it we need additional 
concepts.  

Given two admissible graphs $ G_1, G_2$ on the same vertex set $ V$, let 
$ G_1 \wedge G_2$ be the smallest admissible graph which contains all the edges 
in $ G_1$ and in $ G_2$.  By smallest, we mean the graph with the fewest number of 
edges; and such a graph may not be defined, in which case we take $G_1 \wedge G_2 $ 
to be undefined.  We recursively define $ G_1\wedge \cdots \wedge G_k \coloneqq 
(G_1 \wedge \cdots G _{k-1}) \wedge G_k$.  This wedge product is associative. 

Let $ \mathcal G_0$ be the set admissible graphs on $ V$ which are \emph{not} 
of the form $ G_1\wedge G_2$ for admissible $ G_1, G_2$. These are the `prime' graphs. 
(If $ V$ is of cardinality $ 2$ or $ 3$, every graph is prime.)  
For instance, in the case of $ V= \{v_1,v_2,v_3,v_4\}$ the two graphs 
below are prime. 
\begin{equation*}
\begin{array}{cccccccc} 
v_1 &  & v_2 && v_3 && v_4 
\\
\Box  &   & \Box  && \Box  && \Box 
\\
\bullet & = & \bullet  &    &   \bullet & = & \bullet 
\\
	&   &    &  & 
\end{array}
\quad \textup{and}
\quad 
\begin{array}{cccccccc} 
v_1 &  & v_4 && v_2 && v_3 
\\
\Box  &   & \Box  && \Box  && \Box 
\\
\bullet & = & \bullet  &    &   \bullet & = & \bullet 
\\
	&   &    &  & 
\end{array}
\,. 
\end{equation*}
The only difference between the two is the ordering of the vertices in the top 
row.  There are no coincidences in the third row, and the first row, with the 
$ \Box$s, never has a coincidence.  These two graphs are distinct, and clearly 
members of $ \mathcal G_0$.  Note that their wedge product is 
\begin{equation*}
\begin{array}{cccccccc} 
v_1 &  & v_2 && v_3 && v_4 
\\
\Box  &   & \Box  && \Box  && \Box 
\\
\bullet & = & \bullet  & =   &   \bullet & = & \bullet 
\\
	&   &    &  & 
\end{array}
\end{equation*}

Now define $ \mathcal G_k$ to be those graphs which are equal to 
a wedge product $ G_1 \wedge \cdots \wedge G_k$, with $ G_j\in \mathcal G_0$, 
and moreover, $ k$ is the smallest integer for which this is true. 
Clearly, we only need to consider $ k\le q$.

Then, by the inclusion exclusion principle, 
\begin{equation}\label{e.INexclude}
\operatorname {Prod} (\operatorname {NSD} (V)) 
=\sum _{k=0} ^{q} (-1) ^{k}  \sum _{G\in \mathcal G _{k}}
\operatorname {Prod} (X (G))\,. 
\end{equation}
The number of admissible graphs 
on a set of vertices $ V$ is at most $ 2 ^{\abs{ V}}  \abs{ V}!< 2 ^{\abs{ V}} \abs V ^{\abs 
V}$.
So that using (\ref{e.equiv}) clearly implies Theorem~\ref{t.NSD}.

\subsection*{Norm Estimates for Admissible Graphs}

We begin this section with a further reduction to \emph{connected} admissible graphs. 
Let us write $ G\in \operatorname {BG} (C_0, C_1,C_2, C_3 ,\eta )$ if the 
estimates (\ref{e.equiv}) holds. 
(`$ \operatorname {BG}$' for `Beck Gain.') 
We need to see that all admissible graphs are in $ \operatorname {BG} (C_0, C_1,C_2,
C_3 ,\eta)$ 
for non-negative, finite choices of the relevant constants.

\begin{lemma}\label{l.holder} Let $ C_0, C_1,C_2, C_3, \eta$ 
be non-negative constants. Suppose that $ G$ is an admissible graph, 
and that it can be written as a union  subgraphs $ G_1 ,\dotsc, G_k$ where 
all $ G_j \in \operatorname {BG} (C_0, C_1,C_2, C_3, \eta)$. Then, 
\begin{equation*}
G \in \operatorname {BG} (C_0, C_1, C_2, C_2+C_3, \eta )\,. 
\end{equation*}
\end{lemma}

With this Lemma, we will identify a small class of graphs for which we 
can verify the property (\ref{e.equiv}) directly, and then appeal to this 
Lemma to deduce Theorem~\ref{t.NSD}.   Accordingly, we modify our notation.  
If $ \mathcal G$ is a class of graphs, we write $ \mathcal G\subset \operatorname {BG} (\eta )$ 
if there are constants $ C_0, C_1, C_2, C_3 $ such that 
$ \mathcal G\subset \operatorname {BG} (C_0, C_1, C_2, C_3, \eta )$.

\begin{proof}
We then have  by Proposition~\ref{p.products}
\begin{equation*}
\operatorname {Prod} (X (G)) 
= 
\prod _{j=1} ^{k}
\operatorname {Prod} (X (G_j)) 
\,.  
\end{equation*}
Using H\"older's inequality, we can estimate 
\begin{align*}
\norm \operatorname {Prod} (X (G)) .p. 
&\le \prod _{j=1} ^{k} \norm \operatorname {Prod} (X (G_j)). k p. 
\\
& \le \prod _{j=1} ^{k}  [C_0 (kp) ^{C_1} q ^{C_2} n ^{-\eta } ] ^{\abs{ V_j}}
\\
& \le [C_0 p ^{C_1} q ^{C_2+C_1} n ^{-\eta } ] ^{\abs{ V}}\,. 
\end{align*}
Here, we use the fact that since the graphs are non-empty, we necessarily have 
$ k\le q$. 

\end{proof}

\begin{proposition}\label{p.products} 
Let $ G_1 ,\dotsc, G_p$ be admissible graphs on pairwise disjoint
vertex sets $ V_1 ,\dotsc, V_p $.  Extend these graphs in the natural way to 
a graph $ G$ on the 
vertex set $ V=\bigcup V_t$.  Then, we have 
\begin{equation*}
\operatorname {Prod} (X(G))
=
\prod _{t=1} ^{p} \operatorname {Prod} (X (G_t))\,. 
\end{equation*}
\end{proposition}

\subsection*{Connected Graphs Have the Beck Gain.}

  Let $ \mathcal G _{\textup{connected}} $
be the collection of 
of all admissible  connected 
 graphs on $ V\subset \{1,\dotsc, q\}$.

\begin{lemma}\label{l.twoCliques}    
We have $ \mathcal G _{\textup{connected}}\subset \operatorname {BG} ( \tfrac1 {15} )$. 
\end{lemma}

One can depict small examples of  these graphs as follows.  
\begin{gather*}
\begin{array}{ccc} 
\Box  &   & \Box  
\\
\\
\bullet & = & \bullet 
\end{array}
\, ; \qquad 
\begin{array}{ccccc} 
\Box  &   & \Box  && \Box 
\\
\bullet & = & \bullet 
\\
	&   &  \bullet & = & \bullet 
\end{array}
\, ; \qquad 
\begin{array}{cccccccc} 
\Box  &   & \Box  && \Box  && \Box 
\\
\bullet & = & \bullet  &    &   \bullet & = & \bullet 
\\
	&   &  \bullet & = & \bullet 
\end{array}
\end{gather*}
These are graphs on $ 2$, $ 3$ and $ 4$ vertices  respectively.  We will have 
to pay special attention to the case of $ 2$ and $ 3$ vertices, as 
these cases are not amenable to the general procedure we invoke below. 
It is important to observe that the first coordinates, represented by $ \Box$ 
above, are necessarily distinct, and have the partial order inherited from the 
vertex set $ V$.   Namely, the vertex set $ V \subset \{1 ,\dotsc, q\}$, 
 and $ V$ inherits the order from  the integers.   By the construction 
 of our Riesz product, the first coordinates inherit this same order.

 Unfortunately, even working with this class of admissible graphs, our 
proof is of an \emph{ad hoc} nature, and we won't actually specify a value 
of $ \eta >0$ for which the Lemma above holds.  
 
\subsubsection*{General Remarks on Littlewood Paley Inequality.}
These  remarks are essential to our analysis of this lemma, and  
the Theorem we are proving. 
The vertex set $ V$ is a subset of $ \{1 ,\dotsc, q\}$ and it inherits 
an order from that set.  Moreover, the tuples of $ \mathsf r$ vectors do as well. 
Namely, writing 
\begin{equation*}
V= \{v_1< \cdots <v _\ell \}, 
\end{equation*}
for $ \{r_1 ,\dotsc, r_\ell \} \in X (G)$, we have, by construction, 
$ r _{1,1} < \cdots <r _{\ell ,1}$. This since $ r _{m,1} \in I _{v_m}$, 
where $ I _{m'}$ is the increasing sequence of intervals of length equal to 
$n/q$ that partition $ \{1 ,\dotsc, n\}$. 

Continuing this line of thought, we see that there is a natural way to 
apply the Littlewood Paley inequalities. 
For integer $ b _{\ell }\in I _{\ell }$, let $ X (G;  \vec b _{\ell })$ 
be the tuple of $ \mathsf r$ vectors $ \{\vec r_1 ,\dotsc, \vec r _{\ell }\}$ 
such that  $ r _{\ell ,1}=b _{\ell }$.  We have 
\begin{equation} \label{e.LPN}
\norm \operatorname {Prod} (X (G)) .p. 
\lesssim \sqrt p 
\NOrm \Bigl[ \sum _{b _{\ell } \in I _{ v _\ell }} 
\abs{ \operatorname {Prod} (X (G;  b _{\ell }))} ^2 
\Bigr] ^{1/2} .p. \,. 
\end{equation}
It is tempting to continue this procedure, by applying the Littlewood Paley 
inequality again to the vertex $ v _{\ell -1}$.  Yet---and this in an important point---due 
to the nature of  $ \mathsf r$ functions, this option is blocked to us.  
The vertex $ v _{\ell }$ is in at least one clique $ Q $ 
of, say, color  $ 2$.  We could 
choose a value $ c _{Q}$ for that clique, thereby specifying all 
coordinates of the vector $ \vec r _{\ell }$.  Set  
 $ X (G; b _{\ell } ; c _{Q })$ 
be the tuple of $ \mathsf r$ vectors $ \{\vec r_1 ,\dotsc, \vec r _{\ell-1 }\}$ 
such that   
\begin{equation*}
\{\vec r_1 ,\dotsc, \vec r _{\ell-1 }\,,\ (a _{\ell },b _{\ell }, n - a _{\ell }- b _{\ell } )\}
\in X (G; a _{\ell })\,. 
\end{equation*}
Here, $ X (G;  b _{\ell }; c _{q })$
consists of tuples of length $ \ell -1$, since 
the vector $ \vec r _{\ell }$ is completely specified.  Thus, we see that 
\begin{equation} \label{e.3halfs}
\norm \operatorname {Prod} (X (G)) .p. 
\lesssim n  \sup _{a _{\ell }, c_Q} 
\NOrm \Bigl[ \sum _{b_ \ell } \operatorname {Prod} (X (G;   b _{\ell } ; c_Q)) ^2 
\Bigr] ^{1/2} .p. \,.  
\end{equation}
At this point, the (Hilbert space) Littlewood Paley inequalities will again apply.

We will refer to the notation above.  Keep in mind that   $ \vec b$ is for the 
coordinates specified by a Littlewood Paley inequality;  $ \vec c$ are for the 
coordinates in a coincidence that we use the triangle inequality on. 
We shall return to these themes momentarily.

\begin{proof}[Proof of Lemma~\ref{l.twoCliques}.]  
We begin the proof with a 
discussion of the case of two and three vertices , which will 
not be susceptible to the general methods related to the Littlewood Paley 
inequality outlined above.

\subsubsection*{The Case of Two Vertices .}
Notice that if $ G$ consists of only two vertices , the relevant estimate is 
(\ref{e.Simple2q}). Namely, we have 
\begin{equation*}
\norm \operatorname {Prod} (X (G)) . p. 
\le C p ^{3/2} n ^{3/2} q ^{-1}\,. 
\end{equation*}
 Equivalently, $ G\in \operatorname {BG} (C_0, 3/4,0, 1/4)$.

\subsubsection*{The Case of Three Vertices }

The case of $ G\in \mathcal G_2$ having three vertices  depends critically 
on the same phenomena behind the Beck Gain for graphs on two vertices .  We 
will deduce this case as a corollary to the case of two vertices .

There are two distinct sub-cases.  The more delicate of the two cases 
is as follows.  The graph is depicted as 
\begin{equation}\label{e.badOrder}
\begin{array}{cccccc} 
v_1 && v_2 && v_3   
\\
\Box  &   & \Box  && \Box  
\\
\bullet & = & \bullet  &    &   
\\
	&   &  \bullet & = & \bullet 
\end{array}
\end{equation}
where $ v_1<v_2<v_3$. (The case of $ v_2<v_1<v_3$ is entirely the same, and we don't 
discuss it directly.)  

By our general remarks on the Littlewood Paley inequality, this 
inequality applies in the first coordinate, to the vertex $ v_3$. 
Using the notation in (\ref{e.LPN}), we have 
\begin{equation*}
\norm \operatorname {Prod} (X (G)) .p. 
\lesssim \sqrt p 
\NOrm \Bigl[ \sum _{b_3 \in I _{v_3 }} 
\abs{ \operatorname {Prod} (X (G;   b _{3 }))} ^2 
\Bigr] ^{1/2} .p. \,. 
\end{equation*}
The vectors $ v_2$ and $ v_3$ have a coincidence in the third coordinate. 
Therefore, we specify the value of the coincidence to be $ c_3$ and estimate 
\begin{equation} \label{e.badOrder3/2}
\norm \operatorname {Prod} (X (G)) .p. 
\lesssim  \sqrt p \cdot   n ^{3/2} \cdot  \sup _{b_3,c_3}
\norm 
\operatorname {Prod} (X (G;  b _{3 }; c_3)) .p. \,. 
\end{equation}

Recall that 
$ X (G;  a _{3 }; b_3)$ consists only of pairs of vectors.  This graph 
can be depicted as 
\begin{equation*}
\begin{array}{cccc} 
v_1 && v_2 
\\
\Box  &   & \Box  
\\
\bullet & = & \bullet 
\\
	&   &  c_3
\end{array}
\end{equation*}
But this is the case considered in  (\ref{e.Simple2aC}). From that inequality,
we see that we have the estimate 
\begin{equation*}
\norm \operatorname {Prod} (X (G;  b _{3 }; c_3)) .p. \lesssim \sqrt p n ^{5/4} q ^{-1/2}\,. 
\end{equation*}
Therefore,  from (\ref{e.badOrder3/2}), we see that 
\begin{equation} \label{e.badOrder+}
\norm \operatorname {Prod} (X (G)) .p. 
\lesssim p ^{3/2} n ^{9/4} q ^{-3/2}  \,. 
\end{equation}
Recall that the point of comparison is to $ n ^{3} q ^{-3/2} $, 
and the estimate above is smaller by $ n ^{-1/4}$.
Thus 
the class of graphs given by (\ref{e.badOrder}) 
are contained in $ \operatorname {BG} ( \tfrac1{12}- \epsilon )$. 

\medskip 

The other case is when the graph can be depicted by 
\begin{equation*}
\begin{array}{cccccc} 
v_1 && v_3 && v_2   
\\
\Box  &   & \Box  && \Box  
\\
\bullet & = & \bullet  &    &   
\\
	&   &  \bullet & = & \bullet 
\end{array}
\end{equation*}
where $ v_3$, the maximal index is in both cliques.  This case is much easier, as 
one application of the Littlewood Paley inequality, and the triangle inequality 
will determine the value of both cliques.   It is very easy to see that 
this class of graphs is in $ \operatorname {BG} (1/6)$, and the details 
are omitted. Hence the discussion 
graphs on three vertices , with all cliques of size $ 2$ is complete.

\subsubsection*{A General Estimate}

We now present a general recursive estimate for the $ L ^{p}$ norm of 
$ \operatorname {Prod} (X (G))$, assuming that $ G$ is a connected graph on 
at least four vertices.  Write $ V$ as 
\begin{equation*}
V= \{v_1< \cdots  < v_\ell \}\,.
\end{equation*}

The estimate is obtained recursively.  Along the way we will construct 
two disjoint subsets $ V _{3/2}, V _{1/2} \subset V$.  
$ V _{3/2}$ will be the vertices  to which we apply both the Littlewood Paley 
and triangle inequalities, thus these vertices  contribute  $ n ^{3/2} q ^{-1/2}$ 
to our estimate.  $ V _{1/2}$ will be the vertices  to which we apply only 
the Littlewood Paley inequality, thus these vertices  contribute $ (n/q) ^{1/2} $
 to our estimate.  Those vertices  not in $ V _{3/2} 
 \cup V _{1/2}$ will be those which are determined by earlier steps in the procedure. 
 They contribute nothing to our estimate.  In estimating an $ L ^{p}$ 
 norm, the power of $ p$ is one-half of the number of applications of the 
 Littlewood Paley inequality, namely $ \tfrac 12 \sharp (V _{3/2} \cup V _{1/2})$. 
 
 The purpose of these considerations is to prove the estimate 
 \begin{equation} \label{e.purpose}
\norm 
\operatorname {Prod} (X (G)) .p. 
\le (C \sqrt p) ^{\abs{ V _{3/2}}+ \abs{ V _{1/2}}}
(n/q) ^{ (\abs{ V _{3/2}}+ \abs{ V _{1/2}})/2 } n ^{\abs{ V _{3/2}}} \,. 
\end{equation}

Initialize 
\begin{equation*}
V _{3/2} \leftarrow \emptyset \,, 
\qquad 
V _{1/2} \leftarrow \emptyset \,, 
\quad 
\mathcal Q _{\textup{fixed}} \leftarrow \emptyset \,. 
\end{equation*}
The last collection consists of those cliques which are specified by 
earlier stages of the argument.  

At each stage, we will have an estimate for the form 
\begin{equation}\label{e.Estimate}
\begin{split}
\norm 
\operatorname {Prod} (X (G)) .p. 
&\le (C \sqrt p) ^{\abs{ V _{3/2}}+ \abs{ V _{1/2}}} 
n ^{\abs{ V _{3/2}}} 
\\
& \qquad \times \sup _{\vec c \in \{ 1 ,\dotsc, n \} ^{\mathcal Q _{\textup{fixed}}}}
\NOrm \Bigl[ \sum _{\vec b \in \{1 ,\dotsc, n \} ^ { V _{3/2}\cup V _{1/2}}} 
\operatorname {Prod} (X (G;  \vec b; \vec c)) ^2  \Bigr] ^{1/2} 
.p.
\end{split}
\end{equation}
Here, $ X (G;  \vec b; \vec c)$ denotes those tuples $ \{\vec r_v \mid v\in V\}$ 
such that if $ v\in V _{3/2} \cup V _{1/2}$ then, $ r_{v,1}=b_v$. 
And if $ v$ is in a clique $ Q\in \mathcal Q _{\textup{fixed}}$ of color $ t$, 
then $  r _{v,t}=c _{Q} $. 

\emph{Base Case of the Recursion.} We update $ V _{3/2} \leftarrow \{v _{\ell }\}$, 
since it is the maximal element.  We update $ \mathcal Q _{\textup{fixed}}$ 
to those cliques which contain $ v _{\ell }$.  Then (\ref{e.Estimate}) 
is a consequence of (\ref{e.3halfs}).

\emph{Recursive Case.}  At this point, we have the datum 
$ V _{3/2}$, $ V _{1/2}$, and $ \mathcal Q _{\textup{fixed}}$.  
We also have datum $ \vec b\in  \{1 ,\dotsc, n \} ^ { V _{3/2}\cup V _{1/2}}$, 
and $ \vec c \in\{ 1 ,\dotsc, n \} ^{\mathcal Q _{\textup{fixed}}} $. 
Notice that this datum can completely specify the $ \mathsf r$ vectors 
associated to vertices  not in $ V _{3/2} \cup V _{1/2}$---think of a vertex that is in two cliques 
in $ \mathcal Q _{\textup{fixed}}$. 

The recursion stops if every  vertex $ v_k$  is determined by this datum. Otherwise, let  
$ k$ to be the largest integer such that $ \vec r _{v_k}$ is \emph{not} 
determined by this datum.   If \emph{no clique in $ \mathcal Q _{\textup{fixed}}$ 
contains $ v_k$} update 
\begin{equation*}
V _{3/2} \leftarrow V _{3/2} \cup \{v_k\}\,, 
\end{equation*}
and update $ \mathcal Q _{\textup{fixed}}$ to include those 
cliques which contain $ v_k$.  By application of the Littlewood Paley inequality 
and the triangle inequality, the estimate (\ref{e.Estimate}) continues to 
hold for these updated values.  

If \emph{some clique in $ \mathcal Q _{\textup{fixed}}$ 
contains $ v_k$}, then there can be exactly one clique $ Q _{v_k}$ which does, 
for otherwise $ \vec r _{v_k}$ would be completely specified by these two 
cliques.  Update 
\begin{equation*}
V _{1/2} \leftarrow V _{1/2} \cup \{v_k\}\,, 
\end{equation*}
and update $ \mathcal Q _{\textup{fixed}}$ to include  all 
cliques which contain $ v_k$.  By application of the Littlewood Paley inequality 
and the triangle inequality, the estimate (\ref{e.Estimate}) continues to 
hold for these updated values.  

\smallskip 

Once the recursion stops the inequality (\ref{e.Estimate}) holds.  But note that 
we necessarily have 
\begin{equation*}
\operatorname {Prod} (X (G; \vec b; \vec c)) ^2 \equiv 1\,, 
\end{equation*}
as all $ \mathsf r$ vectors are completely determined by $ \vec b$ and $ \vec c$. 
Therefore, we have proven (\ref{e.purpose}).

\subsubsection*{The Conclusion of the Proof.}
Since $ V _{3/2} $ and $ V _{1/2}$ are disjoint subsets of $ V$, we have proven 
the inequality 
\begin{equation}\label{e.Purpose}
\rho ^{\abs{ V}}\norm 
\operatorname {Prod} (X (G)) .p. 
\le (C \sqrt p) ^{\abs{ V}}
n ^{ \tfrac 32\abs{ V _{3/2}}+ \tfrac 12 \abs{ V _{1/2}} - \abs{ V} } 
\,. 
\end{equation}
And the remaining analysis concerns the exponent on $ n$ above, namely 
we should see that 
\begin{equation}\label{e.eta}
\abs{ V} ^{-1} \bigl[ \tfrac 32\abs{ V _{3/2}}+ \tfrac 12 \abs{ V _{1/2}} - \abs{ V} \bigr]
\le  -\eta \,,
\end{equation}
for a fixed  positive choice of $ \eta  $, and all connected graphs $ G $ 
on at least four vertices .  We would conclude that this collection of 
graphs  is in $ \operatorname {BG} (\eta )$.

It would be helpful to consider a couple of simple cases.  
Consider the graph on five vertices  
\begin{equation} \label{e.5}
\begin{array}{cccccccccc}
v_1 && v_4 && v_2 && v_5 && v_3
\\
\tfrac 32 && 0    && \tfrac 32 && 0  && \tfrac 32 
\\
  &   &   &&   &&    && 
\\
\bullet & = & \bullet  &    &   \bullet & = & \bullet  && 
\\
	&   &  \bullet & = & \bullet   &&      \bullet & = & \bullet 
\end{array}
\end{equation}
Note that we specify a particular order on the vertices  in the top row, 
and indicate the membership of each vertex in $ V _{3/2}$, $ V _{1/2}$, and in $ V_0
\coloneqq V- V _{3/2}- V _{1/2}$. 
Note that the zeros at $ v_4$ and $ v_5$ are forced.  
Consider the graph on six vertices  
\begin{equation} \label{e.6}
\begin{array}{cccccccccccc}
v_1 && v_6 && v_2 && v_5 && v_3  && v_4 
\\
\tfrac 32 && 0    && \tfrac 32 && 0  && \tfrac 32 && \tfrac 12 
\\
  &   &   &&   &&    && 
\\
\bullet & = & \bullet  &    &   \bullet & = & \bullet  &&  \bullet &=& \bullet 
\\
	&   &  \bullet & = & \bullet   &&      \bullet & = & \bullet 
\end{array}
\end{equation}
Here, there is one vertex in $ V _{1/2}$, but of course all vertices  
in $ V _{1/2}$ contribute to the Beck Gain.  But the reader should keep in mind 
that the graphs can in general have a much more complicated structure than 
these two linear examples.

\smallskip 

The extremal cases in the estimate (\ref{e.eta}) are those cases 
in which $ V _{3/2}$ is as large as possible.
 To continue, we note another formula.  
Let $ E (G)$ be the total number of edges in the graph $ G$, 
and let $ E (v)$ be the number of edges in $ G$ with one endpoint of the 
edge being $ v$.

For $ v\in V _{3/2} \cup V _{1/2}$, let $ F (v)$ be the number of edges 
which are specified upon the selection of that vertex in our recursive procedure. 
It is clear that we have $ E (v)=F (v)$ if $ v\in V _{3/2}$.  But also, 
\begin{equation*}
\sum _{v\in V _{3/2}\cup V _{1/2}} F (v)= E (G)\,. 
\end{equation*}
It follows that to maximize the cardinality of $ V _{3/2} $, 
those vertices  must be in small cliques.  There are two different classes 
of graphs which are extremal with respect to these criteria.

The 
first  extremal class consists of  graphs $ G$ with all cliques being of size $ 2$, 
and the number of cliques is $ \abs{ V}+1$, that is the graphs 
are like in (\ref{e.5}) and (\ref{e.6}). 
For such graphs, 
$ \abs{ V _{3/2}} \le \lceil \tfrac 12  \abs{ V} \rceil$,
and if the value is maximal then $ V _{1/2}$ is either $ 0$ if $ \lvert  V\rvert $ is odd, 
and $ 1$ if $ \lvert  V\rvert $ is even.  It is straight forward to see that the 
maximum of (\ref{e.eta}) occurs at $ \abs{ V}=5$, and is $ -\tfrac1 {10}$. 
Here, it is vital that we have already discussed the case of two and 
three vertices !

The second class are graphs on an even number of vertices , with 
half the vertices  in a  clique $ Q$, and each vertex $ v\in Q$ is in 
one other clique of size $ 2$.  One can depict the graph as 
\begin{equation*}
\begin{array}{ccccccccccc}
v_1  && v_2 && v_3 && v_4 && v_5 && v_6
\\
 \ast  &=& \ast &=& \ast  
\\
a && b && c && a && b && c
\end{array}
\end{equation*}
The   vertices  are written in increasing order:
$ v_1< v_2< v_3<v_4<v_5<v_6$.  Note that $ v_1,v_2,v_3$ 
form a single clique of color $ 2$.  
There are three additional cliques of size $ 2$, all of color $ 3$. They are 
 $ \{v_j,v _{j+3}\}$ 
for $ j=1,2,3$. 
  For such a graph, it is clear that $
\abs{ V _{3/2}}=\tfrac12 \abs{ V}$, and $ \abs{ V _{1/2}}=1$.\footnote{If for example 
the maximal vertex $v_6 $ where in the clique of size $ 3$, our algorithm then 
predicts a smaller estimate for the such a graph.}
The
term (\ref{e.eta}) behaves exactly like the first class of extermal  graphs 
on an even number of verticies. Our proof is complete.

\end{proof}

\chapter{Irregularities of Distributions}

\section{Discrepancy}  \label{s.discrep}

We outline the Discrepancy Theory, highlighting its relevance to the Small Ball Problem. 
In $d$ dimensions, one takes $\mathcal A_N$ to be $N$ points in the unit cube, and considers 
the function 
\begin{equation}  \label{e.discrep}
D_N(x)=\sharp \mathcal A_N \cap [\vec 0,\vec x)-N \abs{ [\vec 0,\vec x)}
\end{equation}
Here, $[\vec 0,\vec x)=\prod _{j=1}^d [0,x_j)$, that is a rectangle with antipodal corners being $\vec 0$ and $\vec x$.
We will typically suppress the dependence upon the selection of points $\mathcal A_N$. 
A set of points will be \emph{well distributed} if 
this function is small in some appropriate function space.  Thus, it of interest to understand the 
`min--max' function 
\begin{equation*}
\inf _{\mathcal A_N} \,\norm D_N.L^p([0,1]^d) .\,, \qquad  0<p\le \infty \,. 
\end{equation*}
For the purposes of this note, we will primarily be concerned with 
lower bounds for this quantity, with $ 1\le p\le \infty $. 
Dimension will be held fixed, with 
$ N$ large.  Many variants of this question are interesting;  interested readers is 
encouraged to consult one of the excellent references in this area.

It  turns out that relevant norms of this function must tend to infinity, in dimensions $2$ and higher. 
Using the basic facts of the next section, we can prove the Theorem below, 
which concatenates results of Roth \cite{MR0066435} in the case of $ p=2$. Indeed 
the proof we give below is the `hyperbolic orthogonal function' method he initiated; 
and Schmidt \cite{MR0491574} for other values of $ p$.  
 The end point estimate  below is a consequence of the method, and don't seem 
 to be as well known.

\begin{theorem}\label{t.DP} 
For any collection of points $ \mathcal A_N\subset [0,1] ^{d}$, we have the 
estimates 
\begin{equation}\label{e.DP}
\norm D_N. p. \gtrsim (\log N) ^{(d-1)/2}
\end{equation}
More particularly, we have the endpoint estimate
\begin{equation}\label{e.DPend}
\norm D_N. L (\log L )  ^{(d-1)/2}. 
\gtrsim (\log N) ^{(d-1)/2}
\end{equation}
\end{theorem}

\begin{proof}
As is usual, the proof is by duality, following Roth \cite{MR0066435}, and we 
use the Haar function approach of Schmidt \cite{MR0319933}.

We stick to the hyperbolic setting, with the rationale that extremal point distributions, 
whatever they might be, must have about one point in any rectangle of volume 
about $ 2 ^{-n}$. 

For each $ \vec r\in \mathbb H_{n}$
 construct the $ \mathsf r$ function $ f _{\vec r}$ as in Proposition~\ref{p.rvec}, and set 
 \begin{equation*}
F \coloneqq \sum _{\vec r\in \mathbb H_{n}} f _{\vec r}\,. 
\end{equation*}
By construction we have
\begin{equation*}
n ^{d-1} \lesssim \ip D_N, F, \le \norm D_N .2. \norm F .2. \le \norm D_N.2.\, n ^{(d-1)/2}\,. 
\end{equation*}
This prove (\ref{e.DP}) in the case of $ p=2$, and by extension to all $ p\ge 2$. 
To finish the proof,
recall that $ L (\log L )  ^{(d-1)/2}$ and $ \operatorname {exp} (L ^{2/(d-1)})$ 
are dual spaces, see \S~\ref{s.orlicz}.  Thus, we 
we should observe that 
\begin{equation*}
\norm F. \operatorname {exp} (L ^{2/(d-1)}) . \lesssim n ^{(d-1)/2} \,.
\end{equation*}

But, the square function of $ F$  
\begin{equation*}
\operatorname S (F) \coloneqq 
\Bigl[ \sum _{\vec r\in \mathbb H_{n} } \abs{ f _{\vec r}} ^2  \Bigr] ^{1/2} 
\lesssim  n ^{(d-1)/2}\,. 
\end{equation*}
The last estimate is an $ L ^{\infty }$ estimate. 
Therefore, by Theorem~\ref{t.distributional}, we conclude that 
$\norm F. \operatorname {exp} (L ^{2/(d-1)}) .  \lesssim n ^{(d-1)/2} $. 
This implies the $ L (\log L )  ^{(d-1)/2}$
 endpoint estimate for $ D_N$ in (\ref{e.DPend}). 
 
\end{proof}

While this last Theorem is quite adequate for $ L ^{p}$, the endpoint cases 
of $ L ^{\infty }$ and $ L ^{1}$ are not amenable to the same techniques, 
and the relevant fact is that the $ L ^{\infty }$ bound should be larger. 
In dimension $2$, the end point estimates are known.  At $ L ^{\infty }$, 
it is the Theorem of Schmidt \cite{MR0319933}. 

\begin{schmidt} 
We have the estimates below, valid for all collections $\mathcal A_N\subset [0,1]^2 $:  
\begin{align}  \label{e.schmidt}
\norm D_N .\infty. {}\gtrsim{} \log N . 
\end{align}
\end{schmidt} 
We shall see that this is a rather precise analog of Talagrand's  theorem; the proof we give will share a great deal 
of similarity with the proof of Temlyakov we have described in \S~\ref{s.talagrand}.   

Let us comment that there is an interpolant between the result of Schmidt and the $ L ^{p}$ 
results, provided one uses the scale of exponential Orlicz classes. 

\begin{theorem}\label{t.beyondSchmidt}  
We have the estimates below, valid for all collections $\mathcal A_N\subset [0,1]^2 $:  
\begin{align}  \label{e.schmidt}
\norm D_N . \operatorname {exp} (L ^{p}). {}\gtrsim{} (\log N ) ^{1-1/p} \, , \qquad 2<p<\infty \,. 
\end{align}
\end{theorem}

In dimensions $3$ and higher,  there is the following   improvement on 
J.~Beck's result \cite{MR1032337}, due to Lacey and Bilyk \cite{bl}. 

\begin{theorem}\label{t.beck} 
There is a choice of $0<\eta <\tfrac12 $  for which the following estimate 
holds for all collections $\mathcal A_N\subset [0,1]^3 $:  
\begin{align}  \label{e.beck}
\norm D_N .\infty. {}\gtrsim{} (\log N) ^{1+\eta }\,. 
\end{align}
\end{theorem}

Beck's result is as above, with $ (\log N) ^{\eta }$ replaced by  a 
doubly logarithmic term. 
There is no further result known about the Small Ball Problem, nor the Discrepancy 
Function in higher dimensions.\footnote{The student of the literature will find an article 
published some years ago that claims an extension of Beck's result to higher dimensions. 
While this paper can serve as a useful summary of Beck's argument, 
an early critical Lemma in that paper is in error; a technique to repair the 
error is unknown to me.}

Hal{\'a}sz established the $ L ^{1}$ endpoint estimate for the Discrepancy function 
in two dimensions.  Namely

\begin{halasz} \label{t.halasz}
For any collection of points  $ \mathcal A _N \subset 
[0,1] ^{2}$ of cardinality $ N$ we have 
\begin{equation}\label{e.halasz}
\norm D_N .1. \gtrsim \sqrt {\log N}\,. 
\end{equation}
\end{halasz}

While the $ L ^{\infty }$ case is in close analogy to the Small Ball Conjecture, 
this analogy breaks down in this case.   
We will give Hal{\'a}sz' proof of this result, as well as a new one, which is 
again a duality method, but the construction of the dual function is 
\emph{not} by way of a Riesz product.  See \S~\ref{s.halasz}.

In the reverse direction, concerning point distributions with small Discrepancy 
function, the following is known. 

\begin{theorem}\label{t.<}  In dimension $ d$, there are point distributions $ \mathcal A_N$
with 
\begin{equation*}
\norm D_N . p. \lesssim (\log N) ^{(d-1)/2}\,,  \qquad 0<p<\infty \,. 
\end{equation*}
\end{theorem}

These constructions are delicate, and the product of significant effort over a 
period of decades.  See especially Davenport \cite{MR0082531}, Roth \cites{MR553291
,
MR598865}, and Chen \cite{MR610701}. 
These earlier constructions were random in nature; recently Chen and 
 Skriganov \cites{MR1896098,MR1805869} found subtle deterministic constructions.

On the other hand, Schmidt's result is sharp, for Halton \cite{MR0121961} has constructed 
point sets with Discrepancy function of $ L ^{\infty }$ norm that matches 
his lower bound. 

\begin{halton} For dimension $ d\ge 2$ there are point sets $ \mathcal A_N$ with 
\begin{equation*}
\norm D_N. \infty . \gtrsim (\log N) ^{d-1}
\end{equation*}
\end{halton}

\section{Conjectures for Discrepancy} 

\subsection*{The $ L ^{\infty }$ Conjectures}
In light of the  close connection between the proof of the lower 
bounds in the $ L ^{\infty }$ case  
and the Small Ball Conjecture,  one suspects that 
an extra square root of $n {}\simeq{}\log N$ 
is all that should  be obtainable at the end point estimate at $ L ^{\infty }$ for 
the discrepancy function.

\begin{hyperbolic}
For all choices of $N$ points $\mathcal A_N$ we have 
\begin{equation}  \label{e.discreplower}
\norm D_N.\infty. {}\gtrsim{} (\log N) ^{d/2}. 
\end{equation}
\end{hyperbolic}

What should be clear, in light of the sharpness of the Small Ball Conjecture, 
is that those who hold the conviction that this last conjecture falls short of the 
truth will necessarily seek a proof other than the hyperbolic one.

\begin{Sharp} 
We have the estimate 
\begin{equation}  \label{e.discrepconjecture}
\inf _{\mathcal A_N}\; \norm D_N.\infty.  {}\lesssim{}  (\log N) ^{d/2}.
\end{equation}
\end{Sharp}

In this paper, we emphasize the similarity in proof techniques in the 
Small Ball Problem and the Discrepancy problems.  It would be of interest to 
establish some formal connection between these two problems.

The reader can consult the survey article by Temlyakov \cite{MR1984119} 
for a discussion of the connection between the Discrepancy function in 
$ L ^{\infty }$ and cubature formulas.

One suspects that Theorem~\ref{t.beyondSchmidt} is sharp.  (Compare to \cite{2000b:60195}.)

\begin{conjecture} In dimension $ 2$, one has 
\begin{equation*}
\min _{\mathcal A_N} \norm D_N . \operatorname {exp} (L ^{\alpha }). 
\simeq 
(\log N) ^{1-1/\alpha }\,, \qquad 2\le \alpha <\infty \,. 
\end{equation*}
\end{conjecture}

\subsection*{The $ L ^{1}$ Conjecture}

The other outstanding conjecture concerns the  $ L ^{1}$ norm endpoint.

\begin{ellOne}
In any Dimension $ d$ one has the estimate 
\begin{equation*}
\norm D_N. 1. \gtrsim (\log N) ^{(d-1)/2}\,. 
\end{equation*}
\end{ellOne}

It appears that any improvement in the estimate (\ref{e.DPend}), by 
e.\thinspace g.~replacing the logarithmic Orlicz space by one closer to $ L ^{1}$, 
will generate an interesting new proof technique.

\subsection*{The $ L ^{p}$ Conjecture, for $ 0<p<1$}

One can ask about the size of the Discrepancy Function in $ L ^{p}$, for $ 0<p<1$. 
The absence of duality methods has prevented any progress towards this conjecture. 

\begin{conjecture}  We have the estimate below, for all $ 0<p<1$. 
\begin{equation*}
\norm D_N.p. \gtrsim (\log N) ^{(d-1)/2}\,. 
\end{equation*}
\end{conjecture}

Here, we indicate a result in this direction. 

\begin{theorem}\label{t.<1}  For $ 0<p<1$, and dimension $ d\ge 2$ we have the estimate 
\begin{equation*}
\norm \operatorname M D_N .p. \gtrsim (\log N) ^{(d-1)/2}\,. 
\end{equation*}
 Here, $ \operatorname M$ denotes the strong maximal function in $ d$ dimensions, 
 thus 
 \begin{equation*}
\operatorname M f (x)=\sup _{\textup{$ R$ dyadic}} \mathbf 1 _{_R } (x) \mathbb E (f\,|\, R)\,. 
\end{equation*}
\end{theorem}

\begin{proof}
We are uncertain as to how interesting this is, so our proof is somewhat abbreviated. 
The only real observation to make is that the theory of multi-parameter
Hardy space is relevant.  See \cites{cf1,MR658542}.  In particular, letting $ H ^{p}$ 
denote Hardy space, one has 
\begin{equation*}
\norm f. H ^{p}. \simeq \norm \operatorname  M f.p. \simeq \norm \operatorname S (f).p.\,, 
\quad 0<p\le 1\,. 
\end{equation*}

We apply this to $ D_N$.  Let $ \mathcal G$ be the class of good rectangles, as defined in 
Proposition~\ref{p.rvec}.  We then have 
\begin{equation*}
\norm \operatorname M D_N.p. \simeq \norm \operatorname S (D_N).p. 
\gtrsim \NOrm \Bigl[ \sum _{R\in \mathcal G} \mathbf 1 _{R} \Bigr] ^{1/2} .p.
\end{equation*}
It is an elementary exercise to see that the last term is $ \gtrsim (\log N) ^{(d-1)/2}$. 
\end{proof}

\section{Elementary Propositions} 

 Throughout, we will specify 
$n$ by $2N\le{} 2^n<4N$, so that $ n \simeq \log N$. 
The value of $ n$ plays the same role in this section as it does in our 
discussion of the Small Ball Conjecture.  In this section, we use the 
notation and definitions of \S~\ref{s.definition}.

Recall that $f$ an $\mathsf r$ function if it is equal to 
\begin{equation*}
f=\sum _{R\in \mathcal R _{\vec r}} \varepsilon _{R} h_R
\end{equation*}
where $\varepsilon _{R}\in \{-1,0,1\}$.   Recall that $ \mathcal R _{\vec r}$ 
consists of all dyadic rectangles $ R$ with $ \abs{ R_j} = 2 ^{-r_j}$ for 
all coordinates $ j$. 

\begin{proposition}\label{p.rvec}  For each $\vec r\in \mathbb H_{n}$, there is an $\mathsf r$ function $f _{\vec r} $ with 
\begin{equation*}
\ip D_N, f _{\vec r},\ge{}c_d. 
\end{equation*}
Here $ c_d$ is a dimensional constant. 
\end{proposition}

\begin{proof}
There is a very elementary one dimensional fact: For all dyadic intervals $I$, 
\begin{equation} \label{e.veryelementary}
\mathbb E x \mathbf 1 _{I} (x)  h _{I}(x)=\tfrac 14 \abs{ I}^2.
\end{equation}
This immediately implies that in any dimension 
\begin{equation*}
\mathbb E  h_R(\vec x) \abs{ [0,\vec x)}= 4 ^{-d}\abs{ R}^2.
\end{equation*}

We shall rely upon the construction of the this function $f _{\vec r}$ below.  Recall that $\mathcal A_N$, the distribution of 
$N$ points in the unit cube, is fixed.  Call a cube $R\in \mathcal R _{\vec r}$ \emph{good} if $R$ does {\bf not}
 intersect $\mathcal A_N$, otherwise call it \emph{bad}. 
 Set 
 \begin{equation}  \label{e.f_r}
 f _{\vec r} {} \coloneqq {}\sum _{\substack{R\in \mathcal R _{\vec r}\\ \text {$R$ is good} } }h_R 
 +
 \sum _{\substack{R\in \mathcal R _{\vec r}\\ \text {$R$ is bad} } } 
 \operatorname {sgn} (\ip D_N, h_R,) h_R
 \,. 
 \end{equation}

 Each bad rectangle contains at least one point in $\mathcal A_N$, and  $2^n\ge2N$, so 
 there are at least $N$ good rectangles.  Moreover, since the counting function 
$\sharp \mathcal A_N\cap [0,\vec x)$ is constant over each good rectangle, we have 
\begin{equation*}
\ip D_N, h_R,=N\prod _{j=1}^d \ip x_j, h _{R_j},=N2 ^{-2n-2d}
\gtrsim 2 ^{-n}
\end{equation*}
Hence, we can estimate 
\begin{align*}
\ip D_N, f _{\vec r}, \ge 
\sum _{\substack{R\in \mathcal R _{\vec r}\\ \text {$R$ is good} } }
\ip D_N, h_R, \gtrsim 2 ^{-n} \sharp \{
R\in \mathcal R _{\vec r}\mid \text {$R$ is good} \} 
\gtrsim 1\,. 
\end{align*}
And so our proof is complete. 
\end{proof}

Another proposition of a similar flavor is this. 

\begin{proposition}\label{p.>n} Let $ f _{\vec s}$ be any $ \mathsf r$ function 
with $ \abs{ \vec s}>n$.  We have 
\begin{equation*}
\abs{ \ip D_N, f _{\vec s}, } \lesssim N {2 ^{-\abs{ \vec s}}}\,. 
\end{equation*}

\end{proposition}

\begin{proof}
This is a brute force proof.  Consider the linear part of the Discrepancy function. 
By (\ref{e.veryelementary}), we have 
\begin{equation*}
\abs{ \ip N \prod _{j=1} ^{d} x_j, f _{\vec s}, }\lesssim N 2 ^{\abs{ \vec s}} \,,
\end{equation*}
as claimed. 

Consider the part of the Discrepancy function that arises from the point set.  
Observe that for any point $\vec x_0$ in the point set, we have 
\begin{equation*}
\abs{ \ip \mathbf 1 _{[\vec 0, \vec x)} (\vec x_0), f _{\vec s}, } \lesssim 2 ^{- \abs{ \vec
s}}\,.
\end{equation*}
Indeed, of the different Haar functions that contribute to $ f _{\vec s}$, there 
is at most one with non zero inner product with the function 
$\mathbf 1 _{[\vec 0, \vec x)} (\vec x_0) $ as a function of $ \vec x$.  It could only be the 
one rectangle which contains $ x_0$ in its interior. Thus the inequality above follows. 
Summing it over the $ N$ points in the point set finish the proof of the Proposition. 

\end{proof}

A final, general proposition is relevant. 

\begin{proposition}\label{p.Not2Many} Fix a collection of $ \mathsf r$ functions 
$\{ f _{\vec r} \mid \vec r \in \mathbb H_{n} \}$.  Fix $ \vec s $ with $ \abs{ \vec s}>n$, 
and let $ 3\le k\le \abs{ s}-n+1$. 
Let $ \operatorname {Count} (\vec s ; k) $ be 
{the number of ways to  choose strongly distinct 
$ r_1 ,\dotsc, r_k\in \mathbb H _n$  so that $ \prod _{j=1} ^{k} f _{\vec r}$
is an $ \vec s$ function.}
We have 
\begin{equation} \label{e.Not2Many}
\operatorname {Count} (\vec s ; k) 
\lesssim  (\abs{ \vec s}-n) ^2 \cdot  k ^{3} \cdot  { (\abs{ \vec s}-n) ^{(d-1)} \choose k-3} \,. 
\end{equation}
For $ k=2$ we have 
\begin{equation}\label{e.Not2Many2} 
\operatorname {Count} (\vec s ; 2) 
\lesssim  \abs{ \vec s}-n\,. 
\end{equation}
\end{proposition}

\begin{proof}
This estimate is only of interest for $ \abs{ \vec s}< d n$, 
and is very crude. 
Fix $ \vec s$.   We want to choose strongly distinct 
$ r_1 ,\dotsc, r_k\in \mathbb H _n$ so that for all coordinate $ 1\le t \le d$ we have 
\begin{equation*}
\max \{r _{1,t} ,\dotsc, r _{k,t}\}=s_k\,.
\end{equation*}
(Of course if the $ \vec r _k$ are not strongly distinct, the product need not be 
a $ \mathsf r$ function.)
Observe that 
for given $ \vec s$, there are at most $ \lesssim (\abs{ \vec s}-n) ^{d-1}$ 
vectors $ \vec r\in 
\mathbb H _n$ with $ r _{t}\le s_t$ for all coordinates $ t$. 

Since the product is to be an $ \mathsf r$ function with parameter $ \vec s$, we 
must have either two or three of the chosen $ \vec r$ functions whose parameters 
are maximal, and equal to $ \vec s$.  There are at most $ k ^{3}$ ways to select 
these $ \mathsf r$ functions among the $ k$ terms were are forming the product over. 
And having selected them, there are at most $ (\abs{ s}-k) ^2 $ ways to select 
these $ \mathsf r$ functions.  The remaining $ k-3$ $ \mathsf r$ functions 
can be selected freely.  This gives (\ref{e.Not2Many}).  

The second estimate (\ref{e.Not2Many2}) is easier.

\end{proof}

\section{Proof of Schmidt's Theorem}  \label{s.schmidt}
We prove the Theorem of Schmidt; this section should be compared to \S~\ref{s.talagrand}. 
With the $\mathsf r$ functions as constructed in the the proof 
of Proposition~\ref{p.rvec}, we set 
\begin{equation*}
\Psi {} \coloneqq {} \prod _{\vec r\in \mathcal P_n}(1+ \alpha f _{\vec r} ) . 
\end{equation*}
Here, $0<\alpha<\frac12$, and to be specific, we can choose $\alpha= 2 ^{-6 }$. 
Clearly, this is a non negative function, with $\int \Psi\; dx=1$.  And so we should argue that 
\begin{equation*}
\ip D_N, \Psi, {}\gtrsim{} n. 
\end{equation*}
  
Write the function $\Psi$ as 
\begin{align*}
\Psi&{}=\sum _{k=0} ^{n}  \psi_k
\\ 
\psi_k&{}={} \sum _{\substack{ W\subset \mathcal P_n \\ \sharp W=k}}  \alpha^{\sharp W} \prod _{\vec r\in W} f _{\vec r}
\end{align*}
where we understand that  $ \psi_0=\mathbf 1 _{[0,1]^2}$. 
 
Clearly, $\ip D_N, \psi_0,=0$.  By Proposition~\ref{p.rvec}, we have 
\begin{equation}  \label{e.primary}
\sum _{\vec r\in \mathcal P_n} \ip D_N,\alpha f_{\vec r},\ge{}\alpha n
\end{equation}
For this, recall that we are specializing to the case of dimension $2$.
 
We provide an upper bound on the remaining inner products $\ip D_N,\psi_k,$ for $k\ge2$.\footnote{Note that 
in the small ball problem, this set is not needed!} 
For a subset $W\subset \mathcal P_n$ of cardinality at least $2$.  Then, the product 
\begin{equation*}
\prod _{\vec r\in W} f _{\vec r}
\end{equation*}
is again a sum of Haar functions, by the Product Rule! See Theorem~\ref{p.product}. 
By Proposition~\ref{p.>n}, 
\begin{equation*}
\Abs{ \Ip D_N,  \prod _{\vec r\in W} f _{\vec r} , } {}\lesssim{} N 2 ^{-\abs{ \vec w}}.  
\end{equation*}

Now, for a fixed $k$ and $\vec w$ with $n+k\le{} \abs{ \vec w }\le2n$, we count the number of distinct 
ways of choosing $W$ so that $ \prod _{\vec r\in W} f _{\vec r}$ is a $\vec w$ function. 
The first coordinates of the vectors $\vec r$ must be $k$ distinct integers  in the range 
\begin{equation*}
n-w_2\le{} r_1\le{} w_1.
\end{equation*}
 Moreover, there must be choices of $\vec r\in W$  whose first coordinates are  equal to either endpoint.
 There are clearly at most 
\begin{equation} \label{e.2Not2}
{ \abs{ \vec w}-n-2 \choose k-2 } 
\end{equation}
choices of $W$.

For an integer $n\le{}\omega\le2n$, there are at most $2n$ vectors $\vec w$ with $\abs{ \vec w}=\omega$.  Therefore, 
\begin{align*}
\abs{\ip D_N, \psi_k,} & {}\le{} 2n\alpha^k N \,  \sum _{\omega=n+k} ^ {2n } { \omega-n-2 \choose k-2 } \; 2 ^{- \omega }
\\&{}=  n\alpha^k N 2 ^{-n-k+1} \sum _{\omega=0}^{n-k}  {\omega+k-2\choose k-2 } 2 ^{-\omega}
\end{align*}
This must be summed over $2\le{} k\le{}n$.  This sum is treated by two 
changes of variables. (One is $v=\omega+k$.)
\begin{align*}        
\sum _{k=2}^n  \sum _{\omega=0}^{n-k}  n \alpha^k N 2 ^{-n-k+1}    {\omega+k-2\choose k-2 } 2 ^{-\omega}
&{}={} \alpha^2 2 ^{-n-1}N  \sum _{k=0} ^{n-1} \sum _{\omega=0} ^{n-k} \alpha ^{k} 2 ^{-k-\omega} { \omega+k\choose k} 
\\&{}\le{}n \alpha^2 2 ^{-n-1}N  \sum _{v=0}^n \sum _{k=0} ^{v} \alpha^k 2 ^{-v} {v \choose k} 
\\&{}\le{} n\alpha^2 2 ^{-n-1}N \sum _{v=0} ^{n} 2 ^{-v}(1+\alpha)^v
\\&{}\le{} 4n\alpha^2 
\end{align*}
For $\alpha$ sufficiently small, we see that this estimate is much smaller than the lower bound 
in (\ref{e.primary}), so that our proof is complete.

\bigskip 

The proof of Theorem~\ref{t.beyondSchmidt} is a simple corollary to the proof above. 
Since $ \norm \Psi . \infty . \le 2 ^{n}$, it is clear that we have 
\begin{equation*}
\norm \Psi . L ^{1} (\log L) ^{\alpha }. \lesssim (\log N) ^{\alpha }\,. 
\end{equation*}
Therefore, we can estimate for $ 2<p<\infty $
\begin{equation*}
 {\log N} \lesssim \ip D_N, \Psi , \lesssim \norm D_N . \operatorname {exp} (L ^{p}). 
\cdot 
\norm \Psi . L ^{1} (\log L) ^{1/p}. 
\lesssim 
(\log N) ^{1/p} \norm D_N . \operatorname {exp} (L ^{p}). \,. 
\end{equation*}

\section{Proof of Theorem~\ref{t.beck}}\label{s.blD}

We rely upon \S~\ref{s.Short}. 
We see that for $ \Psi ^{\textup{sd}}$ as defined in (\ref{e.zCsd}), that we have 
$ \norm \Psi ^{\textup{sd}} .1. \lesssim 1$. Moreover, we have 
\begin{equation} \label{e.zCmain}
\ip D_N, \Psi ^{\textup{sd}} _{1}, \gtrsim a q ^{b} n 
\simeq  a n ^{1+\epsilon /4}
\end{equation}
Here, $ q$ is defined as in (\ref{e.q}), and $ 0<a<1$ is a small constant.
Again, $ q ^{b} \simeq n ^{\epsilon /4}$ is the `gain over the trivial estimate.'

\smallskip 

Consider the terms arising from $ \Psi ^{\textup{sd}} _{k}$.  These are products of 
strongly distinct $ \mathsf r$ vectors.  Hence, 
we combine the estimates from Proposition~\ref{p.>n} and Proposition~\ref{p.Not2Many} 
as follows. 
For $ k=2 $ we have 
\begin{align*}
\abs{ \ip D_N, \Psi ^{\textup{sd}} _{k},} 
&\le \sum _{h=1} ^{3n} \sum _{\vec s\,:\, \abs{ \vec s}=n+h} 
\bigl[\frac {q ^{b}} {n} \bigr] ^2 \cdot N 2 ^{-n-h} \cdot 
\operatorname {Count} (\vec s; 2)
\\
& \lesssim 
\bigl[\frac {q ^{b}} {n} \bigr] ^2
\sum _{h=1} ^{3n}  (n+h) ^2 h \cdot  2 ^{-h}
 \\
& \lesssim q ^{2b}= n ^{\epsilon /2}\,. 
\end{align*}
This is much smaller than the main term (\ref{e.zCmain}). 

\smallskip

We treat the terms arising from $ \Psi _k$ for $ k\ge 3$ as follows. 
\begin{align*}
\sum _{k=3} ^{q}\abs{ \ip D_N, \Psi ^{\textup{sd}} _{k},} 
&\le \sum _{k=2} ^{q} \sum _{h=k} ^{ 2n } 
\sum _{\abs{\vec s}= n+h} N 2 ^{-n-h  } \bigl[ \frac {a q ^{b}} { n }
\bigr] ^{k}
 \cdot \operatorname {Count} (\vec s ; k) 
 \\
 & \lesssim
 \sum _{k=3} ^{q} \sum _{h=k} ^{ 2n }  h ^2 2 ^{-h} k ^{3}
  \bigl[ \frac {a q ^{b}} { n  }
\bigr] ^{k}  { h ^2 \choose k-3} 
 \\
 & \le 
 q ^{3}  \bigl[ \frac {a q ^{b}} { n }
\bigr] ^{3} \sum _{h=3} ^{2n}  h ^2 2 ^{-h}\sum _{j=0} ^{q} 
 \bigl[ \frac {a q ^{b}} { n}
\bigr] ^{j}  { h ^2 \choose j} 
\end{align*}
We have crudely estimated a term or two, and reversed the order of summation. 
Observe that $ q= n ^{\epsilon }$ is much smaller than $ n$, so that we can estimate 
\begin{equation*}
\sum _{j=0} ^{q} 
 \bigl[ \frac {a q ^{b}} { n}
\bigr] ^{j}  { h ^2 \choose j} 
\le 
\sum _{j=0} ^{h ^2 } 
 \bigl[ \frac {a q ^{b}} { n}\bigr] ^{j} \cdot \bigl[1- \frac {a q ^{b}} { n}\bigr] ^{h ^2 -j}
  { h ^2 \choose j} \le 1\,. 
\end{equation*}
It follows that 
\begin{align*}
\sum _{k=3} ^{q}\abs{ \ip D_N, \Psi ^{\textup{sd}} _{k},} 
& \lesssim 
q ^{3}  \bigl[ \frac {a q ^{b}} { n  }
\bigr] ^{3} \sum _{h=3} ^{2n}  h ^2 2 ^{-h} \lesssim  q ^{6} \cdot n ^{-3}
\end{align*}
which is again much smaller than the main term (\ref{e.zCmain}). Our proof is complete.

\section{The $ L ^{1}$ bound in dimension $ 2$} \label{s.halasz}

We will indicate two proofs of Hal{\'a}sz' Theorem~\ref{t.halasz}. 
The first is the proof of Hal{\'a}sz. 
Let $ f _{\vec r} $ be the $ \mathsf r$ functions has in Proposition~\ref{p.rvec}. 
Consider the Riesz product 
\begin{equation*}
\Psi \coloneqq \prod _{t=1} ^{n}(1+i \frac {a} {\sqrt n} f _{\vec r})\,.
\end{equation*}
Here, $ 0<a<1$ is a small constant to be chosen.  Because of the imposition of the 
imaginary number, it is is evident that this $ \Psi $ is a bounded complex valued 
function.   But one can argue that 
\begin{equation*}
\ip D_N , \operatorname {Im} (\Psi ), \gtrsim \sqrt n \,. 
\end{equation*}
much as the lines of the argument used to prove Schmidt's theorem. We omit the details.  

\bigskip 

The second proof, is as far as the author knows, is new; as with Ha{\'a}sz' proof, 
it does not admit a straight forward extension to higher dimensions.   We offer 
it as a technically interesting object, as the function we use 
 is \emph{not} a Riesz product, rather it is 
\begin{equation}\label{e.sine}
\Phi \coloneqq \sin \Bigl( \frac  a {\sqrt n} \sum _{\abs{ \vec r}=n} f _{\vec r}  \Bigr)\,.
\end{equation}
As usual, $ 0<a<1$ is a sufficiently small constant. 
And we argue that $ \ip D_N, \Phi , \gtrsim \sqrt n$. 

Recall that the argument of the sine function above has $ \operatorname {exp} (L ^2 )$ 
norm bound independent of $ n$.  Thus, as one may directly check, the Taylor expansion 
of $ \Phi $ is convergent in all $ L ^{p}$.  That is, we may expand 
\begin{equation} \label{e.PhiExpanded}
\Phi = \sum _{k=0} ^{\infty } \frac { (-1) ^{2k-1}} {(2k+1)! n ^{(2k+1)/2}} 
\Bigl(  \sum _{\abs{ \vec r}=n} f _{\vec r}   \Bigr) ^{2k+1}\,. 
\end{equation}
and the sum is convergent in all $ L ^{p}$, $ 1<p<\infty $.  A remarkable fact 
is that this infinite expansion is in fact a finite sum.  To see this, let us observe 
the odd powers above have a simple closed form. 

\begin{lemma}\label{l.kthpower} For  integers $k $
\begin{align*} 
n ^{-(2k+1)/2}
\Bigl(  \sum _{\abs{ \vec r}=n} f _{\vec r}   \Bigr) ^{2k+1}
&= \sum _{\substack{v=1\\ \textup{$ v $ odd } }} ^{\min (n, 2k+1)} 
\frac {(2k+1)!} { 2 ^{v}}  n ^{-v}  G_v
\\
\textup{where} \quad 
G_v & \coloneqq   \sum _{ \vec r_1 ,\dotsc, \vec r_v \ \textup{distinct}} 
\prod _{w=1} ^{v} f _{\vec r_w}\,.
\end{align*}
The last sum is over all distinct $ v$ tuples of 
$ \mathsf  r$ vectors with $ \abs{ \vec r}=n$.  
\end{lemma}

\begin{proof}
Only odd products of $ \vec r$ functions can occur in the expanded product. 
Fix $ v$ odd, and distinct $ \mathsf r $ vectors 
$ \vec r_1 ,\dotsc, \vec r_v $.  It suffices to count the number of ways this 
product can arise from the expanded product.  
But this is  
\begin{equation*}
{ 2k+1 \choose v } \cdot v! \cdot  \frac {(2k+1-v)} {2 ^{v}} n ^{(2k+1-v)/2}
\end{equation*}
Indeed from the terms 
\begin{equation*}
\Bigl(  \sum _{\abs{ \vec r}=n} f _{\vec r}   \Bigr) ^{2k+1}
\end{equation*}
we choose $ v $ terms from which we take one of the pre-specified $\mathsf r $ 
functions $ f _{\vec r_1} ,\dotsc, f _{\vec r_v}$.  These products 
can be specified in one of $ v!$ ways.  

In the remaining $ 2k+1-v$ terms, we divide them into groups of two.  And select 
one of $ n$ $ \mathsf r$ functions for each pair.  This proves the Lemma. 
\end{proof}

Expanding the Taylor series  we see that 
\begin{align} \nonumber 
\Phi & = 
\sum _{ \textup{ $ k$ odd} } 
\frac { (-1) ^{(k+1)/2}} {k!}  n ^{-k/2} 2 ^{-k}
\Bigl[ \sum _{\vec r\in \mathbb H _n} f _{\vec r} \Bigr] ^{k} 
\\ \nonumber 
& = 
\sum _{ \textup{ $ k$ odd} } 
 { (-1) ^{(k+1)/2}}  2 ^{-3k/2} \sum _{\substack{ v=1\\ \textup{ $ v$ odd} }} ^{n} 
 2 ^{v/2}n ^{-v/2} G_v 
\\
& =  \label{e.finalexpand} c
\sum _{\substack{ v=1\\ \textup{ $ v$ odd} }} ^{n}  (-1) ^{(v+1)/2}  2 ^{-v}
n ^{-v/2} G_v \,.  
\end{align}
Here, $ c=(1-4 ^{-3/2}) ^{-1} $.

We turn our attention to the terms in (\ref{e.finalexpand}). 
Now, by construction, we have 
\begin{equation*}
\ip D_N , n ^{-1/2} G_1 , \gtrsim  n ^{-1/2} \sum _{\vec r\in \mathbb H _n} 
\ip D_N, f _{\vec r}, \gtrsim n ^{1/2} \simeq \sqrt {\log N}\,. 
\end{equation*}

As for the terms $ 3 \le v \le n$, note that by Proposition~\ref{p.>n}, 
(\ref{e.2Not2}) and the definition of $ G_v$, we have 
\begin{align*}
\abs{ \ip D_N, G_v , } 
& \lesssim  N \sum _{s=n+v-1} ^{2n} 2 ^{-s} {s-n-1 \choose v-2}\,. 
\end{align*}
And so we  estimate as follows.  Here is convenient that the 
sum is only over odd $v \ge 3 $. 
\begin{align*}
\sum _{\substack{ v=3\\ \textup{$ v$ odd } }} ^{n} 
 2 ^{-v} n ^{-v/2} \abs{ \ip D_N, G_v , } 
 & 
 \lesssim 
N  \sum _{\substack{ v=3\\ \textup{$ v$ odd } }} ^{n} 
\sum _{s=n+v-1} ^{2n} 2 ^{-s-v} n ^{-v/2} {s-n-1 \choose v-2}
\\
& \lesssim N  n ^{-1} 
\sum _{s = n+3} ^{2n} \sum _{v=0} ^{ s-n-1 }  2 ^{-s -v} 
n ^{-v/2} {s-n-1 \choose v} 
\\
& \lesssim  N n ^{-1} \sum _{s=n} ^{2n }  2 ^{-s}(1+1/\sqrt n) ^{s-n+1 } 
\\
& \lesssim n ^{-1}\,. 
\end{align*}
Since this estimate tends to zero with $ n$, this proves our Theorem for sufficiently 
large $ N$.

\chapter{Some Aspects of Harmonic Analysis} \label{background}

\section{Exponential Orlicz Classes} \label{s.orlicz}

Let $ \psi \,:\, \mathbb R \longrightarrow \mathbb R $ be a symmetric 
convex function with $ \psi (x)=0$ iff $ x=0$.  Define the Orlicz norm 
\begin{equation}\label{e.Orlicz}
\norm f . \psi . \coloneqq 
\inf \{C>0 \mid \mathbb E \psi (f/C)\le 1\}\,.
\end{equation}
We take the infimum of the empty set to be $ +\infty $, and denote by $ L ^{\psi  }$ 
to be the collection of functions for which $ \norm f . \psi .<\infty $.  

It is straight forward to see that $ \norm \cdot .\psi . $ is in fact a norm, 
with the triangle inequality following from Jensen's inequality.  If $ \psi (x)= x ^{p}$, 
then $ \norm \cdot .\psi .$ is the usual $ L ^{p}$ norm.  

We are especially interested in the class of $ \psi $ given by 
\begin{equation*}
\psi _{\alpha } (x)= \operatorname e ^{\abs{ x} ^{\alpha }}\, , \qquad \abs{ x} \gtrsim 1 \,.
\end{equation*}
Here, we insist upon equality for $ \abs{ x} $ sufficiently large, depending upon $ x$. 
We will write $ L ^{\psi _\alpha }= \operatorname {exp} (L ^{\alpha })$.  These 
are the exponential Orlicz classes.    

Especially important is the the case of $ \alpha =2$, which is the class 
$ \operatorname {exp} (L ^{2})$, of exponentially square integrable functions, 
of which the Gaussian random variables are a canonical example.  
A function $ f\in \operatorname {exp} (L ^{2})$ is said to be \emph{sub-gaussian.}

Using Stirling's formula, and the Taylor expansion for $ \operatorname e ^{x}$, one 
can check that 

\begin{proposition}\label{p.exp}  We have the equivalence of norms 
\begin{align*}
\norm f. \operatorname {exp} (L ^{\alpha }).  &\simeq 
\sup _{p\ge1} p ^{-1/\alpha } \norm f.p. 
\\
& \simeq 
\sup _{\lambda >0} \lambda ^{-\alpha } \log \mathbb P ( \abs{ f}>\lambda )\,. 
\end{align*}

\end{proposition}

One also has a familiar Lemma for the maximum of  random variables. 

\begin{lemma}\label{l.subG} Let $ X_1 ,\dotsc, X_N$ be random variables 
in $ L ^{\psi }$ of norm at most one.  Then, we have 
\begin{equation*}
\mathbb E \sup _{n\le N} \abs{ X_N} \lesssim \psi ^{-1}  (N)\,. 
\end{equation*}
\end{lemma}

So for  $ X_1 ,\dotsc, X_N \in \operatorname {exp} (L ^{2})$ of norm one, we have 
\begin{equation}\label{e.sqrtLog}
\mathbb E \sup _{n\le N} \abs{ X_N} \lesssim \sqrt{\log N +1}\,. 
\end{equation}
Indeed, we will leave to the reader to  verify  that under the assumptions above 
\begin{equation}\label{e.expsqrtLog}
\norm \sup _{n\le N} \abs{ X_N} \, . \operatorname {exp} (L ^2 ).\lesssim \sqrt{\log N +1} \,. 
\end{equation}

\begin{proof}
By Jensen's inequality 
\begin{align*}
\psi (\mathbb E \sup _{n\le N} \abs{ X_N} ) 
& \le \mathbb E \sup _{n\le N} \psi (\abs{ X_N} ) 
\\
& \le \sum _{n=1} ^{N} \mathbb E \psi (\abs{ X_N} )
\\
& \lesssim N\,. 
\end{align*}
The proof is complete. 
\end{proof}

Another class of relevant spaces are given by the convex functions 
\begin{equation*}
\varphi _{\beta } (x) \coloneqq \abs{ x} (\log 2+\abs{ x})\,. 
\end{equation*}
We denote $ L ^{\varphi _{\beta }}= L (\log L) ^{\beta }$.  
The connection with the exponential Orlicz classes is by way of duality. 
\begin{equation}\label{e.duality}
[ \operatorname {exp} (L ^{\alpha })] ^{ \ast}
=L (\log L) ^{1/\alpha } \,. 
\end{equation}

These spaces are closely associated with the \emph{extrapolation} principle.  

\begin{proposition}\label{p.extra}  Let $ \operatorname T$ be a  linear operator with 
\begin{equation} \label{e.extra}
\norm \operatorname T . L ^{p} ([0,1]^{d})\rightarrow L ^{p} ([0,1]^{d}). 
\lesssim (p-1) ^{\alpha }\,, \qquad 1<p\le 2\,,\ 0<\alpha <1\,. 
\end{equation}
We then have the inequality 
\begin{equation} \label{e.extra+}
\norm \operatorname T f. L ^{1}. \lesssim \norm f. L (log L) ^{\alpha }. \,.
\end{equation}
More generally, 
\begin{equation}\label{e.more}
\norm \operatorname T f. L ^{1} (\log L) ^{\beta }.
\lesssim \norm f. L (log L) ^{\alpha +\beta }. \,, \qquad 0<\beta <\infty\, .
\end{equation}
\end{proposition}

\begin{proof}
Let us consider (\ref{e.extra+}).  This inequality is dual to 
\begin{equation*}
\norm \operatorname T ^{\ast} f. \operatorname {exp} (L) ^{1/\alpha }.
\lesssim \norm f.\infty . \,.
\end{equation*}
But, taking $ f\in L ^{\infty }$, with $ \norm f. \infty .=1$, and using 
(\ref{e.extra}), we have for $ 2<p<\infty $, 
\begin{align*}
\norm \operatorname T ^{\ast} f. p.& \lesssim 
 p ^{\alpha } 
\end{align*}
and so the dual estimate follows Proposition~\ref{p.exp}.  The inequality 
(\ref{e.more}) is entirely similar. 
\end{proof}

\section{Khintchine Inequalities} \label{s.khintchine}

The utility of the exponential Orlicz classes is that they allow a concise 
expression of a range of inequalities.  This is especially relevant to the 
classical Khintchine Inequalities. In other instances we shall see, 
that Orlicz spaces express sharp inequalities forms of different inequalities. 

Let $ \{r_k\mid k\ge 1\}$ be independent, identically distributed random variables, 
with $ \mathbb P (r_1=1)=\mathbb P (r_1=-1)=\tfrac12$.  Such random variables are 
referred to as Rademacher random variables.  They admit different realizations, 
of which the most direct is 
\begin{equation*}
r_k=\operatorname {sgn} (\sin (2 ^{k} \pi x))\,, \qquad 0\le x \le 1\,. 
\end{equation*}
Such random variables are in particular orthogonal, so that we have 
\begin{equation*}
\NOrm \sum _{k} a_k r_k .2.= \Bigl[\sum _k a_k ^2  \Bigr] ^{1/2} \,. 
\end{equation*}
This holds for all finite sequences of constants $ \{a_k\}$.

The Khintchine Inequality says that these sums, in all $ L ^{p}$, are 
controlled by the $ L ^{2}$ norms.  In its sharp form, this inequality states 

\begin{khintchine} For all finite sequences of constants $ \{a_k\}$
\begin{equation}\label{e.khintchine}
\NOrm \sum _{k} a_k r_k .\operatorname {exp} (L ^2 ). 
\lesssim 
\Bigl[\sum _k a_k ^2  \Bigr] ^{1/2} \,. 
\end{equation}
\end{khintchine}

\begin{proof}
The classical proof of this is quite elementary, passing through the Moment 
Generating Function.  We can restrict attention to the case where 
\begin{equation*}
\Bigl[\sum _k a_k ^2  \Bigr] ^{1/2} =1\,. 
\end{equation*}

Consider the moment generating function, given by 
\begin{align*}
\varphi (\lambda )&= \mathbb E \operatorname e ^{\lambda \sum _{k} a_k r_k }\,, \qquad
\lambda >0
\\
& = \prod _{k} \mathbb E \operatorname e ^{\lambda a_k r_k} 
\\
&= \prod _{k} \tfrac12 ( \operatorname e ^{-\lambda a_k} + \operatorname e ^{\lambda a_k})
\\
&\le\prod _{k} \operatorname e ^{\lambda ^2 a_k ^2 } 
\\
& \le \operatorname e ^{\lambda ^2 }
\end{align*}
Here, we have relied statistical independence of the random variables. 
In particular, if $ X, Y $ are independent random variables, then 
\begin{equation*}
\mathbb E X \cdot Y= \mathbb E X \cdot \mathbb E Y\,. 
\end{equation*}
We have also used the  the elementary inequality 
\begin{align} \label{e.SYM}
\tfrac12(
\operatorname e ^{-\mu }+\operatorname e ^{\mu })
=\sum _{j=1} ^{\infty } \frac{\mu ^{2j}} { (2j)!}
\le \operatorname e ^{\mu ^2 }
\,, \qquad \mu \in \mathbb R \,.
\end{align}

Now estimate 
\begin{align*}
\mathbb P \Bigl(  \sum _{k} a_k r_k >t \Bigr) 
\le \varphi (\lambda ) \operatorname e ^{-\lambda t}
\le \operatorname e ^{\lambda ^2 -\lambda t}\,, \qquad \lambda >0\,. 
\end{align*}
The minimum over $ \lambda >0$ of the right hand side occurs at $ \lambda =t/2$, 
giving us the estimate 
\begin{equation*}
\mathbb P \Bigl(  \sum _{k} a_k r_k >t \Bigr) 
\le \operatorname e ^{-t ^{2}/4}\,. 
\end{equation*}
In view  of  the symmetry of the Rademacher random variables and 
Proposition~\ref{p.exp}, this proves the Theorem. 

\end{proof}

\section{Maximal Function Estimates} 

While our primary interest is in the Littlewood Paley Theory, the maximal 
function and its relevant estimates are essential to the subject. 

Define 
\begin{equation}\label{e.DefMaximal} 
\operatorname M f (x) =\sup _{\substack{x\in I\\ I \in \mathcal D }} \mathbb E (f \,|\, I)\,. 
\end{equation}
The principal properties of the Maximal function are 

\begin{theorem}\label{t.maximal} We have the estimates 
\begin{equation}\label{e.maximal}
\begin{split}
\sup _{\lambda } \lambda \mathbb P ( \operatorname M f >\lambda ) \le \norm f.1.
\\
\norm \operatorname M f .p. \lesssim (1 +1/(p-1)) \norm f. p. \,, \qquad 1<p\le \infty \,. 
\end{split}
\end{equation}
\end{theorem}

The left hand side of the first inequality is referred to as the \emph{weak $ L ^{1}$ norm}, 
and we write it as $ \norm \operatorname M f.1,\infty .$. 
More generally, we define 
\begin{equation}\label{e.weakNorm}
\norm f.p, \infty . =\sup _{\lambda >0} \lambda ^{-1} \mathbb P (\abs{ f}>\lambda ) ^{1/p}\,.
\end{equation}
As with the Orlicz norms, in certain instances these norms define sharp inequalities.  

\begin{proof}[Proof of Theorem~\ref{t.maximal}.] 
This is especially easy as we are working with the dyadic maximal function, this 
is especially easy.  We begin with the weak type inequality. 

Fix $ \lambda >0$, and let $ \Lambda $ be the collection of maximal dyadic intervals 
with $ \mathbb E (f \, |\, I)\ge \lambda $.  By maximality these intervals 
are disjoint, so 
\begin{align*}
\lambda \mathbb P ( \operatorname M f >\lambda ) 
& = \sum _{I\in \Lambda } \mathbb E f \mathbf 1 _{I } 
\le \mathbb E f \le \norm f.1.\,. 
\end{align*}

For the proof of the remaining inequalities, one interpolates with the obvious 
$ L ^{\infty }$ bound, as is described in Stein and Weiss \cite{MR0304972}.

\end{proof}

The norm estimate we give above, as $ p \downarrow 1$ is sharp, 
which extrapolates to this estimate

\begin{theorem}\label{t.stein}
We have the estimate 
\begin{equation}\label{e.stein}
\norm \operatorname M f. 1. \lesssim \norm f. L \log L.\,. 
\end{equation}
\end{theorem}

\begin{proof}
This nearly follows from Proposition~\ref{p.extra}, but $ \operatorname M $ 
is not a linear operator.  Yet, the bound for the maximal operator 
in Theorem~\ref{t.maximal} is equivalent to the same bound for the family 
of linear  operators 
\begin{equation*}
\operatorname T f 
\coloneqq \sum _{I \in \mathcal D} \mathbf 1 _{E (I)} \mathbb E (f\,|\, I)\,, 
\end{equation*}
where $ \{ E (I)\mid I\in \mathcal D\}$ is a family of pairwise disjoint 
sets with $ E (I)\subset I$ for all $ I$.   (For a given $ f$, one takes $ E (I)$ 
to be the set of $ x\in I$ for which the supremum in the definition of $ \operatorname M$ 
is achieved at $ I$.) 

These operators, being linear, satisfy the estimate (\ref{e.stein}), by
Proposition~\ref{p.extra}.  Therefore the Lemma follows. 
\end{proof}

There is a striking converse to this last Theorem, 

\begin{theorem}\label{t.STEIN}[E.\thinspace M.\thinspace Stein]
IF $ \operatorname M \abs  f \in L $, then we have $  f\in L \log L$. 
\end{theorem}

\begin{proof}
We can assume that $ f\ge 0$.  Let us first show that 
\begin{equation}\label{e.s1} 
 \lambda ^{-1} \mathbb E f \mathbf 1 _{ \{ f\ge \lambda \} } \lesssim
 \mathbb P (\operatorname M f > \lambda )\,, \qquad \lambda > \mathbb E f\,. 
\end{equation}
Indeed, let $ \mathcal I$ be the collection of maximal dyadic intervals 
with $ \{\operatorname M f > \lambda \}=\bigcup _{I\in \mathcal I} I$. 
Then, if $ x\not\in \{\operatorname M f > \lambda \}$, we must have $ f (x)\le \lambda $
by the Martingale Convergence Theorem.  
In addition,  $ \lambda > \mathbb E f$, so no $ I\in \mathcal I$ can be maximal.  That implies 
that  $ \mathbb E (f\,|\, I)\le 2 \lambda $.
But then, 
\begin{align*}
\mathbb E f \mathbf 1 _{ \{ f\ge \lambda \} } 
& \le \sum _{I\in \mathcal I} \mathbb E f \mathbf 1 _{I} 
\\
& \le2 \lambda  \sum _{I\in \mathcal I} \abs{ I} 
\\
&= 2 \lambda \mathbb P (\operatorname M f > \lambda )\,. 
\end{align*}

Hence, we can estimate 
\begin{align*}
\int _{ \mathbb E f} ^{\infty } \lambda ^{-1} \mathbb E f \mathbf 1 _{ \{f>\lambda \}} \; d \lambda 
& \le 
\int _{ \mathbb E f} ^{\infty } \mathbb P (\operatorname M f > \lambda )\; d \lambda \,,
\end{align*}
and our conclusion follows easily from this. 
\end{proof}

\section{Littlewood Paley Theory} \label{s.LP}

We consider the Haar basis on $ [0,1]$, given by 
$ \{\mathbf 1 _{[0,1]}\}\cup \{h_I \mid I\in \mathcal D\}$, where we 
remind the reader that $ \mathcal D$ consists of the dyadic intervals in $ [0,1]$. 
We also remind the reader that the Haar functions are normalized to have $ L ^{\infty }$
norm one, so that our formulas are different from most of our references.  

It is important to our applications that we consider the Haar basis as 
one for vector valued functions. The vector space should be a Hilbert space 
$ \mathcal H$, and by $ L ^{p} _{\mathcal H}$ we mean the class of measurable 
functions $ f \mid [0,1] \longrightarrow \mathcal H$ such  that 
\begin{equation*}
\mathbb E \abs{ f} _{\mathcal H} ^{p} < \infty \,. 
\end{equation*}

The Haar Square Function is 
\begin{equation*}
\operatorname S (f) \coloneqq 
\Bigl[\abs{ \mathbb E f} ^2 + \sum _{I\in \mathcal D} \frac {\abs{ \ip f,h_I,} ^2 } 
{\abs{ I} ^2 } \mathbf 1 _{I}\Bigr] ^{1/2} \,. 
\end{equation*}
Here, we are taking the Hilbert space norm of those terms that involve $ f$. 
Of course we have $ \norm f.2.= \norm \operatorname S (f).2.$ just by the 
fact that the Haar basis is an orthogonal basis.  

The Littlewood Paley Inequalities are a profound extension of this equality, 
to an approximate version that holds on all $ L ^{p}$, $ 1<p<\infty $. 

\begin{lpi} \label{t.lpi}
For $ 1<p<\infty $ there are absolute constants $ 0<A_p<B_p<\infty $ so that 
\begin{equation}\label{e.LLPP1}
\begin{split}
 \norm f.p. &\le B_p \norm \operatorname S (f).p. \,, \qquad 1<p<\infty 
 \\
 B_p &\lesssim 1+\sqrt p\,.
\end{split}
\end{equation}
In the reverse direction, we have 
\begin{equation} \label{e.LLPP2}
\begin{split}
 A_p \norm \operatorname S (f).p. &\le \norm f .p. \,, \qquad 1<p<\infty 
 \\
 A_p &\simeq  1+1/\sqrt {p-1}\,.
\end{split}
\end{equation}
\end{lpi}

We stress that these results are delicate.  Burkholder \cite{MR976214} 
has shown that the best constants in the inequality above for \emph{general 
martingales} are $ A_p ^{-1}=B_p=\max\{p,q\}-1$.  However, a Haar series 
is not a general martingale; it is dyadic, which  forces conditional symmetry.  
See \cite{MR1018577}. 

The constants above are sharp.  To see that $ B_p \simeq \sqrt p$ is sharp for 
$ p$ large, just use the Central Limit Theorem for Rademacher random variables, 
or the sharpness of the Khintchine Inequality.  
A duality argument shows that one can take $ A_p=B _{p'} ^{-1} $, 
where $ p'$ is the conjugate index to $ p$. 

The inequality (\ref{e.LLPP1}) holds for $ 0<p<2$, but we do not need that 
case, so don't discuss it.

\subsection*{Duality Principle}
With the Littlewood Paley Inequalities, there is an important duality principle 
which permits us to pass from one inequality to another. 
Let us see that we can take 
\begin{equation}\label{e.dual1}
B_p=A_q ^{-1} \,, \qquad \tfrac1p+\tfrac1q=1\,. 
\end{equation}

Assume the inequality $A _{q}    \norm \operatorname S (g).q.\le \norm g.q. $. 
Fix $ f \in L  ^{p}$, and choose $ g\in L ^{q}$ of norm one so that we have 
\begin{align*}
\norm f.p. & = \ip f, g, 
\\
&= \mathbb E f \cdot _{\mathcal H} \mathbb E g +\sum _{I\in \mathcal D} \ip f, h_I, 
\cdot _{\mathcal H} \ip g, h_I, 
\\
& \le \norm \operatorname S (f).p. \cdot \norm \operatorname S (g).q. 
\\
& \le A_q ^{-1} \norm \operatorname S (f).p. 
\end{align*}
So (\ref{e.dual1}) holds.

\subsection*{The Chang Wilson Wolff Inequality}

A key step in the proof of this inequality is to first prove the 
Chang Wilson Wolff inequality, \cite{MR800004}.

\begin{cww} \label{t.cww} We have the estimate below for Hilbert space valued $ f$. 
\begin{equation}\label{e.cww}
\norm f. \operatorname {exp} (L ^2 ). \lesssim \norm \operatorname S (f).\infty . \,. 
\end{equation}
\end{cww}

\begin{proof}
It is immediately clear that if we knew $ \norm f.p. \lesssim \sqrt p \norm 
\operatorname S (f).p.$ for $ p\ge 2$, in the Hilbert space valued case, then 
the inequality (\ref{e.cww}) would follow.

Our strategy is to first prove the inequality (\ref{e.cww}) in the case that the function 
$ f$ is real valued.  From this, we will deduce a quadratic inequality, which will 
prove the Littlewood Paley inequalities for large $ p$, in the Hilbert space 
valued case.  This will complete the proof of the Chang Wilson Wolff inequality as 
we have stated it.

We give the proof of Chang Wilson and Wolff, in the real valued case, 
which they learned from Herman Rubin.  Indeed, this proof can be regarded as the 
conditional version of the proof we have already given of the Khintchine inequalities. 

Let us recall that a sequence of functions $ g_1 ,\dotsc, $ form 
a \emph{martingale} iff for all sequences 
\begin{equation*}
\mathbb E (g _{n+1} \,|\, g_1 ,\dotsc, g_n )= g_n\,. 
\end{equation*}
Here, we are taking the conditional expectation of $ g _{n+1}$ 
with respect to the sigma field generated by $ g_1 ,\dotsc, g_n$.

Let $ \mathcal F_n $ be the sigma field generated by the dyadic intervals of length $ 2
^{-n}$,  so that 
\begin{equation*}
f_n \coloneqq \mathbb E (f\,| \mathcal F_n\,)=\mathbb E f+\sum _{\abs{ I }\ge 2 ^{-n}} 
\frac{\ip f, h_I, } {\abs{ I}} h_I
\end{equation*}
is a dyadic martingale.  We assume that $ \mathbb E f=0$. 

For $ t>0$ we define a new martingale by the formula 
\begin{equation*}
q_n \coloneqq \operatorname e ^{t f_n}
\Bigl[\prod _{j=1} ^{n-1} \mathbb E ( \operatorname e ^{t (f_ {j+1}- f_j)} \,|\, \mathcal F_j) \Bigr]
^{-1} \,.
\end{equation*}

Of course, it is hardly obvious that $ q_n$ is a martingale, and so we check this now. 
Clearly, $ q_n$ is $ \mathcal F_n$ measurable.  We should then check that 
$ \mathbb E (q _{n+1}\,|\, \mathcal F_n)=q_n$. 
\begin{align*}
\mathbb E (q _{n+1}\,|\, \mathcal F_n)
& = 
\mathbb E \Bigl(\operatorname e ^{t f_ {n+1}} 
\Bigl[\prod _{j=1} ^{n} \mathbb E ( \operatorname e ^{t (f_ {j+1}- f_j)} 
\,|\, \mathcal F_j) \Bigr]^{-1} 
\,\big|\, \mathcal F_n\Bigr)
\\
& = \mathbb E (\operatorname e ^{t f_ {n+1}} 
\,|\, \mathcal F_n) \cdot 
\Bigl[\prod _{j=1} ^{n} \mathbb E ( \operatorname e ^{t (f_ {j+1}- f_j)} 
\,|\, \mathcal F_j) \Bigr]^{-1} 
\\
& = \mathbb E (\operatorname e ^{t (f_ {n+1}-f_n)} 
\,|\, \mathcal F_n) \cdot 
\operatorname e ^{t f _{n}} \cdot 
\Bigl[\prod _{j=1} ^{n} \mathbb E ( \operatorname e ^{t (f_ {j+1}- f_j)} 
\,|\, \mathcal F_j) \Bigr]^{-1} 
\\
&=\operatorname e ^{t f _{n}} \cdot 
\Bigl[\prod _{j=1} ^{n-1} \mathbb E ( \operatorname e ^{t (f_ {j+1}- f_j)} 
\,|\, \mathcal F_j) \Bigr]^{-1}
=q_n \,. 
\end{align*}
And  therefore, $ \mathbb E q_n=1$ for all $ n$. 

The fact that we work with a dyadic martingale enters. For we can appeal to (\ref{e.SYM}) 
to see that 
\begin{equation*}
\prod _{j=1} ^{n-1} \mathbb E ( \operatorname e ^{t (f_ {j+1}- f_j)} 
\,|\, \mathcal F_j)
\le \prod _{j=1} ^{n-1} \mathbb E ( \operatorname e ^{t ^2  (f_ {j+1}- f_j) ^2 } 
\,|\, \mathcal F_j)
=\prod _{j=1} ^{n-1}  \operatorname e ^{t ^2 (f_ {j+1}- f_j) ^2 } 
=\operatorname e ^{t ^2 \operatorname S (f) ^2 }\,. 
\end{equation*}

Therefore, under the assumption that $ \norm \operatorname S (f).\infty .\le 1$, we see 
that 
\begin{equation*}
\mathbb E \operatorname  e^{ t f_n-t ^{2}}\le \mathbb E q_n= 1\,.
\end{equation*}
As this holds for all $ n$, we can take $ n\to \infty $.  Therefore, we have for $ \lambda
>0$, 
\begin{equation*}
\mathbb P (f> \lambda  )\le \operatorname e ^{-t \lambda } \mathbb E e ^{t f}
\le \operatorname e ^{-t \lambda +t ^2 }\,. 
\end{equation*}
Taking $ t=\lambda /2$ proves the Chang Wilson Wolff inequality in the case that 
$ f$ is real valued.

\end{proof}

\subsection*{Proof of the Littlewood Paley Inequalities}

The first step is to derive a `Good $ \lambda $ Inequality,' as below.  
This exotic looking inequality, first devised in \cite{MR0440695}, has proven to
be a very powerful technique. 

\begin{good} For $ \lambda >0$ we have the inequality 
\begin{equation}\label{e.goodLambda}
\mathbb P ( \operatorname M f > 2 \lambda \,;\, \operatorname S (f) < \epsilon \lambda  )
\lesssim \operatorname e ^{-c \epsilon ^{-2}} \mathbb P (\operatorname M f >\lambda )
\,, \qquad 0<\epsilon < \tfrac 12 \,. 
\end{equation}
Here $ \operatorname M f$ is the dyadic maximal function, and $ 0<c<1$ 
is an absolute constant.  The point of the estimate is that it holds 
for all $ 0<\epsilon <\tfrac 12$, with the constant on the right tending to 
zero as $ \epsilon \downarrow 0$. 
\end{good}

\begin{proof}
Define a stopping time by 
\begin{equation*}
\tau =\min \Bigl\{n \mid \sum _{j=1} ^{n} (f_j-f _{j-1}) ^2 \ge \epsilon \lambda  \Bigr\}\,. 
\end{equation*}
As is usual, the minimum of the empty set will be taken to be $ +\infty $. 

Let $ f_I= \mathbb P (I) ^{-1} \mathbb E f \mathbf 1 _{I}$ be the average value of 
$ f$ on $ I$.

Let $ \mathcal Q$ be the maximal dyadic intervals with $  f_I
\ge \lambda \mathbb P (I)$, so that 
\begin{equation*}
\{\operatorname M f >\lambda \}=\bigcup _{I\in\mathcal Q}I\,. 
\end{equation*}
On each $ I$ the event $ E_I \coloneqq I\cap \{\operatorname M f > 2 \lambda \,;\, \operatorname S (f) <
\epsilon \lambda\}$. This is the main point: If $ E_I$ is non-empty then 
$ \mathbb E f \mathbf 1 _{I}
\le (1+\epsilon) \lambda \mathbb P (I)$.  Indeed, let $ I'$ denote the dyadic interval 
which contains it and is twice as long. So the average value of $ f$ on $ I'$ is 
less than $ \lambda $.  If our claim is not true, then 
\begin{equation*}
\abs{\ip f, h _{I'},} \ge \epsilon \lambda  \mathbb P (I)\,,
\end{equation*}
contradicting $ E_I$ being non empty.  

Now observe that 
\begin{align*}
\mathbb P (E_I)=\mathbb P ( \operatorname M f >2\lambda\, ;\, 
\tau =\infty )\le \mathbb P ( \operatorname M ( f _{\tau }- f_I)> (1-\epsilon )\lambda )\,.
\end{align*}
Moreover, $ \norm \operatorname S (f _{\tau }).\infty .\le \epsilon \lambda $. 
Therefore, by the Chang Wilson Wolff inequality applied to the renormalized 
martingale $ f _{\tau }- f_I$, 
\begin{equation*}
\mathbb P (E_I) \lesssim \operatorname e ^{ -c \epsilon ^{-2}} \mathbb P (I)\,. 
\end{equation*}
By summing over $ I\in \mathcal Q$ we complete the proof. 
\end{proof}

There is a standard way to pass from the Good $ \lambda $ Inequalities to norm inequalities, 
illustrated by this computation.  Since $ \abs{ f}\le \operatorname M f$, 
it suffices to prove the estimate  $\norm \operatorname M f. p. \lesssim B_p \norm 
\operatorname S (f).p.$.
First observe that 
\begin{align*}
\mathbb P (\operatorname  M f >2\lambda )
&\le \mathbb P ( \operatorname S (f)\le \epsilon \lambda )
+ \mathbb P ( \operatorname M f > 2 \lambda \,;\, \operatorname S (f) < \epsilon \lambda )
\\
& \le \mathbb P ( \operatorname S (f)\le \epsilon \lambda )
+ C \operatorname e ^{-c \epsilon ^{-2}}  
\mathbb P (\operatorname M f >\lambda )\,. 
\end{align*}
Then, we can compute 
\begin{align*}
\norm \operatorname M f. p. ^{p}
&= p 2 ^{p} \int _0 ^{\infty } \lambda ^{p-1} \mathbb P (\operatorname M f >2\lambda )\; d \lambda 
\\
& \le p 2 ^{p}\int_0 ^{\infty }  \mathbb P ( \operatorname S (f)\le \epsilon \lambda ) 
\; d \lambda 
+ 
p 2 ^{p} C \operatorname e ^{-c \epsilon ^{-2}}   
\int _0 ^{\infty } \lambda ^{p-1} \mathbb P (\operatorname M f >\lambda )\; d \lambda
\\
& \le (2/\epsilon ) ^{p} \norm \operatorname S (f).p. ^{p}
+ p 2 ^{p} C \operatorname e ^{-c \epsilon ^{-2}}    
\norm \operatorname M f. p. ^{p}\,. 
\end{align*}
Observe that if we take $ \epsilon \simeq p ^{-1/2}$, we can conclude 
\begin{equation*}
\norm \operatorname M f . p.  ^{p}\lesssim 
(C \sqrt p ) ^{p} \norm \operatorname S (f).p. ^{p}
\end{equation*}
which proves the desired inequality.  

To recap, we have proved (\ref{e.LLPP1}) in the range $ 1<p<\infty $ 
for real valued functions $ f$.  
By the duality principle, this proves (\ref{e.LLPP2}) in the same range.

To deduce the stronger result, for Hilbert space valued functions $ f$, 
we need a different formulation of the Chang Wilson Wolff inequality. 
 Fefferman and Pipher \cite{MR1439553} have devised an elegant proof,  inspired 
by the work of Wilson \cite{MR972707}. Also see \cite{MR850744}

\begin{definition}\label{d.Ap} For $ 1<p<\infty $, a function $ w\ge 0$ on $ [0,1]$, 
say that it is in dyadic $ A_p$ if 
\begin{equation}\label{e.Ap}
\norm w. A_p. \coloneqq \sup _{I\in \mathcal D}  \abs{ I} ^{-1} \norm w \mathbf 1 _{I}.p. 
\cdot \norm w ^{-1} \mathbf 1 _{I} . p/(p-1).  <\infty \,. 
\end{equation}
We are especially interested in the endpoint cases.  To be explicit, these are 
\begin{align*}
\norm w. A_1. &\coloneqq \sup 
_{I\in \mathcal D}  \abs{ I} ^{-1} \norm w \mathbf 1 _{I}.1. 
\cdot [ \inf _{x\in I} w (x)] ^{-1}   <\infty
\\
\norm w. A_\infty . &\coloneqq \sup 
_{I\in \mathcal D}  \sup _{x\in I} w (x)  
\cdot \abs{ I} ^{-1}    \norm w ^{-1}  \mathbf 1 _{I}.1.  <\infty
\end{align*}
\end{definition}

The functions $ w\ge 0$ are `weights' that we use to construct $ L ^{p}(w)$ 
spaces, with norm $ \norm f. L ^{p} (w). ^{p }= \mathbb E f ^{p} \cdot w$. 
By an abuse of notation, we will write this last expectation as 
\begin{equation*}
\mathbb E _{w} f \coloneqq \mathbb E f \cdot w\,. 
\end{equation*}
Likewise $ \mathbb P _{w} (A)=\mathbb E _{w} \mathbf 1 _{A}$. 	
The result we are interested in is:

\begin{theorem}\label{t.FP}  We have the inequality 
\begin{equation}\label{e.FP}
\norm f . L ^{2} (w). \lesssim \norm w. A ^{1}. ^{1/2} \norm \operatorname S (f). 
L ^{2} (w). 
\end{equation}
This holds for all Hilbert space valued $ f$. 
\end{theorem}

There are two key observations about this Theorem.  First, the estimate is quadratic 
in nature, a key reason for passing to this level of generality. 
In particular, in order to establish this Hilbert space valued $ f$, it suffices 
to establish it for real valued $ f$.  Indeed, if $ f$ takes values in a  
Hilbert space, then we can assume that the Hilbert space is $ \ell ^2 $, 
and write $ f= (f_k\mid k\in \mathbb N )$.  Assuming the real valued version, we can 
just sum on $ k$. 
\begin{align*}
\norm f . L ^{2} (w). ^2 &= \sum _{k\in \mathbb N } \norm f _k. L ^{2} (w). ^2 
	\lesssim \sum _{k\in \mathbb N } \norm w. A ^{1}. 
\norm \operatorname S (f _k). L ^{2} (w). ^2\,. 
\end{align*}
So the Hilbert space case is immediate. 

Second, the dependence in terms of the $ A ^{1}$ constant is sharp, which permits the 
deduction of the sharp growth rate in $ L ^{p}$ constants, for $ p>2$.  This is 
a standard argument, following Rubio de Francia.  For $ p>2$, write 
\begin{equation*}
\norm f.p. ^2 \le \mathbb E f^2 \cdot \varphi 
\end{equation*}
for some non-negative $ \varphi $ with $ \norm \varphi .(p/2)'.=1$. 
We dominate $ \varphi $ by an $ A ^{1}$ weight, which is given as follows. 
\begin{equation}\label{e.trick}
v \coloneqq  \sum _{k=0} ^{\infty } (2 \mu ( (p/2)')) ^{k}\operatorname M ^{k} \varphi \,. 
\end{equation}
In this display, $ \operatorname M$ denotes the dyadic maximal function, and 
$ \operatorname M ^{k}$ denotes the $ k$th power of $ \operatorname M$.  We 
interpret the $ 0$th power to be the identity.  The constant $ \mu (q)$
is the norm of the Maximal Function on $ L ^{q}$.  The relevant fact for us here is that 
$ \mu (q) \simeq (q-1) ^{-1} $ as $ q\downarrow 1$.  In particular, $ \mu ((p/2)') 
\simeq p$ as $ p\to \infty $.  It is clear that $ \norm v. (p/2)'. \lesssim 1$. 

Now $ v$ satisfies $ \norm v. A ^{1}. \lesssim p$, since for any dyadic interval $ I$
\begin{align*}
 \mathbb E (v \,|\, I)
& = \sum _{k=0} ^{\infty } 
(2 \mu ( (p/2)')) ^{k} \cdot \abs{ I} ^{-1} \mathbb E \mathbf 1 _{I}
\operatorname M ^{k} \varphi
\\
& \le 2 \mu ((p/2)') \sum _{k=1} ^{\infty } 
(2 \mu ( (p/2)')) ^{k} \inf _{x\in I} \operatorname M ^{k} \varphi
\\
& \lesssim p \inf _{x\in I} v (x)\,. 
\end{align*}
But then, we have 
\begin{align*}
\norm f.p. ^2 &\le \mathbb E _v f ^2  
\\&\lesssim \norm v. A ^{1}. \mathbb E _v \operatorname S (f) ^2 
\\&\lesssim p \mathbb E _v \operatorname S (f) ^2 
\\ & \lesssim p \norm S (f) .p. ^2 \,. 
\end{align*}
So the Littlewood Paley estimates holds for all $ p>2$, in the Hilbert space valued case.

\subsubsection*{Proof of Theorem~\ref{t.FP}.}
We need an additional result on the way in which $ A ^{1}$ weights embed in $ 
A ^{\infty }$ weights. 

\begin{lemma}\label{l.FP}[Lemma 3.6, \cite{MR1439553}.]
Given $ 0<\eta <1$, there is a $ C>0$ so that for all $ w\in A ^{1}$ 
and sets $ E\subset I$ where $ I$ is dyadic, we have 
\begin{equation*}
\mathbb P (E\,|\, I)< \operatorname e ^{-C \norm w. A ^{1}.} 
\quad \textup{implies} \quad 
 \mathbb P _{w} (E\,|\, I)< \eta 
\end{equation*}

\end{lemma}

\begin{proof}
As we work on a probability space, we have the H\"older inequality 
\begin{equation*}
\mathbb E \abs{ f}\le \norm  f.p.\,, \qquad 1<p<\infty \,,
\end{equation*}
as well as the Orlicz variants, $ \mathbb E \abs{ f}\le \norm f. L ^{1} \log L.$. 
It is a key attribute of the weighted theory that one can reverse some of these 
inequalities for weights $ w\in A_p$.  In the case of $ A_1$
the reverse H\"older inequality is
\begin{equation*}
\norm w . L ^{1} (\log L) (I; dx/\abs{ I}). \lesssim \norm w. A ^{1}. 
\norm w . L^{1} (I; dx/\abs{ I}). \,. 
\end{equation*}
This follows immediately from Theorem~\ref{t.STEIN} and the definition of 
$ A ^{1}$. 

But then we can estimate 
\begin{align*}
\mathbb E( w \mathbf 1 _{E} \,|\, I) &\le \norm w . L ^{1} (\log L) (I; dx/\abs{ I}). 
  \norm \mathbf 1 _{E} . \operatorname {exp} (L (I; dx/\abs{ I})). 
  \\ & \lesssim  \norm w . A ^{1}. \norm w. L ^{1} (I; dx/\abs{ I}). 
   \log \mathbb P (E\,|\,I) ^{-1} 
  \\
  \\& 
  \lesssim  \eta   
\norm w . L^{1} (I; dx/\abs{ I}). 
\end{align*}
This proves our Lemma. 

\end{proof}

Recall the Chang Wilson Wolff good $ \lambda $ inequality 
\begin{equation*}
\mathbb P ( \operatorname M f > 2 \lambda \,;\, \operatorname S (f) < \epsilon \lambda  )
\lesssim \operatorname e ^{-c \epsilon ^{-2}} \mathbb P (\operatorname M f >\lambda )
\,, \qquad 0<\epsilon < \tfrac 12 \,. 
\end{equation*}
Taking $ \epsilon \simeq \norm w. A ^{1}. ^{-1/2}$, we can deduce the 
weighted good $ \lambda $ inequality 
\begin{equation*}
\mathbb P_w ( \operatorname M f > 2 \lambda \,;\, \operatorname S (f) < \epsilon \lambda  )
\lesssim \eta  \mathbb P_ (\operatorname M f >\lambda )\,. 
\end{equation*}
And the standard way to prove the $ L ^{2}$ estimate from the good $\lambda  $ 
inequality gives us the inequality 
\begin{equation*}
\norm f. ^2 . L ^{2} (w). \lesssim c \epsilon ^{-2} \norm \operatorname S (f).2. 
\simeq \norm 2. A ^{1}. \norm \operatorname S (f).2.\,, 
\end{equation*}
and so the proof is done.

\subsection*{Weak $ L ^{1}$ Estimate}

At $ L ^{1}$, the equivalence $ \norm f.1. \simeq \norm \operatorname S (f).1.$ 
fails.\footnote{Instead, one has $ \norm \operatorname M f.1. \simeq \norm \operatorname S (f).1.$
where $ \operatorname M$ is Maximal function.  The theory of Hardy space $ H ^{1}$ depends critically 
on this equivalence.}  Nevertheless, there is an endpoint estimate of 
interest to us.  It is 

\begin{weak} We have the inequality 
\begin{equation}\label{e.weak}
\sup _{\lambda >0} \lambda\mathbb P (\operatorname S (f)>\lambda )
\lesssim \norm f.1. \,. 
\end{equation}
We stress that this inequality holds for Hilbert space valued functions $ f$. 
\end{weak}

\begin{remark}\label{r.later} Traditional approaches to these issues treat the 
weak $ L ^{1} $ estimate first, and then interpolate to $ L ^{p}$.  We are interested in the 
sharp constants for the square function, which are not available by way of the 
weak $ L ^{1}$ norm. 
\end{remark}

Central to the proof of this estimate is the \emph{Calder\'on Zygmund Decomposition.}

\begin{decomp}
For $ f\in L ^{1} _{\mathcal H}$  of norm one, and let $ \lambda >0$.  Then, we can write 
$ f=g_1+g_2$ so that  $ \norm g_1.\infty . \le \lambda $, and $ g_2$ 
is supported on disjoint dyadic intervals $ \{I_j\mid j\ge 1\}$, with 
\begin{equation}\label{e.g2}
\Abs{\bigcup I_j} \le \lambda ^{-1} \,, 
\qquad 
\mathbb E (g_2 \,|\, I_j)=0\,. 
\end{equation}
\end{decomp}

\begin{proof}
This is a stopping time argument, but as we work on the dyadic grid, 
the details simplify considerably.  Take $ \{I_j\}$ to be the maximal dyadic 
intervals such that 
\begin{equation*}
\mathbb E (\abs f \,|\, { I_j})\ge \lambda   \,. 
\end{equation*}
Maximality assures us that these intervals are disjoint. 
Since $ \norm f.1.=1$, we have 
\begin{equation*}
\sum _{j} \abs{ I_j}\le \lambda \sum _{j} \mathbb E (\abs f \mathbf 1 _{ I_j})\le 1\,. 
\end{equation*}

Set 
\begin{equation*}
g_1 (x)= 
\begin{cases}
\mathbb E (\abs f \,|\, { I_j})  & x \in I_j\,,\ j\ge 1\,, 
\\
f (x) & \textup{otherwise.}
\end{cases}
\end{equation*}
By the Lebesgue Differentiation Theorem (or Martingale Convergence Theorem), 
$ \norm g_1. \infty . \le \lambda $. 

It is then clear that we have 
\begin{equation*}
g_2 \mathbf 1 _{I_j} = f \mathbf 1 _{I_j} - \mathbb E (\abs f \,|\, { I_j})\,. 
\end{equation*}
Thus, $ g_2$ satisfies all its desired properties. 
\end{proof}

\begin{proof}[Proof of (\ref{e.weak}).] 
Fix $ f\in L ^{1} _{\mathcal H}$ of norm one  and $ \lambda >0$. 
As we work on a probability space, we can further restrict attention to $ \lambda >1$. 
Apply the Calder\'on Zygmund Decomposition,  writing $ f=g_1+g_2$. 

Note that we have 
\begin{equation*}
 \mathbb P ( \operatorname S (f) >2 \lambda ) 
\le  \mathbb P (\operatorname S (g_1)>\lambda )+ 
 \mathbb P (\operatorname S (g_2)>\lambda )\,,
\end{equation*}
so that is suffices to analyze the two terms on the right separately. 

For $ g_1$, we use the $ L ^{2}$ bound for the Square Function
so that 
\begin{align*}
\lambda ^{2}  \mathbb P (\operatorname S (g_1)>\lambda )
&\le \norm \operatorname S (g_1).2.^2 
\\
& \lesssim \norm g_1.2. ^2 
\\
&= 2\int _{0} ^{\lambda } u\mathbb P (\abs{g_1}>u) \; du 
\\
& \le 2\lambda \,. 
\end{align*}
The matches the required bound from (\ref{e.weak}).

The case of $ g_2$ is simpler.  The function $ g_2$ is supported on the dyadic 
intervals $ I_j$, and has mean zero on each dyadic interval.  Thus, if $ J$ 
is any dyadic interval that strictly contains an $ I_j$, we must have 
$ \ip g_2, h_J,=0$.  It follows that the square function of $ g_2$ is supported 
on the $ I_j$, so that 
\begin{equation*}
\mathbb P (\operatorname S (g_2)>0)
\le \sum _{j}\mathbb P (I_j) \le \lambda ^{-1} \,. 
\end{equation*}
Our proof is complete. 

\end{proof}

We need further extensions of the Chang Wilson Wolff inequality, namely these 
extensions, which are essentially known. 

\begin{theorem}\label{t.beyondCWW} For $ \beta \ge 0$ we have 
\begin{equation*}
\norm \operatorname {S} (f) . L ^{1} (\log L) ^{\beta }. 
\lesssim \norm f.  L ^{1} (\log L) ^{\beta +1/2}. 
\end{equation*}
Again, this holds for Hilbert space valued functions $ f$.
\end{theorem}

\begin{proof}
A variant of the duality principle is useful to us.  
We can choose a function $  g$ so that $ \operatorname S (g)$ 
has $ \operatorname {exp} (L ^{1/\beta })$ norm one for which 
\begin{align*}
\norm \operatorname {S} (f) . L ^{1} (\log L) ^{\beta }.
& =  \mathbb E f \cdot _{\mathcal H} \mathbb E g+ \sum _{I\in \mathcal D}
\abs{ I} ^{-1} \ip f, h_I, \cdot _{\mathcal H} \ip g, h_I, 
\\
& = \ip f, g, 
\\
& \le \norm f.  L ^{1} (\log L) ^{\beta +1/2}. 
\cdot \norm g. \operatorname {exp} (L ^{(\beta +1/2) ^{-1} }). 
\end{align*}
Now, by Proposition~\ref{p.exp}, and the sharp Littlewood Paley inequalities, 
\begin{align*}
\norm g. \operatorname {exp} (L ^{(\beta +1/2) ^{-1} }). 
&\simeq  
\sup _{r>2} r^{-(\beta +1/2)} \norm g.r. 
\\
& \lesssim \sup _{r>2} r ^{-\beta } \norm \operatorname S (g). r. 
\\
& \lesssim 1\,. 
\end{align*}
Our proof is complete. 
\end{proof}

\section{Product Theory} \label{s.productTheory}

The product theory is a branch of Harmonic Analysis devoted to a range of issues 
that are effectively analyzed with tensor products of Haar bases.\footnote{A more 
typical description involves questions that are invariant under a  family of dilations 
of two or more dilations. Dilations don't appear in these notes due to the local 
nature of the questions studied.} 

To describe this, again due to the local nature of the questions, we need to 
slightly modify the dyadic intervals. Before, we used $ \mathcal D$ to denote 
the dyadic intervals contained in $ [0,1]$.  Let us set $ \mathcal D_+$ 
to be these dyadic intervals together with the interval $[0,2] $.  Let us 
define the Haar function associated with $ [0,2]$ to be the constant function. 
\begin{equation*}
h _{[0,2]} = \mathbf 1 _{[0,1]}\,. 
\end{equation*}
(We could have taken these steps earlier, but it would have been confusing to do so.) 
Then, $ \{h_I\mid I\in \mathcal D_+\}$ is an orthogonal basis for $ L ^{2} ([0,1])$.

Let us construct the tensor product basis for $ L ^{2} ([0,1] ^{d})$.  The basis 
elements are indexed by $ \mathcal R _{d} \coloneqq \mathcal D _+ ^{d}$, and 
for $ R=R_1 \times \cdots \times R_d\in \mathcal R _{d}$, set 
\begin{equation*}
h _{R} (x_1, ,\dotsc, x_d) = \prod _{s=1} ^{d} h _{R_j} (x_j)\,. 
\end{equation*}
This is an orthogonal basis for $ L ^{2} ([0,1] ^{d})$.  As in the one parameter setting, 
we are interested in the vector valued version of this space.  

The Haar Square Function in this setting is 
\begin{equation}\label{e.Sd}
\operatorname S (f) \coloneqq \Bigl[ \sum _{R\in \mathcal R_d} 
\frac {\abs{\ip f, h_R,} ^2 } { \abs{ R} ^2 } \mathbf 1 _{R}\Bigr] ^{1/2} \,. 
\end{equation}
As in the one parameter setting, it is clear that $ \norm f.2.=\norm \operatorname S (f).2.$,
and there is a deep extension of this equivalence to all $ L ^{p}$.  Again, we are interested
in the version of this result which has the sharp dependence in $ p$. 

\begin{theorem}\label{t.LPd} We have the inequalities below, valid on $ [0,1] ^{d}$. 
\begin{align}\label{e.LLPPd1}
\norm f. p. \lesssim p ^{d/2} \norm \operatorname S (f).p.\,, \qquad  1<p<\infty \,
\\ \label{e.LLPPd2}
\norm \operatorname S(f). p. \lesssim (p-1) ^{-d/2} \norm \operatorname f.p.\,, 
\qquad 1<p< \infty \,.
\end{align}

\end{theorem}

\begin{proof}
The Duality Principle is still in effect, and so it suffices to prove 
one of the set of inequalities above.  We prefer to prove the first inequalities.

The method of proof is a standard iteration of the one parameter inequalities, 
in the vector valued setting, a common technique in the subject, see 
for instance \cites{MR0252961,MR0290095}. 

Observe that the product Square Function is the composition of 
Square Functions applied in each coordinate.  These Square Functions 
are then applied to Hilbert space valued functions.  
In particular, let $ \operatorname S_j$ be the one parameter square function 
applied in the coordinate $ x_j$.  Then, 
\begin{equation*}
\operatorname S = \operatorname S_1 \circ \cdots \circ\operatorname S_d\,. 
\end{equation*}
Note that in applying 
$ \operatorname S_1 ,\dotsc, \operatorname S _{d-1}$, one should interpret it 
as  applied to a
Hilbert space valued functions.  Namely in two dimensions, we interpret 
\begin{equation*}
\operatorname S_1 f (x_1,x_2) 
\coloneqq 
\bigl\{ \frac { \ip {f (x_1,x_2)} , h_{I_1} (x_1), } {\sqrt{\abs{ I_1}}} \mathbf 1 _{I_1} (x_1)
\mid I_1\in \mathcal D_+  \bigr\}\,,
\end{equation*}
and one computes the $ \ell ^{2} (\mathcal D_+)$ norm of this quantity.  
Then, 
\begin{equation*}
\operatorname S_2 \circ \operatorname S_1 f (x_1,x_2) 
\coloneqq 
\bigl\{ \frac { \ip {f(x_1,x_2)} , h_{I_1} (x_1) h_{I_2} (x_2), } {\sqrt{\abs{ I_1} 
\cdot \abs{ I_2}} }  \mathbf 1 _{I_1} (x_1) \mathbf 1 _{I_2} (x_2)
\mid I_1,I_2\in \mathcal D_+  \bigr\}\,,
\end{equation*}
And one computes the $ \ell ^2 {\mathcal D_+ \times \mathcal D_+}$ norm of the 
right hand side.  

It is clear that the Theorem then follows from the one parameter Littlewood Paley 
inequalities.

\end{proof}

\begin{remark}\label{r.alternate}
Alternately, one can use the weighted inequality in Theorem~\ref{t.FP}, applied 
$ d$ times.  Details are left to the reader. 
\end{remark}

We briefly mention some other relevant inequalities.  The weak type estimate is replaced by 

\begin{theorem}\label{t.weakd} We have the inequality below on $ [0,1] ^{d}$. 
\begin{equation} \label{e.weakd}
\norm \operatorname S (f). 1, \infty . \lesssim \norm f . L ^{1} (\log L) ^{d-1} \,. 
\end{equation}
\end{theorem}

The Maximal Function is 
\begin{equation*}
\operatorname M f (x)=\sup _{R\in \mathcal D} \mathbf 1 _{R} (x) \mathbb E (f\,|\, R)\,. 
\end{equation*}
The principal inequalities are below, and in general are sharp.  

\begin{theorem}\label{t.maxd} We have the inequalities 
\begin{align*}
\norm M f. L ^{1,\infty }.& \lesssim \norm f. L ^{1} (\log L) ^{d-1}. \,, 
\\
\norm M f. p. \lesssim [ 1+ (p-1) ^{-1} ] ^{-d} \norm f.p.\,. 
\end{align*}
\end{theorem}

As we don't use this estimate, we do not prove it.

\chapter{Other Applications: Approximation Theory and Probability Theory}

 \section{Mixed Derivatives}  \label{s.realapplications} 
 
 We will take an abbreviated  view of the subject of this chapter, referring the reader 
 to references, especially \cite{MR1005898} for more information. 
 In $d $ dimensions, 
 consider the map 
 \begin{equation*}
 \operatorname {Int}_d f(x_1,\cdots,x_d) {}\coloneqq{}\int_0 ^{x_1 }\!\!\cdots \!\!\int_0 ^{x_d } f(y_1,\cdots,y_d)\; dy\cdots dy_d 
 \end{equation*}
 We consider this as a map from $L ^{p }([0,1]^d) $ into $C([0,1]^d) $. 
 Clearly, the image of $\operatorname {Int}_d $ consists of 
 functions with $ L ^{p}$ integrable mixed partial derivatives.  
 Let us set 
 \begin{equation*}
 \operatorname {Ball}(MW ^{p} ([0,1] ^{d}))
 \coloneqq   \operatorname {Int}_d ( \{f \in L ^{p} ([0,1] ^{d}) 
 \mid \norm f.p. \lesssim 1\})\,.
\end{equation*}
That is, this is the image of the unit ball of $ L ^{p}$. This is the unit ball of the  space of 
functions with mixed derivative in $ L ^{p}$.  Our main theorem Theorem~\ref{t.bl} has 
consequences for the case of $ p=1$, but in this discussion we concentrate of the 
case of $ p=2$, for which we have no new results.

These sets are  compact in  
in $C([0,1]^d) $, and it is of relevance to quantify the compactness.  The traditional 
way to do this is through \emph{entropy numbers}. 
For $0<\epsilon<1 $, set  $N(\epsilon) $ to be the least number 
$N $ of  points $x_1,\cdots,x_N\in C([0,1]^d) $ so that 
\begin{equation*}
 \operatorname {Ball}(MW ^{2} ([0,1] ^{d}))
\subset\bigcup _{n=1}^N x_n+\epsilon B_\infty . 
\end{equation*}
Here, $ B _{\infty }$ is the unit ball of $ C ([0,1] ^{d})$. 
An  upper bound on these numbers is known, 
\begin{equation} \label{e.Ne<}
\log N(\epsilon) {}\lesssim{} \epsilon ^{-1} 
 (-\log \epsilon)^{d-1/2} 
\end{equation}
 And the task at hand is to prove that this estimate is sharp.  
 The case of $ d=2$ below follows from Talagrand \cite{MR95k:60049}.

\begin{conjecture} For $ d\ge 2 $ one has the estimate 
\begin{equation*}
\log N(\epsilon) \simeq  \epsilon ^{-1} 
 (\log 1/\epsilon)^{d-1/2} \,,  \qquad \epsilon \downarrow 0\,. 
\end{equation*}
\end{conjecture}

 How does Small Ball Conjecture enter in? We should use a 
 a `smooth' version of the Small Ball Conjecture. That is,
 in the Small Ball Conjecture, (\ref{small}), 
 one should replace the 
 `rough' Haar functions by smooth variants.  There is no canonical 
 way to do this,\footnote{One can replace the splines below by tensor products of wavelets, or by appropriate 
 hyperbolic trigonometric polynomials.} 
 and so we simply choose the `spline variant'
 of Talagrand \cite{MR95k:60049}.
 For dyadic interval $ R\in \mathcal D$ in dimension $ d$, set $ u _{R}
 = \operatorname {Int}_d h_R$.   
 
 \begin{smooth}  \label{smooth} For all sequences $ \alpha (R)$, we have the estimate 
 below valid for all integers $ n$. 
 \begin{equation}\label{e.smooth}
2 ^{-2n} \sum _{\abs{ R}= 2 ^{-n}} \abs{ \alpha (R)} \lesssim n ^{(d-2)/2}
\NOrm \sum _{\abs{ R}= 2 ^{-n}} \alpha (R) u _{R} . \infty . 
\end{equation}
 \end{smooth}
 
 The power  $ 2 ^{-2n} $ is explained by the fact that the functions $ 
  u_R$ have $ L ^{\infty }$ norm comparable to $ 2 ^{-n}$.  This inequality is 
  true, and proved by Talagrand \cite{MR95k:60049}, but the methods of 
  \S~\ref{s.talagrand} will provide simple proofs of related facts.

 Let us explain how this conjecture provides lower bounds for entropy numbers. 
 Given a choice of signs $\sigma\mid \{R \mid \abs{ R}=2 ^{-n}\}\longrightarrow\{\pm1\} $, we 
 consider the functions 
 \begin{equation*}
 F_\sigma {}\coloneqq{}n ^{(1-d)/2}  \sum _{\abs{ R} = 2 ^{-n}}\sigma(R) \, h_R .  
 \end{equation*}
 Then, the mixed derivative of 
 $  F_\sigma $ has norm about $ 1$.  The point of view is to let $\sigma $ vary to construct  sets of points in $\operatorname 
 {Int}_d (B_2 )$ that are widely separated.

 Suppose that  for two different choices of $\sigma $ and $\sigma' $, 
 we have 
 \begin{equation} \label{e.ze-big} 
 \sum _{\abs{ R} = 2 ^{-n}} \abs{\sigma(R)- \sigma '(R) } {}\gtrsim{}n ^{d-1}2^n . 
 \end{equation}
 Then, Conjecture~\ref{smooth} enters in the following way: 
 \begin{equation} \label{e.zI0} 
 \begin{split}
 \norm \operatorname {Int}_d( F_\sigma-F_{\sigma'}).\infty. 
 &= \NOrm  \sum _{\abs{ R}= 2 ^{-n}} (\sigma (R)- \sigma ' (R))  u _{R} . \infty .
 \\
 &\gtrsim{}n ^{-d+3/2} 2 ^{-2n} 
 	\sum _{\abs{ R} = 2 ^{-n}} \abs{\sigma(R)-\sigma(R') }
\\& 
	{}\gtrsim{}n ^{1/2} 2 ^{-n}  
	\end{split} 
	\end{equation}
Thus,  a collection of $ F_\sigma $ satisfying (\ref{e.ze-big}) 
are uniformly separated in $ L ^{\infty }$ norm. 
	
	Notice that we have reduced the problem to one of finding many  proportional subsets of $\mathcal R_n $ that 
	are essentially disjoint from each other.   This is addressed in a general fashion by this proposition.   
 \begin{proposition}\label{p.big-diff}
   There is a  constant $c>0 $ so that for all integers $m $,  there is a collection of subsets $\mathcal A $ 
of $\{1,\ldots,m\} $ so that 
\begin{align}  
		\label{e.big-diff-1} 
\operatorname {card}(A \triangle A')&{}\ge{} cm,\qquad A\not=A'\in\mathcal A, 
\\  \label{e.big-diff-2} 
\operatorname {card}(\mathcal A)&{}\ge{}\operatorname {exp} (c m) . 
\end{align}
\end{proposition} 

Apply this proposition  the collection of dyadic rectangles $ \{ R \mid \abs{ R}= 2 ^{-n}
\}$.  Let $ \mathcal A$ be the corresponding subsets of this collection, thus for 
$ A,A'\in \mathcal A$ we have $ \abs{ A \triangle A'} \gtrsim n ^{d-1} 2 ^{n}$. 
Let $ A$ also stand for the function 
\begin{equation*}
A (R) \coloneqq 
\begin{cases}
1 & R\in A
\\
-1 & R \not\in A 
\end{cases} \,. 
\end{equation*}
Consider the collection $ \{F_A \mid A\in \mathcal A\}$. 
Any two distinct functions in this collection obey the estimate (\ref{e.zI0}), hence 
it follows that 
\begin{equation*}
\log N ( n ^{1/2}2 ^{-n} ) \gtrsim \log (\sharp \mathcal A )\gtrsim  n ^{ d-1} 2
^{n}\,.
\end{equation*}
Setting $ \epsilon  \simeq n ^{1/2} 2 ^{-n}$, we see that we have 
\begin{equation}\label{e.Ne}
\log N (\epsilon ) \gtrsim \delta ^{-1} (\log 1/\epsilon ) ^{d-1/2}\,, \qquad \epsilon
\downarrow 0\,. 
\end{equation}
This would match the known upper bound, (\ref{e.Ne<}).  Again, this inequality 
is known, and a consequence of Talagrand's work, in dimension $ d=2$.

\subsection*{A Coding Theory Result}

A useful observation is that Proposition~\ref{p.big-diff} is concerned with the central issues of coding theory.  Namely, each 
subset of $\{1,\ldots,m\} $ is identified with a word of length $m $, in an alphabet of two colors.  The condition 
(\ref{e.big-diff-1}) implies that the words differ in a constant times $m $ slots---that is that their Hamming distance 
is proportionally as large as possible.   And the condition (\ref{e.big-diff-1}) assures us that the code has a large capacity. 
Fortunately, we can appeal to a well known result from Coding Theory to address this proposition.

\begin{theorem}\label{t.vg}
   [Varshanmov--Gilbert Bound].   We view $\{0,1\}^n $ as a linear vector space mod $2 $.  It contains   $V $, 
a linear subspace mod $2 $ with  
\begin{equation*}
\abs{x-y} _{\ell^1}\ge{} d \qquad x,y\in V 
\end{equation*}
iff the inequality below holds.
\begin{equation} \label{e.vg-bound} 
\binom {n-1}{n-1}+\cdots+\binom{n-1}{d-2}<2 ^{n-k} 
\end{equation}
\end{theorem}

To prove Proposition~\ref{p.big-diff}, in this Theorem, we take $d=k=\alpha n $, for a small constant $\alpha $ to be chosen.  
Note that the left hand side of (\ref{e.vg-bound}) is at most 
\begin{align*}
d \binom{n-1}{d-2}&{} \le{}\alpha n\binom n {\alpha n} 
\\ &
{} \sim{} \alpha n \frac{ n^n e ^{-n} }{\alpha^n (1-\alpha)^n  n ^n e ^ {-n } } \sqrt{\alpha(1-\alpha)} 
\\ & 
{}={} \alpha \sqrt{\alpha(1-\alpha)} n (\alpha^\alpha (1-\alpha)^{1-\alpha} ) ^{-n} 
\\ & 
<2 ^{(1-\alpha)n}. 
\end{align*}
Here, we are using  Stirling's formula $m!\sim{} m^{m-1/2} e^{-m} $, meaning that the ratio of these two 
terms approaches a non zero constant.  Observe that $\alpha^\alpha\rightarrow0 $ as $\alpha\rightarrow0 $, 
so that we can make a choice of $\alpha $ for which this inequality will be true for all large $n $.

\section{The Brownian Sheet}

\subsection*{General Gaussian Processes}

A \emph{Gaussian process} is a random map $ X_t \mid T \longrightarrow \mathbb R $ 
where $ T$ is some index set, so that for all finite $ S\subset T$ and reals $a_s $, 
\begin{equation*}
\sum _{s\in S} a_s X_s  
\end{equation*}
is a random variable with a Gaussian distribution.  It is a fundamental property that 
a mean zero Gaussian process is characterized by the covariances 
\begin{equation*}
\rho (s,t) \coloneqq \mathbb E X_s \cdot X_t\,. 
\end{equation*}
Throughout, we will be concerned with processes which are almost surely 
have bounded sample paths, namely 
\begin{equation*}
\mathbb P (\sup _{t\in T} \abs{ X_t}<\infty )=1\,. 
\end{equation*}

The \emph{Small Ball Problem} concerns estimates for the probability 
\begin{equation*}
\mathbb P (\sup _{t\in T} \abs{ X_t} < \epsilon )\,, \qquad \downarrow 0\,. 
\end{equation*}
See \cite{MR94j:60078} for a survey on these types of questions.

\smallskip

If one is given a subset $ K$ of a Hilbert space $ \mathcal H$, then one can 
define an associated mean zero Gaussian process $ X_s$ for $ s\in K$ by defining 
\begin{equation*}
\mathbb E X_s \cdot X_t \coloneqq \ip s,t, _{\mathcal H}
\end{equation*}
where the last inner product is the one associated with $ \mathcal H$.  
This is a canonical relationship with profound consequences:  Most Gaussian  processes 
of interest can be described in this manner, and the Hilbert space has 
function theoretic description which in turn reflects the structure of the 
Gaussian process. 

For instance, assume that associated with $ \{X_t\mid t\in T\}$ are covariance kernel functions $ K_t$ 
and measure $ \mu $ on $ T$ so that $ \{K_t \mid t\in T\}\subset L ^2 (T, d\mu )$ 
and 
\begin{equation*}
\mathbb E X_s \cdot X_t= \int _{T} K_s \cdot K_t \; d \mu \,.
\end{equation*}
Let $ \mathcal H _{X}$ be  the $ L ^2 (\mu )$ completion of the set of functions 
$ \{K_t \mid t\in T\}$.  This spaces 
is called the \emph{Reproducing Kernel Hilbert Space} 
associated with the Gaussian process $X_t $.

Following on the work Talagrand, Kuelbs and Li \cite{MR94j:60078} uncovered a close connection 
between the 
the Small Ball Probabilities and the covering numbers associated with the unit 
ball of $ \mathcal H_X$ in the $ L ^{\infty } (T)$ metric.  We will recall this 
result in the particular instances of the Brownian sheet below.

\subsection{The Brownian Sheet}

The Brownian sheet is a canonical Gaussian process indexed by points $s \in[0,1]^d $.  Calling the process $B(s) $, it is 
characterized by requiring it to be a mean zero process with covariance structure 
\begin{equation*}
\mathbb E B(s)B(t)=\prod _{j=1}^d \min(s_j,t_j)
\end{equation*}
Note that this covariance functional is given by 
\begin{equation*}
\mathbb E  B(s)B(t)= \int _{[0,1] ^{d}} \mathbf 1 _{[0,s)} \cdot \mathbf 1 _{[0,t)} 
\; dx\,. 
\end{equation*}

The Reproducing Kernel Hilbert Space associated with the Brownian sheet is 
$ WM ^{2} _{d}$, the Sobolev space of functions with square integrable 
mixed derivatives in dimension $ d$. A particular case of the result  of Kuelbs and Li \cite{MR2001c:60059}
states that

\begin{theorem}\label{t.kli} As $ \epsilon \downarrow 0$ we have 
\begin{equation}  \label{e.li-linde}
\log \mathbb P (\norm B .C([0,1]^d).<\epsilon) 
\simeq
\epsilon ^{-2 }(\log 1/\epsilon) ^{\beta} 
\quad \text{iff} \quad \log N(\epsilon)
{}\simeq{}\epsilon ^{-1} (\log 1/\epsilon) ^{\beta/2} . 
\end{equation}
\end{theorem}

Thus, the  Conjecture~(\ref{smooth}) gives a result on these processes. 
And the form of the relevant conjecture here is as follows. 

\begin{brownian} For dimension $ d\ge 2$ , we have 
\begin{equation*}
\log \mathbb P (\norm B .C([0,1]^d).<\epsilon) 
\simeq
\epsilon ^{-2 }(\log 1/\epsilon) ^{2d-1}  \,, \qquad \epsilon \downarrow 0\,. 
\end{equation*}
\end{brownian}

This is known for $ d=2$.  For all $ d\ge 3$, the upper bound on the Small Ball 
probabilities is known; the issue is to obtain the appropriate lower bound.  
In dimension $ d\ge 3$, the best known lower bounds miss the conjecture above by 
a single power  of $  \log 1/\epsilon  $.

\end{document}